\documentclass{amsart}
\usepackage{amsfonts}
\usepackage{amssymb, amsmath, amsthm}
\usepackage{amscd}
\usepackage{verbatim}
\usepackage{enumerate}

\allowdisplaybreaks
\newtheorem{theorem}{Theorem}[section]
\newtheorem{acknowledgement}[theorem]{Acknowledgement}

\newtheorem{definition}[theorem]{Definition}

\newtheorem{lemma}[theorem]{Lemma}

\newtheorem{proposition}[theorem]{Proposition}

\theoremstyle{remark}
\newtheorem{remark}[theorem]{Remark}

\numberwithin{equation}{section}

\def\typeout#1{\message{^^J}\message{#1}\message{^^J}}
\newif\ifSRCOK \SRCOKtrue
\newcount\PAGETOP
\newcount\LASTLINE
\global\PAGETOP=1
\global\LASTLINE=-1
\def\EJECT{\SRC\eject}
\def\WinEdt#1{\typeout{:#1}}
\gdef\MainFile{\jobname.tex}
\gdef\CurrentInput{\MainFile}
\newcount\INPSP
\global\INPSP=0
\def\SRC{\ifSRCOK  \ifnum\inputlineno>\LASTLINE    \ifnum\LASTLINE<0      \global\PAGETOP=\inputlineno    \fi    \global\LASTLINE=\inputlineno    \ifnum\INPSP=0      \ifnum\inputlineno>\PAGETOP        
      \fi    \else      
    \fi  \fi\fi}
\def\PUSH#1{\SRC\ifnum\INPSP=0 \global\let\INPSTACKA=\CurrentInput \else\ifnum\INPSP=1 \global\let\INPSTACKB=\CurrentInput \else\ifnum\INPSP=2 \global\let\INPSTACKC=\CurrentInput \else\ifnum\INPSP=3 \global\let\INPSTACKD=\CurrentInput \else\ifnum\INPSP=4 \global\let\INPSTACKE=\CurrentInput \else\ifnum\INPSP=5 \global\let\INPSTACKF=\CurrentInput \else               \global\let\INPSTACKX=\CurrentInput \fi\fi\fi\fi\fi\fi\gdef\CurrentInput{#1}\WinEdt{<+ \CurrentInput}\global\LASTLINE=0\ifSRCOK\fi\global\advance\INPSP by 1}
\def\POP{\ifnum\INPSP>0 \global\advance\INPSP by -1  \fi\ifnum\INPSP=0 \global\let\CurrentInput=\INPSTACKA \else\ifnum\INPSP=1 \global\let\CurrentInput=\INPSTACKB \else\ifnum\INPSP=2 \global\let\CurrentInput=\INPSTACKC \else\ifnum\INPSP=3 \global\let\CurrentInput=\INPSTACKD \else\ifnum\INPSP=4 \global\let\CurrentInput=\INPSTACKE \else\ifnum\INPSP=5 \global\let\CurrentInput=\INPSTACKF \else               \global\let\CurrentInput=\INPSTACKX \fi\fi\fi\fi\fi\fi\WinEdt{<-}\global\LASTLINE=\inputlineno\global\advance\LASTLINE by -1\SRC}
\def\INPUT#1{\relax}

\let\originalxxxeverypar\everypar
\newtoks\everypar
\originalxxxeverypar{\the\everypar\expandafter\SRC}
\everymath\expandafter{\the\everymath\expandafter\SRC}
\output\expandafter{\expandafter\SRCOKfalse\the\output}
\newif\ifSRCOK \SRCOKtrue
\newcount\PAGETOP
\newcount\LASTLINE
\global\PAGETOP=1
\global\LASTLINE=-1
\gdef\MainFile{\jobname.tex}
\gdef\CurrentInput{\MainFile}
\newcount\INPSP
\global\INPSP=0
\def\EJECT{\SRC\eject}
\def\WinEdt#1{\typeout{:#1}}
\def\SRC{\ifSRCOK  \ifnum\inputlineno>\LASTLINE    \ifnum\LASTLINE<0      \global\PAGETOP=\inputlineno    \fi    \global\LASTLINE=\inputlineno    \ifnum\INPSP=0      \ifnum\inputlineno>\PAGETOP              \fi    \else          \fi  \fi\fi}
\def\PUSH#1{\SRC\ifnum\INPSP=0 \global\let\INPSTACKA=\CurrentInput \else\ifnum\INPSP=1 \global\let\INPSTACKB=\CurrentInput \else\ifnum\INPSP=2 \global\let\INPSTACKC=\CurrentInput \else\ifnum\INPSP=3 \global\let\INPSTACKD=\CurrentInput \else\ifnum\INPSP=4 \global\let\INPSTACKE=\CurrentInput \else\ifnum\INPSP=5 \global\let\INPSTACKF=\CurrentInput \else               \global\let\INPSTACKX=\CurrentInput \fi\fi\fi\fi\fi\fi\gdef\CurrentInput{#1}\WinEdt{<+ \CurrentInput}\global\LASTLINE=0\ifSRCOK\fi\global\advance\INPSP by 1}
\def\POP{\ifnum\INPSP>0 \global\advance\INPSP by -1  \fi\ifnum\INPSP=0 \global\let\CurrentInput=\INPSTACKA \else\ifnum\INPSP=1 \global\let\CurrentInput=\INPSTACKB \else\ifnum\INPSP=2 \global\let\CurrentInput=\INPSTACKC \else\ifnum\INPSP=3 \global\let\CurrentInput=\INPSTACKD \else\ifnum\INPSP=4 \global\let\CurrentInput=\INPSTACKE \else\ifnum\INPSP=5 \global\let\CurrentInput=\INPSTACKF \else               \global\let\CurrentInput=\INPSTACKX \fi\fi\fi\fi\fi\fi\WinEdt{<-}\global\LASTLINE=\inputlineno\global\advance\LASTLINE by -1\SRC}
\def\INPUT#1{\relax}
\let\OldINCLUDE=\include
\def\include#1{\EJECT\PUSH{#1.tex}\OldINCLUDE{#1}\POP}

\let\originalxxxeverypar\everypar
\newtoks\everypar
\originalxxxeverypar{\the\everypar\expandafter\SRC}
\everymath\expandafter{\the\everymath\expandafter\SRC}
\let\zzzxxxbibliography=\bibliography
\def\bibliography#1{\PUSH{\jobname.bbl}\zzzxxxbibliography{#1}\POP}
\output\expandafter{\expandafter\SRCOKfalse\the\output}
\input{tcilatex}

\begin{document}
\title[Model for biological aggregations]{Global attractor for a nonlocal
model for biological aggregation}
\author{Ciprian G. Gal}
\address{Department of Mathematics, Florida International University, Miami,
FL 33199, USA}
\email{cgal@fiu.edu}
\maketitle

\begin{abstract}
We investigate the long term behavior in terms of global attractors, as time
goes to infinity, of solutions to a continuum model for biological
aggregations in which individuals experience long-range social attraction
and short range dispersal. We consider the aggregation equation with both
degenerate and non-degenerate diffusion in a bounded domain subject to
various boundary conditions. In the degenerate case, we prove the existence
of the global attractor and derive some optimal regularity results.
Furthermore, in the non-degenerate case we give a complete structural
characterization of the global attractor, and also discuss the convergence
of any bounded solutions to steady states. In particular, under suitable
assumptions on the parameters of the problem, we establish the convergence
of the bounded solution $u\left( t\right) $ to a single steady state $%
u_{\ast }$ and the rate of convergence:%
\begin{equation*}
\left\Vert u\left( t\right) -u_{\ast }\right\Vert _{L^{p}\left( \Omega
\right) }\sim \left( 1+t\right) ^{-\rho },\text{ as }t\rightarrow \infty ,
\end{equation*}%
for any $p>1,$ and some $\rho =\rho \left( u_{\ast },p\right) \in \left(
0,1\right) $. Finally, the existence of an exponential attractor is also
demonstrated for sufficiently smooth kernels in the case of non-degenerate
diffusion. Our analysis extends and complements the analysis from \cite{BRB}
and many other fundamental works.
\end{abstract}

\tableofcontents


\section{Introduction}

In this article, we are interested in the asymptotic behavior, as time goes
to infinity, of general aggregation models with (non)degenerate diffusion of
the form%
\begin{equation}
\partial _{t}u+\nabla \cdot (u\overrightarrow{V})=\Delta A\left( u\right) 
\text{, in }\Omega \times \left( 0,\infty \right) ,  \label{be1}
\end{equation}%
where $\overrightarrow{V}=\nabla \mathcal{K}\ast u$, and $\Omega $ is a
bounded domain in $\mathbb{R}^{d},$ $d\geq 1$. The kernel $\mathcal{K}$
incorporates the sensing range and degradation for the particular population
density $u$ under consideration, while the term on the right-hand side of (%
\ref{be1}) models the dispersal mechanism, such as, local repulsion. We aim
to investigate (\ref{be1}) with both no-flux boundary conditions and
Dirichlet boundary conditions. More precisely, we wish to consider the
situation in which the boundary $\partial \Omega $ of the domain $\Omega $
consists of two disjoint open subsets $\Gamma _{N}\neq \varnothing $ and $%
\Gamma _{D}$ (possibly empty), each $\overline{\Gamma }_{i}\backprime \Gamma
_{i}$ ($i\in \left\{ D,N\right\} $) is a $\sigma $-null subset of $\partial
\Omega $ and $\partial \Omega =\overline{\Gamma }_{N}\cup \overline{\Gamma }%
_{D}$; Here, $\sigma $ denotes the restriction to $\partial \Omega $ of the $%
(d-1)$-dimensional Hausdorff measure which coincides with the Lebesgue
surface measure, if we assume that $\partial \Omega $ is at least Lipschitz.
Thus, we consider the following boundary conditions for (\ref{be1}):%
\begin{equation}
\left\{ 
\begin{array}{ll}
(\nabla A\left( u\right) -u\overrightarrow{V})\cdot \overrightarrow{n}=0%
\text{,} & \text{on }\Gamma _{N}\times \left( 0,\infty \right) , \\ 
u=0, & \text{on }\Gamma _{D}\times \left( 0,\infty \right) ,%
\end{array}%
\right.  \label{be2}
\end{equation}%
and initial condition%
\begin{equation}
u_{\mid t=0}=u_{0}\text{ in }\Omega .  \label{be3}
\end{equation}

The physical motivation for taking boundary conditions as in (\ref{be2}) is
clear. For equation (\ref{be1}), the homogeneous boundary condition on $%
\Gamma _{D}$ may be interpreted as if the species suffers extinction if say
the patch $\Gamma _{D}\subset \partial \Omega $ where the individuals live
is toxic. The no-flux boundary condition in (\ref{be2}) says that nothing
can cross the boundary $\Gamma _{N}$ . Dirichlet boundary conditions for the
population density $u$ can also arise for reaction-diffusion systems in the
modelling of competition between two population species whose interaction
occurs mainly in a region where their habitats overlap. This gives rise to
Dirichlet boundary conditions for either species on the whole $\partial
\Omega $ or only on a part $\Gamma _{D}$ of $\partial \Omega $ (see \cite%
{MN, ST} for the biological literature; cf. also \cite{BB, HMM}, for some
mathematical results).

We aim to give some results which allow to deduce the $L^{p}-L^{\infty },$
and then the $L^{\infty }-C^{\alpha }\left( \overline{\Omega }\right) $
smoothing properties for solutions of (\ref{be1})-(\ref{be3}) assuming that
some sort of energy estimate is apriori known in $L^{p}$-norm for some
finite $p$. The main tool will be an iterative argument following a
well-known Alikakos-Moser technique combined with a suitable form of
Gronwall's inequality and a refined ODE\ argument. Our estimates are much
stronger than those obtained in \cite{BS, BRB} since our constants are
uniform with respect to time and the initial data. It is well-known that the
above smoothing properties become essential tools in attractor theory where
they can be used to establish the existence of an absorbing set in $%
C^{\alpha }\left( \overline{\Omega }\right) $-norm if this property can be
deduced easily in $L^{p}$-norm for some finite $p$. Recall that a subset $%
\mathcal{V}\subset \mathcal{H}$, where $\mathcal{H}$ is a topological space
endowed with a given metric, is called \emph{absorbing} if the orbits
corresponding to bounded sets $B$\ of initial data enter into $\mathcal{V}$
after a certain time (which may depend on the set $B$)\ and will stay there
forever. If the space $\mathcal{V}$ is further compactly embedded in $%
\mathcal{H}$, then the existence of the global attractor for a system like (%
\ref{be1})-(\ref{be3}) follows from standard abstracts results (see \cite%
{PZ, T}). The global attractor for (\ref{be1})-(\ref{be3}) encodes all the
information about the the long-time behavior of solutions to (\ref{be1})-(%
\ref{be3}) departing from bounded sets of initial data. Our aim is to
establish such a result for our system. We will consider diffusions like $%
A\left( y\right) \sim y^{m}$, for some $m>m_{\ast }\geq 1$ (the subcritical
case). This scenario is in complete agreement with biological observation,
that only in the case of degenerate diffusion there are solutions\ of (\ref%
{be1})-(\ref{be3}) which have compact support, steep edges, and a constant
internal population density (cf. \cite{TBL}). We shall also emphasize the
role of boundary conditions, see (\ref{be2}), and how they affect various
dissipative estimates. In fact, when $\sigma \left( \Gamma _{D}\right) =0$,
and $m_{\ast }=1+1/\gamma -2/d,$ for any $1\leq \gamma \leq d/2,$ we extend
the results in \cite{BRB} to show the existence of a global attractor,
bounded in $C^{\alpha }\left( \overline{\Omega }\right) ,$ for some $\alpha
\in \left( 0,1\right) ,$ as long as $m>m_{\ast }$. In the critical case $%
m=m_{\ast }$, the same result is valid provided that for $\int_{\Omega
}u_{0}dx=\widetilde{M}$ (see (\ref{CM}) below), we have $\widetilde{M}%
<M_{c}, $ where $M_{c}$ is the critical mass estimated in \cite[Theorem 7
and Proposition 3]{BRB}. However, when $\sigma \left( \Gamma _{D}\right) >0$
it appears that we cannot recover the critical exponent $m_{\ast }$ as
above. We can only show the existence of the global attractor, bounded in $%
C^{\alpha }\left( \overline{\Omega }\right) ,$ provided that $m>m_{\ast }=2$
in this case. This appears to be due to \emph{loss} of the conservation of
mass property in the associated non-degenerate problem (i.e., if $A^{\prime
}\left( y\right) \geq \varepsilon >0$, $\forall y\geq 0$) which becomes
crucial in obtaining uniform apriori estimates for the degenerate problem (%
\ref{be1})-(\ref{be3}). Finally, as in \cite{BRB} we give optimal
assumptions on the interaction kernel $\mathcal{K}$ which include important
cases of interest, such as, the Newtonian and Bessel potentials for $d\geq 2$%
. It is interesting to note that when $\mathcal{K}\left( x\right) =\mathcal{N%
}\left( x\right) $ is Newtonian $\mathcal{N}\left( x\right) =C_{d}\left\vert
x\right\vert ^{-\left( d-2\right) },$ for $d\geq 3$, or $\mathcal{N}\left(
x\right) =-C\log \left( \left\vert x\right\vert \right) ,$ for $d=2$, we
recover in (\ref{be1})-(\ref{be3}) the classical parabolic-elliptic
Patlak--Keller--Segel (PKS) model of chemotaxis:%
\begin{equation}
\left\{ 
\begin{array}{l}
\partial _{t}u+\nabla \cdot (u\nabla c)=\Delta A\left( u\right) \text{, in }%
\Omega \times \left( 0,\infty \right) , \\ 
c=\mathcal{N}\ast u,\text{ in }\Omega \times \left( 0,\infty \right) ,%
\end{array}%
\right.  \label{pePKS}
\end{equation}%
that is, the second equation in (\ref{pePKS}) reads $-\Delta c=u$. We could
not find a proof for the existence of the global attractor, and its
properties for the PKS model (\ref{pePKS})\ in the literature. Our results
cover this important case as well. For the precise statements of the
results, we refer the reader to Sections 3-4. Another basic system of
equations modeling chemotaxis was established in the 1970's by Keller and
Segel \cite{KS1, KS2, KS3}. In these systems, the population concentration $%
u $ satisfies the first equation of (\ref{pePKS}), while the chemotactic
agent $c$ satisfies instead the parabolic equation:%
\begin{equation}
\tau \partial _{t}c-\Delta c=\beta _{1}c+\beta _{2}u,\text{ in }\Omega
\times \left( 0,\infty \right) ,  \label{ppPKS}
\end{equation}%
for some real constants $\beta _{1},\beta _{2}$ and $\tau >0.$ The latter
system is usually referred in the literature as the parabolic-parabolic
Patlak--Keller--Segel model of chemotaxis. The mathematical and biological
literature, concerning primarily with the qualitative properties of the
solutions to the parabolic-parabolic PKS (and some of its generalizations),
is quite extensive and much of the work before 2003 is largely referenced in
the survey papers by Horstmann \cite{H1, H2}. More recent results pertaining
to the long-term behavior of solutions, in terms of global and exponential
attractors, of the parabolic-parabolic PKS model\ can be found in \cite{AA1,
AA2, AA3, ENO}. The issue of the convergence to single stationary states as
time goes to infinity is also addressed in \cite{FLP, ZZ}. Finally, it is
worth observing that the parabolic-elliptic PKS model (\ref{pePKS}) follows
for an appropriate choice of the parameters $\beta _{1},\beta _{2}$, not
only as a formal limit but also rigorously as $\tau \rightarrow 0^{+},$ from
the parabolic-parabolic PKS model \cite{BB2}.

Now, we wish to provide the reader with some further background on the above
system (\ref{be1})-(\ref{be3}). The whole issue of well-posedness of weak
solutions for equation (\ref{be1}) with no-flux boundary conditions on a
bounded convex domain of class $\mathcal{C}^{1}$, and sufficiently smooth
interaction kernels $\mathcal{K}$ (see below for the precise assumptions),
was established in \cite{BS, BRB} (see also \cite{BCM} for related results).
When $\mathcal{K}$ is not smooth enough, finite time blowup of some
solutions can occur (see, for instance, \cite{BB, BCL, BL, BLR, BRB}).
However, it is worth emphasizing that in population dynamics, the non-local
effects are generally modelled with smooth, fast-decaying kernels. When
problem (\ref{be1})-(\ref{be3}) is uniformly parabolic (i.e., $A^{^{\prime
}}\left( y\right) \geq \varepsilon >0$, for all $y\geq 0$), we can improve
our analysis from degenerate diffusion. More precisely, we give a complete
characterization of the global attractor $\mathcal{A}_{\varepsilon }$ in
this case, as the union of all unstable manifolds generated by all
equilibria (steady-state)\ solutions of the non-degenerate aggregation
equation (see Section \ref{non_deg}, Theorem \ref{unst_gl}). At this point
one could argue that the long-time behavior of the system (\ref{be1})-(\ref%
{be3}) with nondegenerate diffusion ($A^{^{\prime }}\left( y\right) \geq
\varepsilon >0$)\ is properly described by the global attractor. However, it
is well-known that the global attractor can present several drawbacks, among
which we can mention that it may only attract the trajectories at a slow
rate, and that it may miss important transient behaviors because the global
attractor consists only of states in a final stage. This phenomenon is
already present for models of pattern formation in chemotaxis (see \cite{TT2}%
). Another suitable object which contains the global attractor, and thus is
more rich in content than the global attractor is the so-called exponential
attractor (see \cite{MZ, TT2}; cf. also Section \ref{non_deg}). In Section %
\ref{non_deg} we show that the problem (\ref{be1})-(\ref{be3}) with
nondegenerate diffusion ($A^{^{\prime }}\left( y\right) \geq \varepsilon >0$%
)\ admits also an \emph{exponential} attractor $\mathcal{M}_{\varepsilon }$
(and as a result, the global attractor $\mathcal{A}_{\varepsilon }$ is \emph{%
finite-dimensional}), globally bounded in $C^{\alpha }\left( \overline{%
\Omega }\right) $, for some $\alpha \in \left( 0,1\right) $, provided we
assume, in addition, that the kernel $\mathcal{K}$ is sufficiently smooth at
the origin, i.e.,%
\begin{equation*}
\mathcal{K}\in W^{2,1}\left( B_{1}\left( 0\right) \right) .
\end{equation*}%
Here $B_{1}\left( 0\right) \subset \mathbb{R}^{d}$ is the ball centered at
the origin and radius equal to one.

In the final Section \ref{sta}, we also discuss the convergence of any
bounded solutions $u\left( t\right) $\ of the non-degenerate aggregation
equation to single steady states, provided that $\Phi $ is a real analytic
function on $\mathbb{R}_{+}$, where $\Phi ^{^{\prime \prime }}\left(
y\right) :=A^{^{\prime }}\left( y\right) /y$ and $\Phi \left( 0\right) =\Phi
^{^{\prime }}\left( 1\right) =0.$ Moreover, we also establish the
convergence rate of the bounded solution $u\left( t\right) $ to a single
steady state $u_{\ast }$:%
\begin{equation*}
\left\Vert u\left( t\right) -u_{\ast }\right\Vert _{L^{p}\left( \Omega
\right) }\sim \left( 1+t\right) ^{-\rho },\text{ as }t\rightarrow \infty ,
\end{equation*}%
for any $p>1,$ and some $\rho =\rho \left( u_{\ast },p\right) \in \left(
0,1\right) $ (see Theorem \ref{conv_equil}). We refer the reader to Section %
\ref{sta} for the precise assumptions, statements and further discussions.

As noted in \cite{TBL}, the system (\ref{be1})-(\ref{be3}) exhibits
interesting coarsening dynamics whose behavior is similar to another
well-known (non-biological) model, the so-called nonlocal Cahn-Hilliard
equation, which also exhibits behaviors in which small localized clumps form
and merge into larger clumps over time. The latter equation has also been
recently studied in \cite{BH, BH2, GL1, GL2, Ha, LP}. We refer the reader
for more related results concerning the nonlocal Cahn-Hilliard equation to 
\cite{GM}, where a complete characterization of the long-term behavior is
also given for this equation. Our proofs for the existence of attractors and
their properties will explore various connections which exist between the
aggregation equation (\ref{be1})\ and the non-local Cahn-Hilliard equation
(see \cite{GM}, and references therein).

\section{Weak solutions}

We begin with some basic notations and preliminaries. Throughout the
section, $C\geq 0$ will denote a \emph{generic} constant, while $Q:\mathbb{R}%
_{+}\rightarrow \mathbb{R}_{+}$ will denote a \emph{generic} increasing
function. All these quantities, unless explicitly stated, are \emph{%
independent} of time, an approximation parameter $\varepsilon >0$ and the
initial data. Further dependencies of these quantities will be specified on
occurrence. In the sequel, our investigation will be mainly divided into two
cases:

\begin{itemize}
\item (i) $\sigma \left( \Gamma _{D}\right) =0$, i.e., $\Gamma _{D}$ is
empty, and

\item (ii) $\sigma \left( \Gamma _{D}\right) >0,$ i.e., $\Gamma _{D}$ is an
open set of positive $\sigma $-measure,
\end{itemize}

\noindent where we recall that $\sigma $ denotes the restriction to $%
\partial \Omega $ of the $(d-1)$-dimensional Hausdorff measure.

Following \cite[Section 2]{BS}, \cite[Definition 1]{BRB} we make the
following assumptions on $A$ and $\mathcal{K}$:

\begin{description}
\item[H1] The domain $\Omega \subset \mathbb{R}^{d}$, $d\geq 1,$ is convex
and of class $\mathcal{C}^{1}.$

\item[H2] $A\in C^{1}[0,\infty ),$ $A\left( 0\right) =0$, and there exists
constants $C_{A},\overline{C}_{A}>0$ such that 
\begin{equation*}
C_{A}y^{m-1}\leq A^{^{\prime }}\left( y\right) \leq \overline{C}_{A}y^{m-1},%
\text{ }\forall y\geq 0,
\end{equation*}%
for some $m>m_{\ast }$, with%
\begin{equation*}
m_{\ast }:=\left\{ 
\begin{array}{ll}
1+\frac{1}{\gamma }-\frac{2}{d}\geq 1, & \text{if }\sigma \left( \Gamma
_{D}\right) =0, \\ 
2, & \text{if }\sigma \left( \Gamma _{D}\right) >0,%
\end{array}%
\right.
\end{equation*}%
for any $1\leq \gamma \leq \frac{d}{2}$. For $\sigma \left( \Gamma
_{D}\right) =0$, $m_{\ast }$ is the critical exponent defined in \cite[%
Definition 5 and Lemma 10]{BRB}.

\item[H3] $\mathcal{K}\in W_{\text{loc}}^{1,1}\left( \mathbb{R}^{d}\right)
\cap C^{3}(\mathbb{R}^{d}\backslash \left\{ 0\right\} )$ satisfies the
following assumptions:

\begin{enumerate}
\item[(i)] $\mathcal{K}$ is symmetric, $\mathcal{K}\left( x\right) =k\left(
\left\vert x\right\vert \right) $ and $k$ is nonincreasing.

\item[(ii)] $k^{^{\prime \prime }}\left( r\right) $ and $k^{^{\prime
}}\left( r\right) /r$ are monotone on $r\in \left( 0,\delta \right) ,$ for
some $\delta >0.$

\item[(iii)] $\left\vert D^{3}\mathcal{K}\left( x\right) \right\vert \leq
C\left\vert x\right\vert ^{-d-1},$ for some $C>0.$
\end{enumerate}
\end{description}

Note that since the function $k$ in condition (H3) is nonincreasing, the
nonlocal term in (\ref{be1}) models attraction. Moreover, these conditions
imply that if $\mathcal{K}$ is singular, the singularity is restricted to
the origin, so that both the Newtonian and Bessel potentials for $d\geq 2$
are included in our analysis. Finally, assumption (H2) implies that our
problem (\ref{be1})-(\ref{be3}) is \emph{subcritical} in the terminology of 
\cite[Definition 6]{BRB}, and note that the total population $u$ is
preserved in time. In particular, there holds%
\begin{equation}
\left\langle u\left( t\right) \right\rangle :=\frac{1}{\left\vert \Omega
\right\vert }\int_{\Omega }u\left( t,x\right) ds=\left\langle
u_{0}\right\rangle ,\text{ for all }t\geq 0.  \label{CM}
\end{equation}%
Here, $\left\vert \Omega \right\vert $ denotes the $d$-dimensional Lebesgue
measure of $\Omega $. This property is lost whenever $\Gamma _{D}$ is
nonempty, $\sigma \left( \Gamma _{D}\right) >0$ and $A$ is non-degenerate
such that $A^{^{\prime }}\left( s\right) \geq \varepsilon >0$.

Let us recall some important properties of the kernel $\mathcal{K}$, proven
in \cite[Section 1.3]{BRB}, which imply a useful number of estimates in weak 
$L^{p,\infty }$-spaces, with quasi-norm 
\begin{equation*}
\left\Vert u\right\Vert _{L^{p,\infty }}=\sup_{\alpha >0}\left( \alpha
^{p}\lambda _{\alpha }\left( u\right) \right) ^{1/p},
\end{equation*}%
where $\lambda _{\alpha }\left( u\right) :=\left\vert \left\{ u>\alpha
\right\} \right\vert $ is the distribution function of $u$ (see \cite[Lemmas
1-3]{BRB}).

\begin{lemma}
\label{kest}Let $\mathcal{K}$ satisfy assumption (H3), (i)-(iii).

(a) There holds: $\nabla \mathcal{K}\in L^{\frac{d}{d-1},\infty }\left(
\Omega \right) $. If $d\geq 3$, then $\mathcal{K}\in L^{\frac{d}{d-2},\infty
}\left( \Omega \right) .$

(b) Let $\overrightarrow{V}=\nabla \mathcal{K}\ast u.$ Then, for any $%
1<p<\infty $, there exists a constant $C=C\left( p\right) >0$ such that%
\begin{equation*}
\left\Vert \nabla \overrightarrow{V}\right\Vert _{L^{p}\left( \Omega \right)
}\leq C\left\Vert u\right\Vert _{L^{p}\left( \Omega \right) }.
\end{equation*}%
(c) Let $d\geq 3.$ Suppose that $\gamma \in \left( 1,d/2\right) .$ Then $%
\mathcal{K}\in L_{\text{loc}}^{\frac{d}{d/\gamma -2},\infty }\left( \Omega
\right) $ if and only if $D^{2}\mathcal{K}\in L_{\text{loc}}^{\gamma ,\infty
}\left( \Omega \right) $. The same holds for $\nabla \mathcal{K}\in L_{\text{%
loc}}^{\frac{d}{d/\gamma -1},\infty }\left( \Omega \right) .$
\end{lemma}

\begin{remark}
The above lemma shows that every admissible kernel $\mathcal{K}$, that
satisfies (H3), (i)-(iii), is at least as integrable as the Newtonian
potential.
\end{remark}

As in \cite[Definition 2.1]{BS} and \cite[Definition 3]{BRB}, our notion of
bounded, nonnegative weak solutions to (\ref{be1})-(\ref{be3}) is as
follows. By $H_{D}^{1}\left( \Omega \right) $ we denote the space of all
functions $u\in H^{1}\left( \Omega \right) $ such that $u=0$ on $\Gamma
_{D}, $ if $\sigma \left( \Gamma _{D}\right) >0$, and by $\left(
H_{D}^{1}\left( \Omega \right) \right) ^{\ast }$ the dual of $%
H_{D}^{1}\left( \Omega \right) .$ Below, the space $H_{D}^{1}\left( \Omega
\right) $ is replaced by $H^{1}\left( \Omega \right) $ whenever $\sigma
\left( \Gamma _{D}\right) =0$ (i.e., if $\Gamma _{D}$ is empty).

\begin{definition}
\label{weaksol}Let $T>0$ be given, but otherwise arbitrary. A function $%
u:\Omega \times \left[ 0,T\right] \rightarrow \lbrack 0,\infty )$ is a weak
solution of (\ref{be1})-(\ref{be3}) if%
\begin{align}
u& \in L^{\infty }\left( \Omega \times \left[ 0,T\right] \right) ,\text{ }%
A\left( u\right) \in L^{2}\left( 0,T;H_{D}^{1}\left( \Omega \right) \right) ,
\label{reg} \\
\partial _{t}u& \in L^{2}\left( 0,T;\left( H_{D}^{1}\left( \Omega \right)
\right) ^{\ast }\right) ,\text{ }u\nabla \mathcal{K}\ast u\in L^{2}\left(
\Omega \times \left[ 0,T\right] \right)  \notag
\end{align}%
and the following identity holds:%
\begin{equation}
\left\langle \partial _{t}u\left( t\right) ,w\right\rangle _{\left(
H_{D}^{1}\left( \Omega \right) \right) ^{\ast },H_{D}^{1}\left( \Omega
\right) }+\int_{\Omega }\nabla A\left( u\left( t\right) \right) \cdot \nabla
w-u\left( t\right) \left( \nabla \mathcal{K}\ast u\left( t\right) \right)
\cdot \nabla wdx=0,  \label{identity}
\end{equation}%
for all $w\in H_{D}^{1}\left( \Omega \right) $, for almost all $t\in \left[
0,T\right] .$
\end{definition}

\begin{remark}
Note that by (\ref{reg}), $u\in C(0,T;\left( H_{D}^{1}\left( \Omega \right)
\right) ^{\ast })$ such that the initial condition in (\ref{be3}) is
understood in the weak sense of $\left( H_{D}^{1}\left( \Omega \right)
\right) ^{\ast }$. In fact, in \cite[Section 2]{BS} it is shown that $u\in
C\left( 0,T;L^{p}\left( \Omega \right) \right) $, for all $1<p<\infty $,
provided that $0\leq u_{0}\in L^{\infty }\left( \Omega \right) ,$ so that
the initial condition is also satisfied in the $L^{p}$-sense.
\end{remark}

We have the following result concerning well-posedness of the system (\ref%
{be1})-(\ref{be3}).

\begin{theorem}
\label{well-posed}Let the assumptions (H1)-(H3) be satisfied in both cases
(i)+(ii) (i.e., $\sigma \left( \Gamma _{D}\right) \geq 0$), and assume $%
0\leq u_{0}\in L^{\infty }\left( \Omega \right) .$ Then there exists a
unique (global) nonnegative solution to problem (\ref{be1})-(\ref{be3}) in
the sense of Definition \ref{weaksol}. Moreover, every weak solution
satisfies the following dissipative inequality:%
\begin{equation}
\mathcal{E}\left( u\left( t\right) \right) +\int_{0}^{t}\int_{\Omega
}u\left( s\right) \left\vert \nabla \Phi ^{^{\prime }}\left( u\left(
s\right) \right) -\nabla \mathcal{K}\ast u\left( s\right) \right\vert
^{2}dxds\leq \mathcal{E}\left( u_{0}\right) ,  \label{en_ineq}
\end{equation}%
for all $t\geq 0$, where%
\begin{equation}
\mathcal{E}\left( u\left( t\right) \right) :=\int_{\Omega }\Phi \left(
u\left( t\right) \right) dx-\frac{1}{2}\int_{\Omega }\int_{\Omega }u\left(
x,t\right) \mathcal{K}\left( x-y\right) u\left( y,t\right) dxdy,
\label{energy}
\end{equation}%
$\Phi $ is strictly convex on $\left( 0,\infty \right) $ such that $\Phi
^{^{\prime \prime }}\left( y\right) :=A^{^{\prime }}\left( y\right) /y$ and $%
\Phi \left( 0\right) =\Phi ^{^{\prime }}\left( 1\right) =0.$
\end{theorem}

\begin{proof}
We only briefly mention the main steps in the proof. The uniqueness of weak
solutions follows from \cite[Theorem 2.4]{BS} (cf. also \cite[Theorem 3]{BRB}
for a more general result) with some minor modifications (see Lemma \ref{cdp}
below) since the boundary terms involving convolutions with the kernel $%
\mathcal{K}$ vanish on $\Gamma _{D}$, if $\sigma \left( \Gamma _{D}\right)
>0.$

\emph{Step 1} (Local existence). The \emph{local} existence result can be
carried out in a standard manner, by first regularizing $A\left( y\right) $
with $A_{\varepsilon }\left( y\right) =A\left( y\right) +\varepsilon y,$ $%
\forall \varepsilon >0,$ $\mathcal{K}$ by a sequence of smooth kernels $%
\mathcal{J}_{\varepsilon }\mathcal{K}$ (where $\mathcal{J}_{\varepsilon }$
is a standard mollifier), and then the initial data $u_{\varepsilon }\left(
0\right) =u_{0,\varepsilon }\in C^{\infty }\left( \Omega \right) \cap
C^{1}\left( \overline{\Omega }\right) ,$ followed by passage to limit as $%
\varepsilon \rightarrow 0^{+}$ in the corresponding regularized problem $%
P_{\varepsilon }$:%
\begin{equation}
\left\{ 
\begin{array}{ll}
\partial _{t}u_{\varepsilon }+\nabla \cdot (u_{\varepsilon }\nabla \mathcal{J%
}_{\varepsilon }\mathcal{K}\ast u_{\varepsilon })=\Delta A_{\varepsilon
}\left( u_{\varepsilon }\right) \text{,} & \text{in }\Omega \times \left(
0,\infty \right) , \\ 
(\nabla A_{\varepsilon }\left( u_{\varepsilon }\right) -u_{\varepsilon
}\nabla \mathcal{J}_{\varepsilon }\mathcal{K}\ast u_{\varepsilon })\cdot 
\overrightarrow{n}=0\text{,} & \text{on }\Gamma _{N}\times \left( 0,\infty
\right) , \\ 
u_{\varepsilon }=0, & \text{on }\Gamma _{D}\times \left( 0,\infty \right) ,%
\end{array}%
\right.  \label{approx}
\end{equation}%
(see, \cite[Theorem 2.13]{BS}, \cite[Section 3]{BRB}). This is achieved by
proving a series of uniform (in $\varepsilon >0$) estimates for the
approximate solutions $u_{\varepsilon }\left( t\right) ,$ in particular, by
establishing the uniform bound in $L^{\infty }\left( \Omega \right) $-norm:%
\begin{equation}
\sup_{t\geq 0}\left\Vert u_{\varepsilon }\left( t\right) \right\Vert
_{L^{\infty }\left( \Omega \right) }\leq C,  \label{lsup}
\end{equation}%
for some positive constant $C$ which is independent of $\varepsilon >0$ (and
even the times $t,T$). This norm gives further uniform estimates (with
respect to $\varepsilon >0$) for solutions in the normed spaces of (\ref{reg}%
). More precisely, we can obtain the following bounds:%
\begin{align}
\left\Vert A\left( u_{\varepsilon }\right) \right\Vert _{L^{2}\left(
0,T;H_{D}^{1}\left( \Omega \right) \right) }& \leq C, \\
\left\Vert \partial _{t}u\right\Vert _{L^{2}\left( 0,T;\left(
H_{D}^{1}\left( \Omega \right) \right) ^{\ast }\right) }& \leq C,\text{ } \\
\left\Vert u\nabla \mathcal{K}\ast u\right\Vert _{L^{2}\left( \Omega \times %
\left[ 0,T\right] \right) }& \leq C,
\end{align}%
for some positive constant $C$ independent of $\varepsilon $. Moreover, for
every $\varepsilon >0$ it holds%
\begin{equation}
\left\Vert u_{\varepsilon }\right\Vert _{L^{2}\left( 0,T;H_{D}^{1}\left(
\Omega \right) \right) }\leq \frac{C}{\sqrt{\varepsilon }}.  \label{h1norm}
\end{equation}%
Thus, on account of the above uniform estimates the sequence of solutions $%
u_{\varepsilon }$ is precompact in $L^{p}\left( \Omega \times \left[ 0,T%
\right] \right) ,$ for all $1\leq p<\infty $, (cf. \cite[Lemma 9]{BRB}, \cite%
{BS}), and we can pass to the limit as $\varepsilon \rightarrow 0$ in the
nonlinear terms in (\ref{approx}). In particular, it can be shown that%
\begin{equation*}
A\left( u_{\varepsilon }\left( t\right) \right) \rightarrow A\left( u\left(
t\right) \right) \text{ weakly in }L^{2}\left( 0,T;H_{D}^{1}\left( \Omega
\right) \right)
\end{equation*}%
and 
\begin{equation}
u_{\varepsilon }\nabla \mathcal{J}_{\varepsilon }\mathcal{K}\ast
u_{\varepsilon }\rightarrow u\nabla \mathcal{K}\ast u\text{ weakly in }%
L^{2}\left( \Omega \times \left[ 0,T\right] \right)  \label{finp}
\end{equation}%
(see, \cite[Theorem 1]{BRB}). Let us now comment how to get the estimates (%
\ref{lsup})-(\ref{h1norm}) in each of the following cases: $\sigma \left(
\Gamma _{D}\right) =0$ or $\sigma \left( \Gamma _{D}\right) >0.$

\begin{itemize}
\item \textbf{Case (i)}: $\sigma \left( \Gamma _{D}\right) =0$ (i.e., $%
\Gamma _{D}$ is empty). The argument for deducing (\ref{lsup}) relies on
showing that the $L^{1}$-norm of $u_{\varepsilon }\left( t\right) $ (see (%
\ref{CM})) apriori controls the $L^{p}$-norm for $1<p<\infty $, and then
that the latter norm controls the $L^{\infty }$-norm of $u_{\varepsilon
}\left( t\right) $, see (\ref{lsup}) (cf. also \cite[Lemma 8]{BRB}). The
Lemma \ref{kest} is crucial for the proof of the uniform bound (\ref{lsup}).
Under the assumptions (H1)-(H3) the bounds (\ref{lsup})-(\ref{h1norm}) were
obtained in \cite[Section 3]{BRB} (cf. also \cite{BS}).

\item \textbf{Case (ii)}: $\sigma \left( \Gamma _{D}\right) >0$. The
arguments from \cite{BRB, BS} leading to estimates (\ref{lsup})-(\ref{lsup})
on intervals $\left( 0,T\right) ,$ for some $T=T\left( p\right) >0,$ seem to
break down if $\Gamma _{D}$ is nonempty and $\sigma \left( \Gamma
_{D}\right) >0$ since there is no conservation of mass (\ref{CM}) for (\ref%
{approx}) in this case; hence, the $L^{1}$-norm of $u_{\varepsilon }\left(
t\right) $ is not controlled apriori. However, even in this case, the same
proof of obtaining local-in-time estimates for $u_{\varepsilon }\left(
t\right) $ in \cite[Lemma 8]{BRB} applies, but we need to use the following
inequality%
\begin{equation}
\left\Vert u_{\varepsilon }\left( t\right) \right\Vert _{L^{p}\left( \Omega
\right) }^{p}\leq 2^{p}\left( \left\Vert \left( u_{\varepsilon }\left(
t\right) -k\right) _{+}\right\Vert _{L^{p}\left( \Omega \right)
}^{p}+k^{p}\left\vert \Omega \right\vert \right) ,\text{ for any }k>0,
\label{ukest}
\end{equation}%
instead of the usual one \cite[Ineq. (25)]{BRB}. Then, by arguing as in \cite%
[Ineq. (27) and Lemma 8, Step 2]{BRB} (note that boundary terms involving
convolutions with the kernel $\mathcal{K}$ vanish on $\Gamma _{D}$, if $%
\sigma \left( \Gamma _{D}\right) >0$; also now the last constant on the
right-hand side of \cite[Ineq. (27)]{BRB} depends only on $k,p,$ and not on $%
L^{1}$-norm of $u_{\varepsilon }\left( 0\right) $), we also obtain a uniform
in $\varepsilon >0,$ $L^{\infty }$-bound for $u_{\varepsilon }\left(
t\right) ,$ $t\in \left[ 0,T\right] ,$ for some $T=T\left( p\right) >0$.
\end{itemize}

By virtue of these observations, the proof of the \emph{local} existence
argument in both cases \textbf{(i)}+\textbf{(ii)} goes exactly as in \cite[%
Section 3]{BRB} (cf. also \cite[Section 2.2]{BS}). A proof of the energy
inequality (\ref{en_ineq}) can be found in \cite[Lemma 5.1]{BS}, \cite[%
Proposition 1]{BRB}. Continuation of the weak solutions is a straightforward
consequence of the local existence theory (in Step 1) and the proof of \cite[%
Theorem 4]{BRB}.

\emph{Step 2} (Global existence).

\begin{itemize}
\item First, let $\sigma \left( \Gamma _{D}\right) =0$. Since the problem is
subcritical by assumption (H2), the global existence result follows from 
\cite[Lemma 10 and Remark 9]{BRB} which shows that (\ref{lsup}) is also
satisfied uniformly with respect to $t,T$.

\item When $\sigma \left( \Gamma _{D}\right) >0$, it suffices to establish a
uniform (in $\varepsilon >0$ and time $t,T>0$) $L^{p}$-estimate for $%
u_{\varepsilon }\left( t\right) $ and to exploit an argument similar to \cite%
[Lemma 8, Step 2]{BRB} to deduce the uniform (with respect to $%
t,T,\varepsilon $) bound (\ref{lsup}). However, the key difference with
respect to the case $\sigma \left( \Gamma _{D}\right) =0$ is that mass is 
\emph{not} conserved for (\ref{approx}); hence, the $L^{1}$-norm of $%
u_{\varepsilon }\left( t\right) $ is not controlled apriori and we need a
different argument to deduce the uniform bound in $L^{p}$-norm for $%
u_{\varepsilon }\left( t\right) $, see Proposition \ref{lpbound} and Remark %
\ref{rem} below. The latter bound requires that the critical exponent $%
m_{\ast }=2$ in (H2).
\end{itemize}

The proof of Theorem \ref{well-posed} is now complete.
\end{proof}

\begin{remark}
\label{esse}In order to establish the energy identity in (\ref{en_ineq}), it
would suffice to show that $\Phi ^{^{\prime }}\left( u\right) \in
L^{2}\left( 0,T;H_{D}^{1}\left( \Omega \right) \right) $ so we can test
equation (\ref{identity}) with $\Phi ^{^{\prime }}\left( u\right) $ and $%
\mathcal{K}\ast u$, respectively. At the moment, this seems unreachable for
our problem with a degenerate diffusion function $A\left( u\right) \sim
u^{m},$ like in (H2). However, note that $\Phi \left( u\right) \in
L^{2}\left( 0,T;H_{D}^{1}\left( \Omega \right) \right) ,$ $\Phi \in
L^{\infty }\left( 0,T;L^{1}\left( \Omega \right) \right) $ since $\Phi
\left( u\right) \sim A\left( u\right) $, by (H2). When the diffusion $A$ in (%
\ref{be1}) is non-degenerate (say, like in problem (\ref{approx})), we can
establish the energy identity (see the proof of Theorem \ref{well-posed}).
Indeed, in this case we can easily check that $\Phi ^{^{\prime }}\left(
u_{\varepsilon }\right) \in L^{2}\left( 0,T;H_{D}^{1}\left( \Omega \right)
\right) ,$ $\forall \varepsilon >0,$ on account of (\ref{lsup}) and (\ref%
{h1norm}). Thus, the key multiplication of the corresponding weak
formulation associated with (\ref{approx}) (see (\ref{identity})) with $\Phi
^{^{\prime }}\left( u_{\varepsilon }\right) $ and $\mathcal{K}\ast
u_{\varepsilon }\in L^{2}\left( 0,T;H_{D}^{1}\left( \Omega \right) \right) $%
, respectively, is allowed. Exploiting, for instance, \cite[Lemma 2.6]{BS},
and owing to the convexity of $\Phi _{\varepsilon }$ we get equality in (\ref%
{en_ineq}) for the energy $\mathcal{E}_{\varepsilon }\left( u_{\varepsilon
}\left( t\right) \right) ,$ associated with (\ref{approx}), which is defined
by%
\begin{equation*}
\mathcal{E}_{\varepsilon }\left( u_{\varepsilon }\left( t\right) \right)
:=\int_{\Omega }\Phi _{\varepsilon }\left( u_{\varepsilon }\left( t\right)
\right) dx-\frac{1}{2}\int_{\Omega }\int_{\Omega }u_{\varepsilon }\left(
x,t\right) \mathcal{K}\left( x-y\right) u_{\varepsilon }\left( y,t\right)
dxdy
\end{equation*}%
($\Phi _{\varepsilon }$ is the same function as above, but with $A$ replaced
by $A_{\varepsilon }$).
\end{remark}

The next lemma gives the (H\"{o}lder) continuity of solutions with respect
to the initial data in $\left( H_{D}^{1}\left( \Omega \right) \right) ^{\ast
}$ and $\left( H^{1}\left( \Omega \right) \right) ^{\ast }$, respectively.

\begin{lemma}
\label{cdp}Let the hypotheses of Theorem \ref{well-posed} be satisfied. Let $%
u_{1}\left( t\right) $ and $u_{2}\left( t\right) $ be any two weak solutions
of (\ref{be1})-(\ref{be3}) corresponding to any two initial data $u_{01}$
and $u_{02}$, respectively. If $\sigma \left( \Gamma _{D}\right) =0$, we
further take $\left\langle u_{01}\right\rangle =M_{1}$ and $\left\langle
u_{02}\right\rangle =M_{2}$, for some $M_{1},M_{2}\geq 0$. For all $t\in %
\left[ 0,T\right] $, the following estimates hold:%
\begin{equation}
\left\Vert u_{1}\left( t\right) -u_{2}\left( t\right) \right\Vert _{\left(
H_{D}^{1}\left( \Omega \right) \right) ^{\ast }}\leq C\left( \left\Vert
u_{10}-u_{20}\right\Vert _{\left( H_{D}^{1}\left( \Omega \right) \right)
^{\ast }}\right) ^{e^{-Ct}},\text{ if }\sigma \left( \Gamma _{D}\right) >0,
\label{lip1}
\end{equation}%
and%
\begin{equation}
\left\Vert u_{1}\left( t\right) -u_{2}\left( t\right) \right\Vert _{\left(
H^{1}\left( \Omega \right) \right) ^{\ast }}\leq C\left( \left\Vert
u_{10}-u_{20}\right\Vert _{\left( H^{1}\left( \Omega \right) \right) ^{\ast
}}+\left\vert \left\langle u_{01}-u_{20}\right\rangle \right\vert \right)
^{e^{-Ct}},  \label{lip2}
\end{equation}%
if $\sigma \left( \Gamma _{D}\right) =0.$
\end{lemma}

\begin{proof}
We briefly explain how to get (\ref{lip2}); the estimate (\ref{lip1}) is
similar. Let $u\left( t\right) :=u_{1}\left( t\right) -u_{2}\left( t\right) $
and observe that (\ref{CM}) yields $\left\langle u\left( t\right)
\right\rangle =M_{1}-M_{2}=:M_{12}$, for all $t\geq 0$, but $M_{12}\neq 0$,
in general. As in \cite[Theorem 2.4]{BS}, \cite[Section 2]{BRB}, consider
the Neumann problem%
\begin{equation}
\left\{ 
\begin{array}{ll}
\Delta \phi \left( t\right) =u\left( t\right) -M_{12}\text{ in }\Omega , & 
\\ 
\nabla \phi \left( t\right) \cdot \overrightarrow{n}=0\text{ on }\Gamma , & 
\end{array}%
\right.  \label{Np}
\end{equation}%
for which $\left\langle \phi \left( t\right) \right\rangle =0$. Notice that
since $\left\langle u\left( t\right) \right\rangle =M_{12}$, (\ref{Np}) has
a solution. Consider the operator $A_{N}=-\Delta _{N}$, with domain $D\left(
A_{N}\right) =\{\varphi \in H^{2}\left( \Omega \right) :\left( \nabla
\varphi \cdot \overrightarrow{n}\right) _{\mid \Gamma }=0\}$. Clearly, $\phi
\left( 0\right) =A_{N}^{-1}\left( u\left( 0\right) -M_{12}\right) $ and
recall that, due to a Poincare inequality,%
\begin{equation*}
\left\Vert u\right\Vert _{\left( H^{1}\left( \Omega \right) \right) ^{\ast
}}^{2}=\left\Vert A_{N}^{-1/2}\left( u-\left\langle u\right\rangle \right)
\right\Vert _{L^{2}\left( \Omega \right) }^{2}+\left\langle u\right\rangle
^{2}.
\end{equation*}%
Next, we see that that $\partial _{t}\phi \left( t\right) $ also satisfies
(in the generalized sense) the problem%
\begin{equation}
\left\{ 
\begin{array}{ll}
\Delta \partial _{t}\phi \left( t\right) =\partial _{t}u\left( t\right) 
\text{ in }\Omega , &  \\ 
\nabla \partial _{t}\phi \left( t\right) \cdot \overrightarrow{n}=0\text{ on 
}\Gamma , & 
\end{array}%
\right.
\end{equation}%
and 
\begin{equation*}
\phi \in L^{\infty }\left( \Omega \times \left[ 0,T\right] \right) \cap
C\left( 0,T;H_{D}^{1}\left( \Omega \right) \right) ,\nabla \phi \in
L^{\infty }\left( \Omega \times \left[ 0,T\right] \right) .
\end{equation*}%
Thus, arguing in a standard way as in \cite[(6)-(12)]{BS}, we obtain for%
\begin{equation*}
\eta \left( t\right) :=\left\Vert u\left( t\right) \right\Vert _{\left(
H^{1}\left( \Omega \right) \right) ^{\ast }}^{2}=\left\Vert \nabla \phi
\left( t\right) \right\Vert _{L^{2}\left( \Omega \right) }^{2}+\left(
M_{12}\right) ^{2},
\end{equation*}%
the following inequality:%
\begin{align}
\frac{1}{2}\frac{d}{dt}\eta \left( t\right) & =\left\langle \nabla \phi
\left( t\right) ,\partial _{t}\nabla \phi \left( t\right) \right\rangle
\label{ineqdiff} \\
& =-\left\langle \partial _{t}u\left( t\right) ,\phi \left( t\right)
\right\rangle =I_{1}+I_{2}+I_{3},  \notag
\end{align}%
where%
\begin{align}
I_{1}& :=\int_{\Omega }\nabla A\left( u_{1}\left( t\right) \right) -\nabla
A\left( u_{2}\left( t\right) \right) \cdot \nabla \phi \left( t\right) dx,
\label{Is} \\
I_{2}& :=-\int_{\Omega }\left( u\left( t\right) \right) \nabla \mathcal{K}%
\ast u_{1}\left( t\right) \cdot \nabla \phi \left( t\right) dx,  \notag \\
I_{3}& :=-\int_{\Omega }u_{2}\left( t\right) \nabla \mathcal{K}\ast \left(
u\left( t\right) \right) \cdot \nabla \phi \left( t\right) dx.  \notag
\end{align}%
Since $A$ is increasing and bounded (i.e., $\left\Vert A^{^{\prime }}\left(
u_{i}\right) \right\Vert _{L^{\infty }\left( \Omega \right) }\leq C\,,$
since $u_{i}$ is bounded), from (\ref{Np}) we have%
\begin{align}
I_{1}& =-\int_{\Omega }\left( A\left( u_{1}\right) -A\left( u_{2}\right)
\right) \left( u_{1}-u_{2}\right) dx+M_{12}\int_{\Omega }A\left(
u_{1}\right) -A\left( u_{2}\right) dx  \label{estI1} \\
& \leq C\left( M_{1}-M_{2}\right) ^{2}.  \notag
\end{align}%
To bound the $I_{2},I_{3}$ integral terms, we integrate by parts and proceed
as in \cite{BS, BRB}. We deduce%
\begin{eqnarray}
I_{2} &=&-\sum_{i,j}\int_{\Omega }\partial _{ij}\phi \partial _{j}\mathcal{K}%
\ast u_{1}\partial _{i}\phi dx-\sum_{i,j}\int_{\Omega }\partial _{i}\phi
\partial _{jj}\mathcal{K}\ast u_{1}\partial _{i}\phi dx  \label{estu1} \\
&&+\sum_{i,j}\int_{\Gamma }\left[ \left( \partial _{j}\mathcal{K}\ast
u_{1}\right) n_{j}\right] \left\vert \partial _{i}\phi \right\vert
^{2}d\sigma .  \notag
\end{eqnarray}%
The last term on the right-hand side of (\ref{estu1}) is nonpositive since $%
\Omega $ is convex and $\mathcal{K}$ is radially decreasing (i.e., as in 
\cite{BS}, we have $\left( \nabla \mathcal{K}\ast u_{1}\right) \cdot 
\overrightarrow{n}\leq 0$ on $\Gamma _{N}$). Integration by parts in the
first term gives the bound:%
\begin{equation*}
\sum_{i,j}\int_{\Omega }\partial _{ij}\phi \partial _{j}\mathcal{K}\ast
u_{1}\partial _{i}\phi dx\leq -\frac{1}{2}\int_{\Omega }\left( \Delta 
\mathcal{K}\ast u_{1}\right) \left\vert \nabla \phi \right\vert ^{2}dx,
\end{equation*}%
which together with (\ref{estu1}) entails%
\begin{equation}
I_{2}\leq C\int_{\Omega }\left\vert D^{2}\mathcal{K}\ast u_{1}\right\vert
\left\vert \nabla \phi \right\vert ^{2}dx,  \label{I2bound}
\end{equation}%
for some positive constant $C,$ which only depends implicitly on the
uniformly controlled $L^{p}$-norms of $u_{1}$ and $u_{2}$. Arguing as in 
\cite[(17)]{BRB}, by using Lemma \ref{kest}, (b) and since $\left(
M_{12}\right) ^{2}\leq \eta \left( t\right) ,$ we deduce for any $p\geq 2,$
that%
\begin{equation}
I_{2}\leq Cp\eta \left( t\right) ^{1-1/p}\text{ and }I_{3}\leq C\left\Vert
u_{2}\right\Vert _{L^{\infty }\left( \Omega \right) }\eta \left( t\right) .
\label{estu2}
\end{equation}%
Consequently, we obtain the differential inequality%
\begin{equation}
\frac{d}{dt}\eta \left( t\right) \leq Cp\eta \left( t\right) ^{1-1/p}+C\eta
\left( t\right) ,\text{ }\forall t\in \left[ 0,T\right] .  \label{decay1bd}
\end{equation}%
As in \cite[Theorem 3.3]{CVZ}, the idea is to fix $p=p\left( \eta \left(
t\right) \right) $ in an optimal way. To this end, choose a sufficiently
large constant $\widetilde{C}>0$ such that $p=\log (\widetilde{C}/\eta
\left( t\right) )\geq 2$, for all $t\in \left[ 0,T\right] .$ First, recall
that $\eta \left( t\right) $ is bounded on $\,t\in \left[ 0,T\right] $ and
that the first term on the right-hand side of (\ref{decay1bd}) dominates the
second one since we need the estimate for $\eta \left( t\right) $ small
only. Thus, the second term on the right-hand side is not essential and we
can derive the following differential inequality%
\begin{equation}
\frac{d}{dt}\eta \left( t\right) \leq C\eta \left( t\right) \log \left( 
\frac{\widetilde{C}}{\eta \left( t\right) }\right) ,\text{ }\forall t\in %
\left[ 0,T\right] .  \label{decay2b}
\end{equation}%
Here, we have used the elementary inequality $\eta \left( t\right)
^{-1/p}\leq C$, for the $p$-chosen above. Integrating (\ref{decay2b}) with
respect to $t\in \left( \delta ,s\right) $, we obtain%
\begin{equation}
\eta \left( s\right) \leq \widetilde{C}\left[ \frac{\eta \left( \delta
\right) }{\widetilde{C}}\right] ^{e^{-C\left( s-\delta \right) }},\text{ }%
\forall s\in (0,T].  \label{decay3b}
\end{equation}%
Passing to the limit as $\delta \rightarrow 0$ in (\ref{decay3b}), recalling
that $\eta \left( t\right) $ is continuous, we get the desired estimate (\ref%
{lip2}).
\end{proof}

\section{Optimal regularity and the global attractor}

\label{optim}

In this section, we derive several uniform estimates for the solutions of
the problem (\ref{be1})-(\ref{be2}) which are necessary for the study of the
asymptotic behavior as time goes to infinity. In a first step, we obtain
dissipative estimates for solutions in the spaces $L^{p},$ $L^{\infty }$ and 
$C^{\alpha },$ $\alpha >0$, uniformly with respect to time and the initial
data. Incidently, the estimates derived below allow one also to obtain \emph{%
optimal} regularity results for the weak solutions on $\Omega \times \lbrack
\tau ,\infty ),$ for every $\tau >0$, associated with the system (\ref{be1}%
)-(\ref{be3}). Finally, the apriori estimates will be deduced by a formal
argument, which can be justified rigourously by means of the approximation
procedure devised in \cite[Section 2]{BS}, \cite{BRB} by means of (\ref%
{approx}). Regardless of the type of approximation procedure being used the
regularity properties%
\begin{equation}
u_{\varepsilon }\in L^{\infty }\left( 0,T;L^{\infty }\left( \Omega \right)
\right) \cap L^{2}\left( 0,T;H_{D}^{1}\left( \Omega \right) \right) ,\quad
\forall \varepsilon ,T>0  \label{suf}
\end{equation}%
are \emph{essential} in order to rigorously perform these computations. To
this end, we shall only perform our (formal)\ computations to the original
system (\ref{be1})-(\ref{be3}), for the sake of simplicity.

The (uniform) dissipative $L^{p}$-estimate when $\sigma \left( \Gamma
_{D}\right) >0$ is different than the estimate when $\sigma \left( \Gamma
_{D}\right) =0$ (in this case, it was obtained in \cite[Lemma 10 and Remark 9%
]{BRB}), since in the former case there is no conservation of mass in (\ref%
{approx}). It is given by the following

\begin{proposition}
\label{lpbound}Let the assumptions (H1)-(H3) be satisfied and assume that $%
\sigma \left( \Gamma _{D}\right) >0$ (Case (ii)). Then there exists
constants $C_{\ast },C_{+}>0,$ and $\mu >1$, independent of time and the
initial data, such that every weak solution of (\ref{be1})-(\ref{be3})
satisfies%
\begin{equation}
\left\Vert u\left( t\right) \right\Vert _{L^{p}\left( \Omega \right)
}^{p}\leq \left( \frac{C_{\ast }}{C_{+}}\right) ^{\frac{1}{\mu }}+\left[
C_{+}\left( \mu -1\right) t\right] ^{-\frac{1}{\mu -1}},\text{ for all }%
t\geq \tau >0,  \label{lpest}
\end{equation}%
provided that $p>m+2d/\left( d-1\right) $. The constants $C_{\ast },C_{+}$
and $\mu $ can be computed explicitly in terms of the physical parameters of
the problem.
\end{proposition}

\begin{proof}
We first begin by noting that $\nabla \mathcal{K}\in L^{q,\infty }\left(
\Omega \right) ,$ for any $q>d/\left( d-1\right) ,$ by Lemma \ref{kest},
(a). This yields%
\begin{equation}
\nabla \mathcal{K}\mathbf{1}_{B_{1}\left( 0\right) }\in L^{1}\left( \mathbb{R%
}^{d}\right) \text{ and }\nabla \mathcal{K}\mathbf{1}_{\mathbb{R}%
^{d}\backslash B_{1}\left( 0\right) }\in L^{q}\left( \mathbb{R}^{d}\right) ,
\label{propk}
\end{equation}%
for any $q>d/\left( d-1\right) $. Next, let $\varphi \in C^{1}\left( 
\overline{\Omega }\right) \cap C^{2}\left( \Omega \right) $ and $\lambda $
denote the principal eigenfunction and eigenvalue of 
\begin{equation}
-\Delta \varphi =\lambda \varphi \text{ in }\Omega ,  \label{sp1}
\end{equation}%
with the boundary condition%
\begin{equation}
\nabla \varphi \cdot \overrightarrow{n}=0\text{ on }\Gamma _{N},  \label{sp2}
\end{equation}%
such that%
\begin{equation}
\nabla \varphi \cdot \overrightarrow{n}+\varphi =0\text{ on }\Gamma _{D}.
\label{sp3}
\end{equation}%
Clearly, $\lambda >0$ and by the maximum principle (see, e.g., \cite{Ev}) $%
\varphi >0$ in $\overline{\Omega }$ since $\sigma \left( \Gamma _{D}\right)
>0$. W.l.o.g., we may assume that $\left\Vert \varphi \right\Vert
_{L^{1}\left( \Omega \right) }=1.$

Set now $a\left( y\right) :=A^{^{\prime }}\left( y\right) .$ Testing
equation (\ref{identity}) by $pu^{p-1}\varphi ,$ $p>1,$ we obtain%
\begin{align}
\frac{d}{dt}\int_{\Omega }u^{p}\left( t\right) \varphi dx& =-p\int_{\Omega
}a\left( u\left( t\right) \right) \nabla u\left( t\right) \cdot \nabla
\left( u^{p-1}\left( t\right) \varphi \right) dx  \label{est1} \\
& +p\int_{\Omega }u\left( t\right) \left( \nabla \mathcal{K}\ast u\left(
t\right) \right) \nabla \left( u^{p-1}\left( t\right) \varphi \right) dx 
\notag \\
& =I_{1}+I_{2}+I_{3}+I_{4},  \notag
\end{align}%
where%
\begin{equation*}
\begin{array}{ll}
I_{1}:=-p\left( p-1\right) \int_{\Omega }a\left( u\right) u^{p-2}\left\vert
\nabla u\right\vert ^{2}\varphi \left( x\right) dx, &  \\ 
I_{2}:=-p\int_{\Omega }a\left( u\right) u^{p-1}\nabla u\cdot \nabla \varphi
\left( x\right) dx, &  \\ 
I_{3}:=p\left( p-1\right) \int_{\Omega }u\left( \nabla \mathcal{K}\ast
u\right) \cdot \nabla u\left( u^{p-2}\varphi \left( x\right) \right) dx, & 
\\ 
I_{4}:=p\int_{\Omega }u\left( \nabla \mathcal{K}\ast u\right) \cdot \nabla
\varphi \left( x\right) u^{p-1}dx. & 
\end{array}%
\end{equation*}%
First, from (H2) it is easy to see that%
\begin{equation}
I_{1}\leq -Cp\left( p-1\right) \int_{\Omega }u^{p-3+m}\left\vert \nabla
u\right\vert ^{2}\varphi \left( x\right) dx\leq 0.  \label{est1n}
\end{equation}%
On the other hand, setting $\widetilde{a}\left( u\right)
=\int_{0}^{u}a\left( y\right) y^{p-1}dy\geq \frac{C_{A}}{m+p-1}u^{m+p-1}$,
and using the definition of $\varphi $ from (\ref{sp1})-(\ref{sp3}), we have%
\begin{eqnarray}
I_{2} &=&-p\int_{\Omega }\nabla \widetilde{a}\left( u\right) \cdot \nabla
\varphi \left( x\right) dx=p\int_{\Omega }\widetilde{a}\left( u\right)
\Delta \varphi \left( x\right) dx  \label{est2nn} \\
&&-p\int_{\Gamma _{D}}\widetilde{a}\left( u\right) \nabla \varphi \cdot 
\overrightarrow{n}d\sigma -p\int_{\Gamma _{N}}\widetilde{a}\left( u\right)
\nabla \varphi \cdot \overrightarrow{n}d\sigma  \notag \\
&=&-p\lambda \int_{\Omega }\widetilde{a}\left( u\right) \varphi \left(
x\right) dx  \notag \\
&\leq &-C\lambda \int_{\Omega }u^{m+p-1}\varphi dx.  \notag
\end{eqnarray}%
since $\widetilde{a}\left( 0\right) =0$ ($u=0,$ a.e. on $\Gamma _{D}\times
\left( 0,\infty \right) $) and $\varphi $ satisfies (\ref{sp2}). Moreover,
we can estimate the integrals $I_{1}$, $I_{2}$ using H\"{o}lder and Young
inequalities as follows:%
\begin{align}
I_{3}& =p\left( p-1\right) \int_{\Omega }\left( u^{\frac{p-m+1}{2}}\sqrt{%
\varphi }\left( \nabla \mathcal{K}\ast u\right) \right) \cdot \left( \sqrt{%
\varphi }u^{\frac{p-3+m}{2}}\nabla u\right) dx  \label{est3n} \\
& \leq p\left( p-1\right) \left( \int_{\Omega }u^{p-3+m}\left\vert \nabla
u\right\vert ^{2}\varphi \left( x\right) dx\right) ^{1/2}\left( \int_{\Omega
}\left\vert \nabla \mathcal{K}\ast u\right\vert ^{2}u^{p-m+1}\varphi \left(
x\right) dx\right) ^{1/2}  \notag \\
& \leq \eta p\left( p-1\right) I_{1}+C_{\eta }\int_{\Omega }\left\vert
\nabla \mathcal{K}\ast u\right\vert ^{2}u^{p-m+1}\varphi \left( x\right) dx,
\notag
\end{align}%
for every $\eta >0.$ We must once again absorb the last term on the
right-hand side of (\ref{est3n}) into $I_{1}$. For any $q>d/\left(
d-1\right) $ and $s>1$ such that $1/s+1/q=1,$ Holder's inequality yields%
\begin{align}
C_{\eta }\int_{\Omega }\left\vert \nabla \mathcal{K}\ast u\right\vert
^{2}u^{p-m+1}dx& \leq C_{\eta }\left\Vert \nabla \mathcal{K}\ast
u\right\Vert _{L^{2q}\left( \Omega \right) }^{2}\left( \int_{\Omega
}u^{\left( p-m+1\right) s}dx\right) ^{1/s}  \label{est4n} \\
& \leq C_{\eta }\left( \mathcal{K}\right) \left\Vert u\right\Vert
_{L^{2q}\left( \Omega \right) }^{2}\left\Vert u\right\Vert _{L^{\left(
p-m+1\right) s}\left( \Omega \right) }^{p-m+1},  \notag
\end{align}%
since, from (\ref{propk}), there holds%
\begin{align}
\left\Vert \nabla \mathcal{K}\ast u\right\Vert _{L^{2q}\left( \Omega \right)
}& \leq \left\Vert \nabla \mathcal{K}\mathbf{1}_{B_{1}\left( 0\right)
}\right\Vert _{L^{1}}\left\Vert u\right\Vert _{L^{2q}\left( \Omega \right)
}+\left\Vert \nabla \mathcal{K}\mathbf{1}_{\mathbb{R}^{d}\backslash
B_{1}\left( 0\right) }\right\Vert _{L^{2q}}\left\Vert u\right\Vert
_{L^{1}\left( \Omega \right) }  \label{propkest} \\
& \leq C\left( \mathcal{K},\Omega \right) \left\Vert u\right\Vert
_{L^{2q}\left( \Omega \right) }.  \notag
\end{align}%
Thus, choosing $q,s>1$ in (\ref{est4n}) in an optimal way such that $%
2q=\left( p-m+1\right) s>2d/\left( d-1\right) ,$ we further obtain in (\ref%
{est3n}) by virtue of (\ref{est4n}) and Young's inequality, that%
\begin{align}
I_{3}& \leq \eta p\left( p-1\right) I_{1}+C_{\eta }\left( \mathcal{K},\Omega
\right) \left\Vert u\right\Vert _{L^{p-m+3}\left( \Omega \right) }^{p-m+3}
\label{est5n} \\
& \leq \eta p\left( p-1\right) I_{1}+\overline{\eta }\int_{\Omega
}u^{m-1+p}\varphi dx+C_{\eta ,\overline{\eta }}\left( \mathcal{K},\Omega 
\text{,}p,m,\varphi \right) ,  \notag
\end{align}%
for every $\overline{\eta }>0,$ since $p+m-1>p-m+3$ (recall that $m>2$). It
follows by choosing sufficiently small $\eta \leq C/2$ in (\ref{est3n}), and 
$\overline{\eta }\leq (C/2)\lambda $ in (\ref{est5n}) that the integral term 
$I_{3}$ can be completely absorbed into $I_{1}$. The term $I_{4}$ can be
bounded exactly the same way. We have%
\begin{equation*}
I_{4}=p\int_{\Omega }u\left( \nabla \mathcal{K}\ast u\right) \cdot \nabla
\varphi \left( x\right) u^{p-1}dx\leq \epsilon \int_{\Omega
}u^{p+m-1}\varphi \left( x\right) dx+C_{\epsilon }\left( \mathcal{K},\Omega 
\text{,}p,m,\varphi \right) ,
\end{equation*}%
for every $\epsilon >0$ and some positive constant $C_{\epsilon }$ which
depends on $\varphi \in C^{1}\left( \overline{\Omega }\right) .$ Summing up,
from (\ref{est1}) we deduce 
\begin{equation}
\frac{d}{dt}\int_{\Omega }u^{p}\left( t\right) \varphi dx+C\lambda
\int_{\Omega }u^{p+m-1}\left( t\right) \varphi \left( x\right) dx\leq
C\left( \mathcal{K},\Omega \text{,}p,m,\varphi \right) .  \label{est5}
\end{equation}%
Finally, set $\mu :=\left( m-1+p\right) /p>1$ and 
\begin{equation*}
Z\left( t\right) :=\int_{\Omega }u^{p}\left( t\right) \varphi dx.
\end{equation*}%
By Jensen's inequality, (\ref{est5}) yields the following inequality:%
\begin{equation}
\frac{d}{dt}Z\left( t\right) +C_{+}\left( \lambda \right) \left( Z\left(
t\right) \right) ^{\mu }\leq C_{\ast }.  \label{p9}
\end{equation}%
We can now use the Gronwall's inequality (see, e.g., \cite[Chapter III,
Lemma 5.1]{T}), applied to (\ref{p9}) to deduce the desired claim in (\ref%
{lpest}). The proof of the proposition is complete.
\end{proof}

\begin{remark}
\label{rem}In the case when $u_{0}\in L^{p}\left( \Omega \right) $, $Z\left(
0\right) =\lim_{t\rightarrow 0^{+}}Z\left( t\right) $ is finite, so a
similar argument to \cite[Chapter III, Lemma 5.1]{T} gives%
\begin{equation}
Z\left( t\right) \leq \max \left\{ \left\Vert u_{0}\right\Vert _{L^{p}\left(
\Omega \right) },\left( \frac{C_{\ast }}{C_{+}}\right) ^{\frac{1}{\mu }%
}\right\} ,\text{ }\forall t\geq 0.  \label{lp2}
\end{equation}%
It is worth emphasizing again that the left-hand side of the inequality in
(H2) (i.e., $C_{A}y^{m-1}\leq A^{^{\prime }}\left( y\right) $, $\forall
y\geq 0$) is enough to establish the above assertion. Moreover, the above
estimate (\ref{lp2}) directly implies global well-posedness in the
subcritical case for problem (\ref{be1})-(\ref{be3}) for as long as $\sigma
\left( \Gamma _{D}\right) >0$ (Case (ii))$.$ Unfortunately, we are not able
to argue as in the proof of \cite[Lemma 10]{BRB} to obtain the desired $%
L^{p} $-estimate (similar to the case when $\sigma \left( \Gamma _{D}\right)
=0$, see (\ref{est8})\ below) since we do not know how to get apriori
control over the $L^{1}$-norm of $u\left( t\right) $ for problem (\ref%
{approx}). We emphasize again that when $\sigma \left( \Gamma _{D}\right) =0$%
, mass is conserved in both (\ref{approx}) and (\ref{be1})-(\ref{be3}) so
that $\left\Vert u\left( t\right) \right\Vert _{L^{1}\left( \Omega \right)
}\leq M\left\vert \Omega \right\vert $, for any $M>0$ such that $%
\left\langle u_{0}\right\rangle \leq M.$
\end{remark}

Next, we establish a crucial result which allows to deduce a dissipative $%
L^{\infty }$-estimate uniform with respect to the initial data, and which is
necessary for the attractor theory. For this result, we consider both cases $%
\sigma \left( \Gamma _{D}\right) \geq 0.$

\begin{lemma}
\label{bounded}Assume that there exists $\epsilon >0$ and $y_{\epsilon }>0$
such that $A^{^{\prime }}\left( y\right) \geq \epsilon $ for all $y\geq
y_{\epsilon }.$ Let $v$ be a solution of the following degenerate problem:%
\begin{equation}
\partial _{t}v\left( t\right) +\nabla \cdot (v\left( t\right) 
\overrightarrow{W}\left( t\right) )=\Delta A\left( v\left( t\right) \right) ,%
\text{ in }\Omega \times \left( 0,\infty \right) ,  \label{be4}
\end{equation}%
such that $v$ satisfies (\ref{be3}) and 
\begin{equation}
\left\{ 
\begin{array}{ll}
(\nabla A\left( v\right) -v\overrightarrow{W})\cdot \overrightarrow{n}=0%
\text{,} & \text{on }\Gamma _{N}\times \left( 0,\infty \right) , \\ 
v=0, & \text{on }\Gamma _{D}\times \left( 0,\infty \right) .%
\end{array}%
\right.  \label{be5}
\end{equation}%
If 
\begin{equation*}
\left\Vert \overrightarrow{W}\left( t\right) \right\Vert _{\left( L^{\infty
}\left( \Omega \right) \right) ^{d}}\leq \overline{C},\text{ for all }t>0,
\end{equation*}%
then there exists $\mu \sim (\tau ^{^{\prime }}-\tau )>0$, and a positive
constant $\mathcal{C}=\mathcal{C}\left( \overline{C},\mu ,\epsilon
,y_{\epsilon }\right) $ such that%
\begin{equation}
\sup_{t\geq \tau ^{^{\prime }}>0}\left\Vert v\left( t\right) \right\Vert
_{L^{\infty }\left( \Omega \right) }\leq \mathcal{C}Q(1+\sup_{t\geq \tau
>0}\left\Vert v\left( t\right) \right\Vert _{L^{1}\left( \Omega \right) }),
\label{boundlinf}
\end{equation}%
for all $\tau ^{^{\prime }}>\tau >0.$ The constant $\mathcal{C=C}\left( \mu
\right) \sim \mu ^{-l}$ ($l>0$) is bounded if $\mu $ is bounded away from
zero.
\end{lemma}

\begin{remark}
Lemma \ref{bounded} implies the boundedness of the function $v\left(
t\right) $ provided that $\overrightarrow{W}\left( t\right) $ is a bounded
vector on $[\tau ,\infty )$. It establishes the $L^{1}$-$L^{\infty }$
smoothing property for the solutions of (\ref{be4})-(\ref{be5}). This result
is analogous to the result obtained by Kowalczyk \cite{Kbis}, but it is more
sharp since the function and constant on the right-hand side of (\ref%
{boundlinf}) do \emph{not} depend on the $L^{\infty }$-norm of the initial
data (see, \cite[Lemma 4.1]{Kbis}). This is very useful if we want to
produce uniform estimates with respect to time and the initial data.
\end{remark}

\begin{proof}
\emph{Step 1 }(The local relation). First, we recall the following estimate
which can be obtained exactly as in \cite[Lemma 4.1, (5.1)-(5.4)]{Kbis}.
Indeed, setting $v_{m}=\left( \left\vert v\right\vert -l\right) _{+}$, for
any $l\geq y_{\epsilon }$, multiplying equation (\ref{be4}) by the $p$-th
power of $v_{l},$ $p>1$, and integrating by parts using (\ref{be2}), it
follows, after a suitable rescaling of the time variable $\overline{t}%
=\epsilon t$, that%
\begin{equation}
\frac{d}{d\overline{t}}\int_{\Omega }\overline{v}_{l}^{p+1}\left( \overline{t%
}\right) dx+\frac{2p}{p+1}\int_{\Omega }\left\vert \nabla \overline{v}_{l}^{%
\frac{p+1}{2}}\left( \overline{t}\right) \right\vert ^{2}dx\leq Cp\left(
p+1\right) \left( \int_{\Omega }\overline{v}_{l}^{p+1}\left( \overline{t}%
\right) dx+1\right) ,  \label{apx1}
\end{equation}%
for some constant $C>0$ independent of $p$, but which depends on $l$ and $%
\epsilon ,$ and where $\overline{v}_{l}\left( \overline{t},x\right)
=v_{l}\left( \overline{t}/\epsilon ,x\right) ,$ $\overline{t}\in \left(
0,\epsilon t\right) .$ In fact, $C\sim 1/\epsilon $ as $\epsilon \rightarrow
0.$

Secondly, set $p_{k}=2^{k}-1,$ $k\geq 0$, and define%
\begin{equation}
\mathcal{Y}_{k}\left( \overline{t}\right) :=\int_{\Omega }\overline{v}%
_{l}^{p_{k}+1}\left( \overline{t}\right) dx,  \label{def}
\end{equation}%
for all $k\geq 0$. Let $\overline{t},\mu $ be two positive constants such
that $\overline{t}-\mu /p_{k}>0$, and whose values will be chosen later. We
claim that there holds%
\begin{equation}
\mathcal{Y}_{k}\left( \overline{t}\right) \leq N_{k}\left( \overline{t},\mu
\right) :=C\left( \mu ,\epsilon \right) \left( p_{k}\right) ^{\gamma
}(\sup_{s\geq \overline{t}-\mu /p_{k}}\mathcal{Y}_{k-1}\left( s\right)
+1)^{2},\text{ }\forall k\geq 1,  \label{claim2}
\end{equation}%
where $C,$ $\gamma $ are positive constants independent of $k.$ The constant 
$C=C\left( \mu ,\epsilon \right) $ is bounded if $\mu $ is bounded away from
zero.

We will now prove (\ref{claim2}) when $2<d.$ The case $d\leq 2$ requires
only minor modifications by using a suitable Sobolev embedding. The argument
we follow is similar, for instance, to \cite{GM} (cf. also \cite{G0}). For
each $k\geq 0$, we define%
\begin{equation*}
r_{k}:=\frac{d\left( p_{k}+1\right) -\left( d-2\right) \left( 1+p_{k}\right) 
}{d\left( p_{k}+1\right) -\left( d-2\right) \left( 1+p_{k-1}\right) },\text{ 
}s_{k}:=1-r_{k}.
\end{equation*}%
We aim to estimate the term on the right-hand side of (\ref{apx1}) in terms
of the $L^{1+p_{k-1}}\left( \Omega \right) $-norm of $\overline{v}_{m}.$
First, H\"{o}lder and Sobolev inequalities (with the equivalent norm of
Sobolev spaces in $H^{1}\left( \Omega \right) \subset L^{p_{s}}\left( \Omega
\right) $, $p_{s}=2d/\left( d-2\right) $) yield%
\begin{align}
\int_{\Omega }\overline{v}_{l}^{p_{k}+1}dx& \leq \left( \int_{\Omega }\left( 
\overline{v}_{l}\right) ^{\frac{\left( p_{k}+1\right) d}{d-2}}dx\right)
^{s_{k}}\left( \int_{\Omega }\left( \overline{v}_{l}\right)
^{1+p_{k-1}}dx\right) ^{r_{k}}  \label{ee5b} \\
& \leq C\left( \int_{\Omega }\left\vert \nabla \overline{v}_{l}^{\frac{%
\left( p_{k}+1\right) }{2}}\right\vert ^{2}dx+\int_{\Omega }\left( \overline{%
v}_{l}\right) ^{1+p_{k}}dx\right) ^{\overline{s}_{k}}  \notag \\
& \times \left( \int_{\Omega }\left( \overline{v}_{l}\right)
^{1+p_{k-1}}dx\right) ^{r_{k}},  \notag
\end{align}%
with $\overline{s}_{k}:=s_{k}d/\left( d-2\right) \in \left( 0,1\right) $.
Applying Young's inequality on the right-hand side of (\ref{ee5b}), we get%
\begin{equation}
\int_{\Omega }\overline{v}_{l}^{1+p_{k}}dx\leq \frac{1}{4}\int_{\Omega
}\left\vert \nabla \overline{v}_{l}^{\frac{p_{k}+1}{2}}\right\vert
^{2}dx+Q_{\gamma _{1}}\left( p_{k}\right) \left( \int_{\Omega }\overline{v}%
_{l}^{p_{k-1}+1}dx\right) ^{\theta _{k}},  \label{bulk}
\end{equation}%
for some positive constant $\gamma _{1}$ independent of $p_{k},$ and where%
\begin{equation*}
\theta _{k}:=\frac{r_{k}}{1-\overline{s}_{k}}\geq 1
\end{equation*}%
(in fact, straightforward computations give $\theta _{k}\equiv 2$ since $%
p_{k}=2^{k}-1$, for all $k$).\ Note that $Q_{\gamma _{1}}\left( s\right)
\sim s^{\gamma _{1}}$ as $s\rightarrow \infty $. Inserting estimate (\ref%
{bulk}) on the right-hand side of (\ref{apx1}), we obtain the following
inequality:%
\begin{equation}
\frac{d}{d\overline{t}}\mathcal{Y}_{k}\left( \overline{t}\right) +\frac{3}{4}%
\int_{\Omega }\left\vert \nabla \overline{v}_{l}^{\frac{p_{k}+1}{2}%
}\right\vert ^{2}dx\leq C\left( p_{k}\right) ^{\gamma }\left( \mathcal{Y}%
_{k-1}\left( \overline{t}\right) +1\right) ^{2},  \label{claim}
\end{equation}%
where $C,$ $\gamma $ are positive constants independent of $k.$

We are now ready to prove (\ref{claim2}) using (\ref{claim}). To this end,
let $\zeta \left( s\right) $ be a positive function $\zeta :\mathbb{R}%
_{+}\rightarrow \left[ 0,1\right] $ such that $\zeta \left( s\right) =0$ for 
$s\in \left[ 0,\overline{t}-\mu /p_{k}\right] ,$ $\zeta \left( s\right) =1$
if $s\in \left[ \overline{t},+\infty \right) $ and $\left\vert d\zeta
/ds\right\vert \leq p_{k}/\mu $, if $s\in \left( \overline{t}-\mu /p_{k},%
\overline{t}\right) $. We define $Z_{k}\left( s\right) =\zeta \left(
s\right) \mathcal{Y}_{k}\left( s\right) $ and notice that%
\begin{equation*}
\frac{d}{ds}Z_{k}\left( s\right) \leq \frac{p_{k}}{\mu }\mathcal{Y}%
_{k}\left( s\right) +\zeta \left( s\right) \frac{d}{ds}\mathcal{Y}_{k}\left(
s\right) .
\end{equation*}%
Combining this estimate with (\ref{claim}), (\ref{bulk}) and noticing that $%
Z_{k}\leq \mathcal{Y}_{k}$, we deduce the following estimate for $Z_{k}$:%
\begin{equation}
\frac{d}{ds}Z_{k}\left( s\right) +C\left( \mu \right) p_{k}Z_{k}\left(
s\right) \leq N_{k}\left( \overline{t},\mu \right) ,\text{ for all }s\in %
\left[ \overline{t}-\mu /p_{k},+\infty \right) ,  \label{e12}
\end{equation}%
for some positive constant $C=C\left( \mu \right) $ independent of $k$.
Integrating (\ref{e12}) with respect to $s$ from $\overline{t}-\mu /p_{k}$
to $\overline{t}$ and taking into account the fact that $Z_{k}\left( 
\overline{t}-\mu /p_{k}\right) =0,$ we obtain that%
\begin{equation*}
\mathcal{Y}_{k}\left( \overline{t}\right) =Z_{k}\left( \overline{t}\right)
\leq N_{k}\left( \overline{t},\mu \right) \left( 1-e^{-C\mu }\right) ,
\end{equation*}%
which proves the claim (\ref{claim2}).

\emph{Step 2} (The iteration procedure). Let now $\overline{\tau }^{^{\prime
}}>\overline{\tau }>0$ be given with $\overline{\tau }>0$ such that%
\begin{equation}
\sup_{s\geq \overline{t}_{1}=\overline{\tau }}\left( \mathcal{Y}_{0}\left(
s\right) +1\right) \leq \sup_{s\geq \overline{t}_{1}=\overline{\tau }}\left(
\left\Vert \overline{v}\left( s\right) \right\Vert _{L^{1}\left( \Omega
\right) }+1+l\left\vert \Omega \right\vert \right) \leq C_{\blacklozenge
}\left( l\right) ,  \label{l1}
\end{equation}%
and define $\mu =(\overline{\tau }^{^{\prime }}-\overline{\tau }),$ $%
\overline{t}_{0}=\overline{\tau }^{^{\prime }}$ and $\overline{t}_{k}=%
\overline{t}_{k-1}-\mu /p_{k},$ $k\geq 1$. Using (\ref{claim}), we have%
\begin{equation}
\sup_{\overline{t}\geq \overline{t}_{k-1}}\mathcal{Y}_{k}\left( \overline{t}%
\right) \leq C\left( p_{k}\right) ^{\gamma }(\sup_{s\geq \overline{t}_{k}}%
\mathcal{Y}_{k-1}\left( s\right) +1)^{2},\text{ }k\geq 1.  \label{e13}
\end{equation}%
We can iterate in (\ref{e13}) with respect to $k\geq 1$ and obtain that%
\begin{align}
\sup_{\overline{t}\geq \overline{t}_{k-1}}\mathcal{Y}_{k}\left( \overline{t}%
\right) & \leq \left( Cp_{k}^{\gamma }\right) \left( Cp_{k-1}^{\gamma
}\right) ^{2}\left( Cp_{k-2}^{\gamma }\right) ^{2^{2}}\cdot ...\cdot \left(
Cp_{1}^{\gamma }\right) ^{2^{k}}(C_{\blacklozenge })^{2^{k}}  \label{e14} \\
& \leq C^{A_{k}}2^{B_{k}\gamma }\left( C_{\blacklozenge }\right) ^{2^{k}}, 
\notag
\end{align}%
where%
\begin{equation}
A_{k}:=1+2+2^{2}+...+2^{k},  \label{ak2}
\end{equation}%
\begin{equation}
B_{k}:=k+2\left( k-1\right) +2^{2}\left( k-2\right) +...+2^{k}.  \label{bk2}
\end{equation}%
We can easily show that $\left\{ A_{k}\right\} $ and $\left\{ B_{k}\right\} $
satisfy 
\begin{equation}
A_{k}\leq C\left( 2^{k}\right) \text{ and }B_{k}\leq C\left( 2^{k}\right) ,
\label{abk}
\end{equation}%
for some positive constant $C$ independent of $k,$ $j$ and $\mu $ (see,
e.g., \cite{GM}). We can take the $1+p_{k}=2^{k}$-root on both sides of (\ref%
{e14}) and let $k\rightarrow +\infty $. We deduce%
\begin{equation}
\sup_{\overline{t}\geq \overline{t}_{0}=\overline{\tau }^{^{\prime
}}}\left\Vert \overline{v}_{l}\left( \overline{t}\right) \right\Vert
_{L^{\infty }\left( \Omega \right) }\leq \lim_{k\rightarrow +\infty }\sup_{%
\overline{t}\geq \overline{t}_{0}}\left( \mathcal{Y}_{k}\left( \overline{t}%
\right) \right) ^{1/\left( 1+p_{k}\right) }\leq C\left( \mu ,\epsilon
,y_{\epsilon }\right) \left( C_{\blacklozenge }\right) ,  \label{linf}
\end{equation}%
for some positive constant $C$ independent of $\overline{t},$ $k$, $%
\overline{v}_{l}$ and initial data. Next, we notice that, for any $p>1,$ in
view of $\sqrt[p]{a+b}\leq \sqrt[p]{a}+\sqrt[p]{b},$ there holds%
\begin{equation*}
\left\Vert \overline{v}\left( \overline{t}\right) \right\Vert _{L^{p}\left(
\Omega \right) }\leq 2^{\frac{p-1}{p}}\left( \left\Vert \overline{v}%
_{l}\left( \overline{t}\right) \right\Vert _{L^{p}\left( \Omega \right)
}+l\left\vert \Omega \right\vert ^{1/p}\right) ;
\end{equation*}%
thus, as $p\rightarrow \infty $, we have from (\ref{linf}),%
\begin{equation}
\sup_{\overline{t}\geq \overline{\tau }^{^{\prime }}}\left\Vert \overline{v}%
\left( \overline{t}\right) \right\Vert _{L^{\infty }\left( \Omega \right)
}\leq 2\left( C\left( \mu ,\epsilon ,y_{\epsilon }\right) \left(
C_{\blacklozenge }\right) +l\right) ,  \label{linffinal}
\end{equation}%
for any $l\geq n_{\epsilon }.$ Rescaling back the time variables ($\overline{%
t}=t\epsilon ,$ $\overline{\tau }^{^{\prime }}=\tau ^{^{\prime }}\epsilon ,$ 
$\overline{\tau }=\tau \epsilon $) into (\ref{linffinal}) and (\ref{l1}) and
taking $l=y_{\epsilon }$, we easily obtain the desired inequality (\ref%
{boundlinf}). The proof is finished.
\end{proof}

We can now show the following.

\begin{theorem}
\label{holderthm}Let $0\leq u_{0}\in L^{\infty }\left( \Omega \right) $ and
assume (H1)-(H3) hold for both cases (i)+(ii). If $\sigma \left( \Gamma
_{D}\right) =0$, let $M\geq 0$ be given such that $\left\langle
u_{0}\right\rangle \leq M$. Every weak solution of (\ref{be1})-(\ref{be3})
is globally H\"{o}lder continuous in the cylinders $[\tau ,\tau ^{^{\prime
}}]\times \overline{\Omega }$, for all $\tau ^{^{\prime }}>\tau >0$. In
particular, the following estimate holds:%
\begin{equation}
\sup_{t\geq 1}\left( \left\Vert A\left( u\right) \right\Vert _{L^{2}\left( 
\left[ t,t+1\right] \times \Omega \right) }+\left\Vert u\right\Vert
_{C^{\alpha /2,\alpha }\left( \left[ t,t+1\right] \times \overline{\Omega }%
\right) }\right) \leq C,\text{ if }\sigma \left( \Gamma _{D}\right) >0.
\label{holderbound}
\end{equation}%
Finally, for every bounded subset $B=B\left( M\right) \subset L^{1}\left(
\Omega \right) $, there exists a time $t_{+}=t_{+}\left( B\right) >0$ such
that%
\begin{equation}
\sup_{t\geq t_{+}}\left( \left\Vert A\left( u\right) \right\Vert
_{L^{2}\left( \left[ t,t+1\right] \times \Omega \right) }+\left\Vert
u\right\Vert _{C^{\alpha /2,\alpha }\left( \left[ t,t+1\right] \times 
\overline{\Omega }\right) }\right) \leq C,\text{ if }\sigma \left( \Gamma
_{D}\right) =0,  \label{holderbound2}
\end{equation}%
for some $\alpha >0$ and some constant $C>0$ independent of the initial data
and time.
\end{theorem}

\begin{proof}
It suffices to prove the estimates (\ref{holderbound})-(\ref{holderbound2}).

\emph{Case 1.} Consider first $\sigma \left( \Gamma _{D}\right) >0$. In this
case, we can take $p>d$ and $p>m+2d/\left( d-1\right) $ as large as we want
in Proposition \ref{lpbound} so that%
\begin{equation}
\sup_{t\geq \frac{1}{2}}\left\Vert \nabla \mathcal{K}\ast u\left( t\right)
\right\Vert _{L^{p}\left( \Omega \right) }\leq C  \label{est6}
\end{equation}%
since $\mathcal{K}\in L^{d/\left( d-1\right) ,\infty }\left( \mathbb{R}%
^{d}\right) $, see (\ref{propk}), (\ref{propkest}). Thus, exploiting first
Lemma \ref{kest}-(b), we have for $\overrightarrow{V}=\nabla \mathcal{K}\ast
u,$%
\begin{equation}
\sup_{t\geq \frac{1}{2}}\left\Vert \nabla \overrightarrow{V}\right\Vert
_{L^{p}\left( \Omega \right) }\leq C\left( p\right) \sup_{t\geq \frac{1}{2}%
}\left\Vert u\left( t\right) \right\Vert _{L^{p}\left( \Omega \right) }\leq
C.  \label{est6bis}
\end{equation}%
By Morrey's inequality (see, e.g., \cite{Ev}),%
\begin{equation}
\sup_{t\geq \frac{1}{2}}\left\Vert \overrightarrow{V}\right\Vert _{L^{\infty
}\left( \Omega \right) }\leq C.  \label{est6tris}
\end{equation}%
Invoking now the crucial Lemma \ref{bounded}, we deduce on account of (\ref%
{lpest}), that%
\begin{equation}
\sup_{t\geq 1}\left\Vert u\left( t\right) \right\Vert _{L^{\infty }\left(
\Omega \right) }\leq C.  \label{est7}
\end{equation}%
We can now test equation (\ref{identity}) with $A\left( u\left( t\right)
\right) $, then integrate over $\left( t,t+1\right) $. Setting $\overline{A}%
\left( u\right) :=\int_{0}^{u}A\left( s\right) ds,$ and exploiting (\ref%
{est7}) we obtain for every $t\geq 1,$%
\begin{align*}
& \int_{\Omega }\left( \overline{A}\left( u\left( t+1\right) \right) -%
\overline{A}\left( u\left( t\right) \right) \right)
dx+\int_{t}^{t+1}\left\Vert \nabla A\left( u\left( s\right) \right)
\right\Vert _{L^{2}\left( \Omega \right) }^{2}ds \\
& =\int_{t}^{t+1}\int_{\Omega }\nabla A\left( u\left( s\right) \right) \cdot
\left( u\nabla \mathcal{K}\ast u\right) dxds \\
& \leq \frac{1}{2}\int_{t}^{t+1}\left\Vert \nabla A\left( u\left( s\right)
\right) \right\Vert _{L^{2}\left( \Omega \right) }^{2}ds+C.
\end{align*}%
Thus, we easily deduce the first part of (\ref{holderbound}). For every $%
1<p<\infty ,$ we now have%
\begin{equation}
\Delta \mathcal{K}\ast u\in L^{p}\left( \Omega \times \lbrack 1,\infty
)\right) \text{ and }\overrightarrow{V}=\nabla \mathcal{K}\ast u\in
L^{\infty }\left( \Omega \times \lbrack 1/2,\infty )\right)  \label{vel}
\end{equation}%
since $u$ is bounded. Global H\"{o}lder continuity results for (\ref{be1})-(%
\ref{be2}) have been proven for more general classes of degenerate
quasilinear equations in divergence form \cite{DF, Do, Iv, K, PV}. Due to (%
\ref{vel}), the second estimate in (\ref{holderbound}) (respectively, (\ref%
{holderbound2})) is a simple corollary of these results when $m\geq m_{\ast
} $, see, e.g., \cite{Du} (cf. also \cite{Chen, DF, PV}).

\emph{Case 2.} The case $\sigma \left( \Gamma _{D}\right) =0$ is similar.
Following \cite[Lemma 10]{BRB} and Lemma \ref{kest}, (c), we see that for
any $m\in \left[ 1,m_{\ast }\right] $ such that%
\begin{equation}
\underset{z\rightarrow \infty }{\lim \inf }A^{^{\prime }}\left( z\right)
z^{1-m}>0  \label{Acond}
\end{equation}%
and $\overline{q}=\left( 2-m\right) /\left( 2-m_{\ast }\right) \geq 1$, with 
$u_{k}:=\left( u-k\right) _{+}\in L^{\overline{q}}\left( \Omega \right) $,
the following estimate holds (cf. \cite[pg. 1703]{BRB}):%
\begin{equation}
\frac{d}{dt}\left\Vert u_{k}\left( t\right) \right\Vert _{L^{p}\left( \Omega
\right) }^{p}\leq -\eta \left\Vert u_{k}\left( t\right) \right\Vert
_{L^{p}\left( \Omega \right) }^{p}+C\left( M,\Omega ,k,p,\left\Vert
u_{k}\left( t\right) \right\Vert _{L^{\overline{q}}\left( \Omega \right)
}\right) ,  \label{est8}
\end{equation}%
for all $t\geq 0,$ for some positive constant $\eta >0$. It remains to note
that by assumption (H2), condition (\ref{Acond}) is already satisfied with $%
m=m_{\ast }$, so we can choose $\overline{q}=1$ in (\ref{est8}). Gronwall's
inequality (see, e.g., \cite{T}) applied to (\ref{est8}) yields, on account
of (\ref{ukest}), that%
\begin{equation}
\left\Vert u\left( t\right) \right\Vert _{L^{p}\left( \Omega \right)
}^{p}\leq \left\Vert \left( u_{0}-k\right) _{+}\right\Vert _{L^{p}\left(
\Omega \right) }^{p}e^{-\eta t}+\overline{C}\left( M,\Omega ,k,p\right) ,
\label{est9}
\end{equation}%
since $\left\Vert u\left( t\right) \right\Vert _{L^{1}\left( \Omega \right)
}=\left\vert \Omega \right\vert \left\langle u\left( t\right) \right\rangle
\leq M\left\vert \Omega \right\vert $ (mass is conserved). Hence, for every
subset $B=B\left( M\right) \subset L^{1}\left( \Omega \right) $, and setting 
$R_{\ast }=\overline{C}\left( M,\Omega ,k,p\right) +1,$ we can easily find a
time $t_{+}=t_{+}\left( M\right) >0$ such that (\ref{est9}) implies%
\begin{equation}
\sup_{t\geq t_{+}}\left\Vert u\left( t\right) \right\Vert _{L^{p}\left(
\Omega \right) }^{p}\leq R_{\ast }.  \label{est10}
\end{equation}%
After establishing the bound (\ref{est9}), the same estimates (\ref{est6})-(%
\ref{est6tris}) hold provided that $\nabla \mathcal{K}\in L^{q,\infty
}\left( \Omega \right) ,$ for $q>d/\left( d-1\right) $ is satisfied. Thus,
we can again reach the uniform $L^{\infty }$-estimate (\ref{est7}) now on
the time intervals $[t_{+},\infty ).$ This is enough to get (\ref%
{holderbound2}) once more by applying the results in \cite{Du, PV} (see 
\emph{Case 1} above) and to conclude the proof of the theorem.
\end{proof}

Define, for some given $M\geq 0$,%
\begin{equation*}
\mathcal{Z}_{DN}:=\left\{ 
\begin{array}{ll}
\left\{ u_{0}\in L^{\infty }\left( \Omega \right) :u_{0}\geq 0\text{, }%
\left\langle u_{0}\right\rangle \leq M\right\} , & \text{if }\sigma \left(
\Gamma _{D}\right) =0, \\ 
\left\{ u_{0}\in L^{\infty }\left( \Omega \right) :u_{0}\geq 0\right\} & 
\text{if }\sigma \left( \Gamma _{D}\right) >0.%
\end{array}%
\right.
\end{equation*}%
In the previous section we have proved that system (\ref{be1})-(\ref{be3})
generates a continuous semigroup $S(t)$ on the phase space $\mathcal{Z}%
_{DN}, $ endowed with the metric topology of $\left( H_{D}^{1}\left( \Omega
\right) \right) ^{\ast },$ via%
\begin{equation}
S\left( t\right) :\left[ \mathcal{Z}_{DN}\right] _{\ast }\rightarrow \left[ 
\mathcal{Z}_{DN}\right] _{\ast },\text{ }S\left( t\right) u_{0}=u\left(
t\right) ,\text{ }t\geq 0,  \label{dyn}
\end{equation}%
where $u(t)$ is a unique weak solution of (\ref{be1})-(\ref{be3}) (see
Theorem \ref{well-posed}). Here $\left[ \mathcal{Z}_{DN}\right] _{\ast }$
denotes the closure of $\mathcal{Z}_{DN}$ in the metric of $\left(
H^{1}\left( \Omega \right) \right) ^{\ast }.$ We devote next our attention
to the study of the long-time behavior of trajectories of the semigroup in
terms of global attractors. We need to recall the following definition.

\begin{definition}
\label{def_gl}We say that $\mathcal{A}\subset \mathcal{Z}_{DN}$ is the
global attractor for the dynamical system $\left( S\left( t\right) ,\mathcal{%
Z}_{DN}\right) $ if the following conditions are satisfied:

(i) the set $\mathcal{A}$ is a compact subset of the phase space $\mathcal{Z}%
_{DN}$;

(ii) it is strictly invariant, that is, $S(t)\mathcal{A}$ $=\mathcal{A}$,
for all $t\geq 0$;

(iii) for every bounded subset $B\subset \mathcal{Z}_{DN}$,%
\begin{equation}
dist_{\left[ \mathcal{Z}_{DN}\right] _{\ast }}\left( S\left( t\right) B,%
\mathcal{A}\right) \rightarrow 0\text{, as }t\rightarrow \infty ,
\label{gl_conv}
\end{equation}%
that is, $\mathcal{A}$ attracts the images of all bounded subsets of $%
\mathcal{Z}_{DN}$ as time goes to infinity. Here,%
\begin{equation*}
dist_{\left[ \mathcal{Z}_{DN}\right] _{\ast }}\left( X,Y\right) =\sup_{x\in
X}\inf_{y\in Y}\left\Vert x-y\right\Vert _{\left[ \mathcal{Z}_{DN}\right]
_{\ast }}.
\end{equation*}
\end{definition}

The first main result of this section states the existence of such an
attractor for problem (\ref{be1})-(\ref{be3}).

\begin{theorem}
\label{gl_att}Let the assumptions of Theorem \ref{holderthm} hold. Then the
dynamical system $\left( S\left( t\right) ,\mathcal{Z}_{DN}\right) $
associated with problem (\ref{be1})-(\ref{be3}) possesses a global attractor 
$\mathcal{A}_{DN}$ in the phase-space $\mathcal{Z}_{DN}$, which is globally
bounded in $C^{\alpha }\left( \overline{\Omega }\right) ,$ $\alpha \in
\left( 0,1\right) $ and has the following structure:%
\begin{equation}
\mathcal{A}_{DN}=\Xi _{\mid t=0},  \label{f2.5}
\end{equation}%
where $\Xi $ is the set of all bounded solutions of (\ref{be1})-(\ref{be3}),
defined for all $t\in \mathbb{R}$, such that%
\begin{equation}
\left\Vert \Xi \right\Vert _{C^{\alpha /2,\alpha }\left( \mathbb{R}\times 
\overline{\Omega }\right) }\leq C,  \label{f2.6}
\end{equation}%
for some positive constant $C.$
\end{theorem}

\begin{proof}
As usual, we must check that $S\left( t\right) $ possesses a (pre)compact
absorbing set in $\mathcal{Z}_{DN}$, and that it is closed. The first
assertion is an immediate corollary of Theorem \ref{holderthm}. Indeed, for
every bounded subset $B\subset \mathcal{Z}_{DN}$, there exists a time $%
t_{+}=t_{+}\left( B\right) >0$, such that $S(t)B\subset \mathcal{X}^{\alpha
} $, for all $t\geq t_{+}$, where%
\begin{equation}
\mathcal{X}^{\alpha }:=\left\{ u\in C^{\alpha }\left( \overline{\Omega }%
\right) :\left\Vert u\right\Vert _{C^{\alpha }\left( \overline{\Omega }%
\right) }\leq C\right\}  \label{absset}
\end{equation}%
($C$ is the same constant as in (\ref{holderbound})-(\ref{holderbound2})).
Moreover, owing to Lemma \ref{cdp}, we also have that $S\left( t\right) $ is
a closed semigroup in the sense of \cite{PZ}. Thus, due to the abstract
theorem on the attractor's existence \cite[Corollary 6]{PZ}, this semigroup
possesses a global attractor $\mathcal{A}_{DN}$, bounded in $C^{\alpha
}\left( \overline{\Omega }\right) $. The characterization (\ref{f2.5}) and
estimate (\ref{f2.6}) follow from Theorem \ref{holderthm} and the same
abstract results. Lemma \ref{gl_att} is proved.
\end{proof}

\begin{remark}
Of course, the global attractor $\mathcal{A}_{DN}$ in Theorem \ref{gl_att}
depends on $M>0$ since the constant $C>0$ in (\ref{absset}) does.
Furthermore, the entering time $t_{+}>0$ is even independent of the set $%
B\subset $ $\mathcal{Z}_{DN}$, when $\sigma \left( \Gamma _{D}\right) >0$,
while it only depends on $M>0$ whenever $\sigma \left( \Gamma _{D}\right)
=0. $
\end{remark}

Finally, we can extend the above result in the critical case $m=m_{\ast }$
and $\sigma \left( \Gamma _{D}\right) =0.$

\begin{theorem}
Let assumptions (H1), (H3) and (H2) with $m=m_{\ast }$ be satisfied, and
consider the case when $\Gamma _{D}$ is empty. Let $M_{c}>0$ be the critical
mass estimated in \cite[Proposition 2 and 3]{BRB}. Then for every $M<\frac{%
M_{c}}{\left\vert \Omega \right\vert }$ such that $\left\langle
u_{0}\right\rangle \leq M$, the dynamical system $\left( S\left( t\right) ,%
\mathcal{Z}_{DN}\right) $ possesses a global attractor $\mathcal{A}%
_{M_{c}}\subset \mathcal{Z}_{DN}$, bounded in $C^{\alpha }\left( \overline{%
\Omega }\right) ,$ $\alpha \in \left( 0,1\right) .$
\end{theorem}

\begin{proof}
As in \cite{BRB}, for every $\int_{\Omega }u_{0}dx=\widetilde{M}<M_{c}$,
problem (\ref{be1})-(\ref{be3}) is globally well-posed (see, in particular, 
\cite[Theorems 7 and 9]{BRB}). The assertion in the theorem follows by
arguing exactly as in the proof of Theorem \ref{holderthm} (Case 2) since (%
\ref{Acond}) is still satisfied when $m=m_{\ast }.$
\end{proof}

\section{The non-degenerate case and smooth kernels}

\label{non_deg}

In this section, we are interested in the model proposed in \cite{BCM} which
also takes into account stochastic fluctuations based on a finite number of
individuals subject to long range attraction and short range repulsion. In
this case, the density $u\left( t\right) $ satisfies%
\begin{equation}
\partial _{t}u=\text{div}\left( \nabla A\left( u\right) -u\left( \nabla 
\mathcal{K}\ast u\right) \right) +\varepsilon \Delta u\text{, in }\Omega
\times \left( 0,\infty \right) .  \label{bed1}
\end{equation}%
The additional parameter $\varepsilon >0$ models classical Brownian random
dispersal in equation (\ref{be1}), which can be seen as the limit of (\ref%
{bed1}) as $\varepsilon \rightarrow 0$, see \cite{BCM}. For the sake of
simplicity, we will \emph{only} consider the no-flux boundary condition for (%
\ref{bed1}),%
\begin{equation}
(\nabla \left( A\left( u\right) +\varepsilon u\right) -u\left( \nabla 
\mathcal{K}\ast u\right) )\cdot \overrightarrow{n}=0\text{, on }\Gamma
\times \left( 0,\infty \right) ,  \label{bed2}
\end{equation}%
and, as before, the initial condition%
\begin{equation}
u_{\mid t=0}=u_{0}\text{ in }\Omega .  \label{bed3}
\end{equation}%
Thus, everywhere in this section, the boundary $\Gamma _{D},$ where $u=0,$\
is assumed to be empty (hence, $\Gamma \equiv \partial \Omega $).

As in Section 3, we have the following result whose proof is straightforward
(see \cite{BRB, BS}).

\begin{theorem}
\label{well-nondeg}Let (H1), (H3) be satisfied and assume%
\begin{equation}
\underset{z\rightarrow \infty }{\lim \inf }A^{^{\prime }}\left( z\right)
z^{1-m_{\ast }}>0,  \label{H2bis}
\end{equation}%
where $m_{\ast }\geq 1$ is the same critical exponent used in Section 3. If $%
m_{\ast }=1$, we further assume that each $M\geq 0$ satisfies $M<\frac{M_{c}%
}{\left\vert \Omega \right\vert }$, where $\left\langle u_{0}\right\rangle
\leq M$. Then there exists a unique (global) nonnegative solution $u\left(
t\right) =u_{\varepsilon }\left( t\right) $\ to problem (\ref{bed1})-(\ref%
{bed3}), which belongs to (\ref{reg}), and, in addition,%
\begin{equation*}
u\left( t\right) \in L^{2}\left( 0,T;H^{1}\left( \Omega \right) \right) .
\end{equation*}%
Each weak solution satisfies%
\begin{equation}
\left\langle \partial _{t}u\left( t\right) ,w\right\rangle _{\left(
H^{1}\left( \Omega \right) \right) ^{\ast },H^{1}\left( \Omega \right)
}+\int_{\Omega }\nabla A_{\varepsilon }\left( u\left( t\right) \right) \cdot
\nabla w-u\left( t\right) \left( \nabla \mathcal{K}\ast u\left( t\right)
\right) \cdot \nabla wdx=0,  \label{iden-bed}
\end{equation}%
for all $w\in H^{1}\left( \Omega \right) $, for almost all $t\in \left[ 0,T%
\right] ,$ where $A_{\varepsilon }\left( y\right) :=A\left( y\right)
+\varepsilon y,$ $\varepsilon >0.$
\end{theorem}

\begin{remark}
\label{esse2}The weak solution of (\ref{bed1})-(\ref{bed2}) satisfies the
energy identity%
\begin{equation*}
\mathcal{E}_{\varepsilon }\left( u\left( t\right) \right)
+\int_{0}^{t}\int_{\Omega }u\left( s\right) \left\vert \nabla \Phi
_{\varepsilon }^{^{\prime }}\left( u\left( s\right) \right) -\nabla \mathcal{%
K}\ast u\left( s\right) \right\vert ^{2}dxds=\mathcal{E}_{\varepsilon
}\left( 0\right) ,
\end{equation*}%
for all $t\geq 0$, see Remark \ref{esse}. The critical mass $M_{c}>0$ in the
statement of Theorem \ref{well-nondeg} is estimated in \cite[Proposition 2
and 3]{BRB}.
\end{remark}

In view of estimate (\ref{holderbound2}), the following result can be proven
for (\ref{bed1})-(\ref{bed3}).

\begin{lemma}
\label{h1lemma} Let the assumptions of Theorem \ref{att_nondeg} be
satisfied. Then, for every $\tau >0$, there exists a constant $%
C_{M,\varepsilon ,\tau }>0$ such that%
\begin{equation}
\sup_{t\geq \tau }\left[ \left\Vert u\left( t\right) \right\Vert
_{H^{1}\left( \Omega \right) }+\left\Vert \partial _{t}u\right\Vert
_{L^{2}\left( \left[ t,t+1\right] \times \Omega \right) }\right] \leq
C_{M,\varepsilon ,\tau }.  \label{h1est}
\end{equation}%
Moreover, for any bounded set $\mathcal{B}\subset \mathcal{Z}_{DN}$, there
exists a time $t_{\ast }=t_{\ast }\left( \mathcal{B}\right) >0$ such that $%
S_{\varepsilon }\left( t\right) \mathcal{B}\subset H^{1}\left( \Omega
\right) ,$ for all $t\geq t_{\ast }.$
\end{lemma}

\begin{proof}
It suffices to show (\ref{h1est}) with $\tau =t_{+}+1$ (see (\ref{est9})-(%
\ref{est10})), where $t_{+}>0$ is the time given in Theorem \ref{holderthm}
(which still applies for $A_{\varepsilon }\left( y\right) =A\left( y\right)
+\varepsilon y$). Set $\overline{A}_{\varepsilon }\left( y\right) =$ $%
\int_{0}^{y}A_{\varepsilon }\left( s\right) ds.$ Testing equation (\ref%
{iden-bed}) with $w=A_{\varepsilon }\left( u\left( t\right) \right) $, and
integrating over $\Omega ,$ we deduce%
\begin{align}
& \frac{d}{dt}\left( \overline{A}_{\varepsilon }\left( u\left( t\right)
\right) ,1\right) _{L^{1}\left( \Omega \right) }+\int_{\Omega }\left\vert
\nabla A_{\varepsilon }\left( u\left( t\right) \right) \right\vert ^{2}dx
\label{h2est} \\
& =\int_{\Omega }u\left( t\right) \left( \nabla \mathcal{K}\ast u\left(
t\right) \right) \cdot \nabla A_{\varepsilon }\left( u\left( t\right)
\right) dx  \notag \\
& \leq \frac{1}{2}\int_{\Omega }\left\vert \nabla A_{\varepsilon }\left(
u\left( t\right) \right) \right\vert ^{2}dx+C\left\Vert u\left( t\right)
\right\Vert _{L^{\infty }}^{2}\left\Vert \nabla \mathcal{K}\ast u\left(
t\right) \right\Vert _{L^{2}}^{2}  \notag
\end{align}%
Integrating this inequality from $t$ to $t+1$, and using the fact that $u$
is bounded according to (\ref{holderbound2}), we obtain%
\begin{equation}
\int_{t}^{t+1}\int_{\Omega }\left\vert \nabla A_{\varepsilon }\left( u\left(
s\right) \right) \right\vert ^{2}dx\leq C,  \label{h3est}
\end{equation}%
for all $t\geq t_{+}\left( M\right) ,$ for some positive constant $C$
independent of time and the initial data.

In order to rigorously prove (\ref{h1est}), one must proceed as for the
problem (\ref{approx}). More precisely, recalling that problem (\ref{bed1})-(%
\ref{bed3}) is uniformly parabolic, one has to employ another regularization
scheme in which $A_{\varepsilon }$ is approximated by a sequence of
functions $\left( A_{\varepsilon }\right) _{\epsilon }\in C^{\infty }\left( 
\mathbb{R}_{+}\right) ,$ the data $u_{\mid t=0}=u_{0\epsilon }\in C^{\infty
}\left( \Omega \right) \cap C^{1}\left( \overline{\Omega }\right) $ is such
that $u_{0\epsilon }\rightarrow u_{0}$ in the $L^{\infty }$-metric, and $%
\mathcal{K}$ is replaced by a sequence of smooth kernels $\mathcal{K}%
_{\epsilon }\in C^{\infty }$ such that $\mathcal{K}_{\epsilon }\rightarrow 
\mathcal{K}$ in $W_{\text{loc}}^{1,1}\left( \mathbb{R}^{d}\right) ,$ as $%
\epsilon \rightarrow 0.$ The procedure ensures that the approximate
solutions $\left( u_{\varepsilon }\right) _{\epsilon }$ are smooth enough so
that all the computations below can be performed rigorously. Thus, in what
follows we shall again proceed formally (it will be easy to see that all the
constants in the estimates are independent of $\epsilon \rightarrow 0^{+}$).
To this end, testing equation (\ref{iden-bed}) with $w=\partial _{t}\left(
A_{\varepsilon }u\left( t\right) \right) ,$ and using the fact that $%
A_{\varepsilon }^{^{\prime }}\left( y\right) \geq \varepsilon >0$, we deduce%
\begin{equation}
2\varepsilon \int_{\Omega }\left\vert \partial _{t}u\right\vert ^{2}dx+\frac{%
d}{dt}\int_{\Omega }\left\vert \nabla A_{\varepsilon }\left( u\right)
\right\vert ^{2}dx=2\int_{\Omega }u\left( \nabla \mathcal{K}\ast u\right)
\cdot \nabla \partial _{t}A_{\varepsilon }\left( u\right) dx.  \label{h4est}
\end{equation}%
Since%
\begin{align*}
\int_{\Omega }u\left( \nabla \mathcal{K}\ast u\right) \cdot \nabla \partial
_{t}A_{\varepsilon }\left( u\right) dx& =\frac{d}{dt}\int_{\Omega }u\left(
\nabla \mathcal{K}\ast u\right) \cdot \nabla A_{\varepsilon }\left( u\right)
dx \\
& -\int_{\Omega }\partial _{t}u\left( \nabla \mathcal{K}\ast u\right) \cdot
\nabla A_{\varepsilon }\left( u\right) dx \\
& -\int_{\Omega }u\left( \nabla \mathcal{K}\ast \partial _{t}u\right) \cdot
\nabla A_{\varepsilon }\left( u\right) dx,
\end{align*}%
relation (\ref{h4est}) together with basic Holder and Young inequalities
imply%
\begin{align}
& \varepsilon \int_{\Omega }\left\vert \partial _{t}u\right\vert ^{2}dx+%
\frac{d}{dt}\int_{\Omega }\left\vert \nabla A_{\varepsilon }\left( u\right)
\right\vert ^{2}dx  \label{h5est} \\
& \leq \frac{d}{dt}\int_{\Omega }2u\left( \nabla \mathcal{K}\ast u\right)
\cdot \nabla A_{\varepsilon }\left( u\right) dx+c\int_{\Omega }\left\vert
\nabla A_{\varepsilon }\left( u\right) \right\vert ^{2}dx,  \notag
\end{align}%
for some positive constant $c>0$ independent of time and the initial data.
Next, we multiply (\ref{h5est}) by $e^{c\left( t-s\right) }$ for $s\in
\left( t,t+1\right) $, to obtain%
\begin{equation}
\frac{d}{ds}e^{c\left( t-s\right) }\int_{\Omega }\left\vert \nabla
A_{\varepsilon }\left( u\right) \right\vert ^{2}dx\leq e^{c\left( t-s\right)
}\frac{d}{ds}\int_{\Omega }2u\left( \nabla \mathcal{K}\ast u\right) \cdot
\nabla A_{\varepsilon }\left( u\right) dx.  \label{h6est}
\end{equation}%
Integrating (\ref{h6est}) between $s$ and $t+1$ gives%
\begin{align}
& e^{-c}\int_{\Omega }\left\vert \nabla A_{\varepsilon }\left( u\left(
t+1\right) \right) \right\vert ^{2}dx-e^{c\left( t-s\right) }\int_{\Omega
}\left\vert \nabla A_{\varepsilon }\left( u\left( s\right) \right)
\right\vert ^{2}dx  \label{h7est} \\
& \leq \int_{s}^{t+1}e^{c\left( t-\tau \right) }\frac{d}{d\tau }\int_{\Omega
}2u\left( \tau \right) \left( \nabla \mathcal{K}\ast u\left( \tau \right)
\right) \cdot \nabla A_{\varepsilon }\left( u\left( \tau \right) \right)
dxd\tau .  \notag
\end{align}%
Notice that we can split the integral on the right-hand side of (\ref{h7est}%
) as follows:%
\begin{align*}
& \int_{s}^{t+1}e^{c\left( t-\tau \right) }\frac{d}{d\tau }\int_{\Omega
}2u\left( \tau \right) \left( \nabla \mathcal{K}\ast u\left( \tau \right)
\right) \cdot \nabla A_{\varepsilon }\left( u\left( \tau \right) \right)
dxd\tau \\
& =e^{c\left( t-\tau \right) }\int_{\Omega }2u\left( \tau \right) \left(
\nabla \mathcal{K}\ast u\left( \tau \right) \right) \cdot \nabla
A_{\varepsilon }\left( u\left( \tau \right) \right) dx\mid _{s}^{t+1} \\
& -\int_{s}^{t+1}\left( -c\right) e^{c\left( t-\tau \right) }\int_{\Omega
}2u\left( \tau \right) \left( \nabla \mathcal{K}\ast u\left( \tau \right)
\right) \cdot \nabla A_{\varepsilon }\left( u\left( \tau \right) \right)
dxd\tau \\
& =:I_{6}+I_{7}.
\end{align*}%
Next, 
\begin{align*}
I_{6}& =e^{c\left( t-\tau \right) }\int_{\Omega }2u\left( \tau \right)
\left( \nabla \mathcal{K}\ast u\left( \tau \right) \right) \cdot \nabla
A_{\varepsilon }\left( u\left( \tau \right) \right) dx\mid _{s}^{t+1} \\
& =e^{-c}\int_{\Omega }2u\left( t+1\right) \left( \nabla \mathcal{K}\ast
u\left( t+1\right) \right) \cdot \nabla A_{\varepsilon }\left( u\left(
t+1\right) \right) dx \\
& -e^{c\left( t-s\right) }\int_{\Omega }2u\left( s\right) \left( \nabla 
\mathcal{K}\ast u\left( s\right) \right) \cdot \nabla A_{\varepsilon }\left(
u\left( s\right) \right) dx,
\end{align*}%
which can be further bounded, exploiting standard Holder and Young
inequalities, by%
\begin{align}
I_{6}& \leq \frac{e^{-c}}{2}\int_{\Omega }\left\vert \nabla A_{\varepsilon
}\left( u\left( t+1\right) \right) \right\vert ^{2}dx+C\left\Vert u\left(
t+1\right) \right\Vert _{L^{\infty }\left( \Omega \right) }^{2}
\label{h8est} \\
& +\int_{\Omega }\left\vert \nabla A_{\varepsilon }\left( u\left( s\right)
\right) \right\vert ^{2}dx+C\left\Vert u\left( s\right) \right\Vert
_{L^{\infty }\left( \Omega \right) }^{2}.  \notag
\end{align}%
Moreover, we have%
\begin{eqnarray}
I_{7} &=&\int_{s}^{t+1}ce^{c\left( t-\tau \right) }\int_{\Omega }2u\left(
\tau \right) \left( \nabla \mathcal{K}\ast u\left( \tau \right) \right)
\cdot \nabla A_{\varepsilon }\left( u\left( \tau \right) \right) dxd\tau
\label{h9est} \\
&\leq &C\int_{s}^{t+1}\left( \int_{\Omega }\left\vert u\left( \tau \right)
\right\vert ^{2}dx+\int_{\Omega }\left\vert \nabla A_{\varepsilon }\left(
u\left( \tau \right) \right) \right\vert ^{2}dx\right) d\tau .  \notag
\end{eqnarray}%
Thus, on account of (\ref{h8est})-(\ref{h9est}) and estimate (\ref{supbound}%
), inequality (\ref{h7est}) becomes%
\begin{align}
& e^{-c}\int_{\Omega }\left\vert \nabla A_{\varepsilon }\left( u\left(
t+1\right) \right) \right\vert ^{2}dx-e^{c\left( t-s\right) }\int_{\Omega
}\left\vert \nabla A_{\varepsilon }\left( u\left( s\right) \right)
\right\vert ^{2}dx  \label{h10est} \\
& \leq \frac{e^{-c}}{2}\int_{\Omega }\left\vert \nabla A_{\varepsilon
}\left( u\left( t+1\right) \right) \right\vert ^{2}dx+\int_{\Omega
}\left\vert \nabla A_{\varepsilon }\left( u\left( s\right) \right)
\right\vert ^{2}dx  \notag \\
& +C\int_{s}^{t+1}\left( \int_{\Omega }\left\vert u\left( \tau \right)
\right\vert ^{2}dx+\int_{\Omega }\left\vert \nabla A_{\varepsilon }\left(
u\left( \tau \right) \right) \right\vert ^{2}dx\right) d\tau +C,  \notag
\end{align}%
for all $t\geq t_{+}.$ Therefore,%
\begin{align*}
& \frac{e^{-c}}{2}\int_{\Omega }\left\vert \nabla A_{\varepsilon }\left(
u\left( t+1\right) \right) \right\vert ^{2}dx \\
& \leq e^{c\left( t-s\right) }\int_{\Omega }\left\vert \nabla A_{\varepsilon
}\left( u\left( s\right) \right) \right\vert ^{2}dx+\int_{\Omega }\left\vert
\nabla A_{\varepsilon }\left( u\left( s\right) \right) \right\vert ^{2}dx \\
& +C\int_{s}^{t+1}\left( \int_{\Omega }\left\vert u\left( \tau \right)
\right\vert ^{2}dx+\int_{\Omega }\left\vert \nabla A_{\varepsilon }\left(
u\left( \tau \right) \right) \right\vert ^{2}dx\right) d\tau +C.
\end{align*}%
Integrating this inequality from $t$ to $t+1$ with respect to $s$, and
recalling (\ref{supbound}) we obtain%
\begin{align}
& \frac{e^{-c}}{2}\int_{\Omega }\left\vert \nabla A_{\varepsilon }\left(
u\left( t+1\right) \right) \right\vert ^{2}dx  \label{h11est} \\
& \leq \int_{t}^{t+1}e^{c\left( t-s\right) }\int_{\Omega }\left\vert \nabla
A_{\varepsilon }\left( u\left( s\right) \right) \right\vert
^{2}dxds+\int_{t}^{t+1}\int_{\Omega }\left\vert \nabla A_{\varepsilon
}\left( u\left( s\right) \right) \right\vert ^{2}dx  \notag \\
& +C\int_{t}^{t+1}\int_{s}^{t+1}\left( \int_{\Omega }\left\vert u\left( \tau
\right) \right\vert ^{2}dx+\int_{\Omega }\left\vert \nabla A_{\varepsilon
}\left( u\left( \tau \right) \right) \right\vert ^{2}dx\right) d\tau +C 
\notag \\
& \leq C\int_{t}^{t+1}\int_{\Omega }\left\vert \nabla A_{\varepsilon }\left(
u\left( s\right) \right) \right\vert ^{2}dx+C\int_{t}^{t+1}\left( \tau
-t\right) \int_{\Omega }\left\vert \nabla A_{\varepsilon }\left( u\left(
\tau \right) \right) \right\vert ^{2}dxd\tau +C.  \notag
\end{align}%
By virtue of (\ref{h3est}), (\ref{h11est}) yields%
\begin{equation}
\int_{\Omega }\left\vert \nabla A_{\varepsilon }\left( u\left( t+1\right)
\right) \right\vert ^{2}dx\leq C,\text{ }\forall t\geq t_{+}.  \label{h12est}
\end{equation}%
Since $A_{\varepsilon }^{^{\prime }}\left( y\right) \geq \varepsilon ,$ for
all $y$, and integrating (\ref{h5est}) over $\left( t,t+1\right) $ once
more, (\ref{h12est}) entails the desired estimate (\ref{h1est}). The proof
is finished.
\end{proof}

As in Section 3, we can prove the following result for (\ref{bed1})-(\ref%
{bed3}).

\begin{theorem}
\label{att_nondeg}Let the assumptions of Theorem \ref{well-nondeg} be
satisfied. The dynamical system $\left( S_{\varepsilon }\left( t\right) ,%
\mathcal{Z}_{DN}\right) $ ($\left\{ S_{\varepsilon }\left( t\right) \right\}
_{t\geq 0}$ defined as in (\ref{dyn})) possesses a global attractor $%
\mathcal{A}=\mathcal{A}_{\varepsilon ,M}$ in the sense of Definition \ref%
{def_gl}, such that $\mathcal{A}$ is globally bounded in $C^{\alpha }\left( 
\overline{\Omega }\right) \cap H^{1}\left( \Omega \right) ,$ for some $%
\alpha \in \left( 0,1\right) .$
\end{theorem}

\begin{proof}
As in the proof of Theorem \ref{gl_att}, we must check that $S_{\varepsilon
}\left( t\right) $ is closed and that it admits a compact absorbing set in $%
\mathcal{Z}_{DN}$.

\emph{Step 1} (Closedness of $S_{\varepsilon }$). The proof of this step is
essentially the same as in Lemma \ref{cdp}, where everywhere in the
estimates we must replace $A$ by the function $A_{\varepsilon }.$ This only
affects the estimate for (\ref{estI1}), which now reads%
\begin{align}
I_{1}& =-\int_{\Omega }\left( A_{\varepsilon }\left( u_{1}\right)
-A_{\varepsilon }\left( u_{2}\right) \right) \left( u_{1}-u_{2}\right)
dx+M_{12}\int_{\Omega }A_{\varepsilon }\left( u_{1}\right) -A_{\varepsilon
}\left( u_{2}\right) dx  \label{esteps} \\
& =-\varepsilon \left\Vert u_{1}-u_{2}\right\Vert _{L^{2}\left( \Omega
\right) }^{2}-\int_{\Omega }\left( A\left( u_{1}\right) -A\left(
u_{2}\right) \right) \left( u_{1}-u_{2}\right) dx  \notag \\
& +\varepsilon \left( M_{12}\right) ^{2}+M_{12}\int_{\Omega }A\left(
u_{1}\right) -A\left( u_{2}\right) dx  \notag \\
& \leq -\varepsilon \left\Vert u_{1}-u_{2}\right\Vert _{L^{2}\left( \Omega
\right) }^{2}+C\left( 1+\varepsilon \right) \left( M_{1}-M_{2}\right) ^{2}. 
\notag
\end{align}%
Thus, the same inequality in (\ref{lip2}) is valid for any two weak
solutions $u_{1}\left( t\right) ,u_{2}\left( t\right) $ of problem (\ref%
{bed1})-(\ref{bed3}) corresponding to the initial data $u_{10},u_{20}$.

\emph{Step 2} (Smoothing effect). This step requires only minor
modifications in the proof of Theorem \ref{holderthm}. Indeed, Lemma \ref%
{bounded} also applies to the function $A_{\varepsilon }$ and the results
in, e.g., \cite[Corollary 4.2]{Du} can be still applied to obtain global
Holder continuity of the weak solutions. In particular, each weak solution
of (\ref{bed1})-(\ref{bed3}) satisfies the $L^{1}$-$C^{\alpha }\cap H^{1}$
smoothing property, and there exists a time $t_{+}=t_{+}\left( M\right) >0$
such that%
\begin{equation}
\sup_{t\geq t_{+}}\left\Vert u\right\Vert _{C^{\alpha /2,\alpha }\left( 
\left[ t,t+1\right] \times \overline{\Omega }\right) \cap L^{\infty }\left(
[t_{+};\infty );H^{1}\left( \Omega \right) \right) }\leq C_{M},
\label{supbound}
\end{equation}%
for some positive constant $C_{M},$ independent of time and the initial
data. Hence, a compact in $\mathcal{Z}_{DN}$ absorbing set like in (\ref%
{absset}) can be sought, which is enough to apply \cite[Corollary 6]{PZ}
once again. The proof is finished.
\end{proof}

Next, taking advantage of the asymptotic smoothness of $S_{\varepsilon
}\left( t\right) $ we can also show that (\ref{bed1})-(\ref{bed3}) has a
gradient structure. To this end, we define the $\omega $-limit set of a
trajectory $u\left( t\right) $ of (\ref{bed1})-(\ref{bed3}), starting from%
\begin{equation*}
u_{0}\in \mathcal{Z}_{DN}=\{u_{0}\in L^{\infty }\left( \Omega \right)
:u_{0}\geq 0,\left\langle u_{0}\right\rangle \leq M\},
\end{equation*}%
as follows:%
\begin{equation}
\omega \left( u_{0}\right) =\left\{ u_{\ast }\in \mathcal{Z}_{DN}:\exists
t_{n}\rightarrow \infty \text{ such that }\lim_{n\rightarrow \infty
}\left\Vert u\left( t_{n}\right) -u_{\ast }\right\Vert _{L^{\infty }\left(
\Omega \right) }=0\right\} .  \label{oml}
\end{equation}%
We prove that $\omega \left( u_{0}\right) $ consists of stationary
solutions, satisfying the system%
\begin{equation}
\left\{ 
\begin{array}{ll}
\text{div}\left( u_{\ast }\nabla \left( \Phi _{\varepsilon }^{^{\prime
}}\left( u_{\ast }\right) -\mathcal{K}\ast u_{\ast }\right) \right) =0, & 
\text{in }\Omega , \\ 
u_{\ast }\nabla \left( \Phi _{\varepsilon }^{^{\prime }}\left( u_{\ast
}\right) -\mathcal{K}\ast u_{\ast }\right) \cdot \overrightarrow{n}=0\text{,}
& \text{on }\Gamma , \\ 
\left\langle u_{\ast }\right\rangle :=\frac{1}{\left\vert \Omega \right\vert 
}\int_{\Omega }u_{\ast }\left( x\right) ds=\left\langle u_{0}\right\rangle ,
& 
\end{array}%
\right.  \label{stat}
\end{equation}%
where we recall that $\Phi _{\varepsilon }^{^{\prime \prime }}\left(
y\right) =A_{\varepsilon }^{^{\prime }}\left( y\right) /y$ such that $\Phi
_{\varepsilon }^{^{\prime }}\left( 1\right) =\Phi _{\varepsilon }\left(
0\right) =0$ (see Remark \ref{esse}).

The following proposition justifies to call%
\begin{equation}
\mathcal{E}_{\varepsilon }\left( u\left( t\right) \right) =\int_{\Omega
}\Phi _{\varepsilon }\left( u\left( t\right) \right) dx-\frac{1}{2}%
\int_{\Omega }\int_{\Omega }u\left( x,t\right) \mathcal{K}\left( x-y\right)
u\left( y,t\right) dxdy  \label{ener_nond}
\end{equation}%
an energy functional for (\ref{bed1})-(\ref{bed3}).

\begin{proposition}
\label{omegalimit}Let $u\left( t\right) =S_{\varepsilon }\left( t\right)
u_{0},$ $u_{0}\in \mathcal{Z}_{DN}$, be the (unique) global solution of the
non-degenerate aggregation equation (\ref{bed1})-(\ref{bed3}). Then the
following assertions are true:

(i) The function $\mathcal{E}_{\varepsilon }\left( u\left( t\right) \right) $
is differentiable a.e. on $[\tau ,\infty ),$ for every $\tau >0$, and%
\begin{equation}
\frac{d}{dt}\mathcal{E}_{\varepsilon }\left( u\left( t\right) \right)
=-\int_{\Omega }u\left( t\right) \left\vert \nabla \Phi _{\varepsilon
}^{^{\prime }}\left( u\left( t\right) \right) -\nabla \mathcal{K}\ast
u\left( t\right) \right\vert ^{2}dx,  \label{nonincr}
\end{equation}%
for a.e. $t>0.$

(ii) The function $\mathcal{E}_{\varepsilon }\left( u\left( t\right) \right) 
$ is nonincreasing, and there exists a positive constant $C_{\ast },$
depending only on $\mathcal{K}$, $\Omega ,$ $\varepsilon $ and $\left\Vert
u\right\Vert _{L^{\infty }\left( \Omega \right) }$ such that%
\begin{equation}
\mathcal{E}_{\varepsilon }\left( u\left( t\right) \right) \geq -C_{\ast },
\label{bounded_bel}
\end{equation}%
for a.e. $t\geq 0$. Moreover, $\lim_{t\rightarrow \infty }\mathcal{E}%
_{\varepsilon }\left( u\left( t\right) \right) =\inf_{t>0}\mathcal{E}%
_{\varepsilon }\left( u\left( t\right) \right) =\mathcal{E}_{\varepsilon
,\infty }\in \mathbb{R}$ exists.

(iii) For any $u_{0}\in \mathcal{Z}_{DN}$, $\omega \left( u_{0}\right) $ is
a compact, connected invariant set, and every $u_{\ast }\in \omega \left(
u_{0}\right) $ is a solution of the stationary problem (\ref{stat}).

(iv) Every $u_{\ast }\in \omega \left( u_{0}\right) $ is a critical point of 
$\mathcal{E}_{\varepsilon }\left( u\left( t\right) \right) $ in (\ref%
{ener_nond}), i.e, $\mathcal{E}_{\varepsilon }^{^{\prime }}\left( u_{\ast
}\right) =0.$
\end{proposition}

\begin{proof}
First, note that by the assumption (H2) and (\ref{H2bis}), we have $\Phi
\left( y\right) \sim A\left( y\right) \sim y^{m}$ with $m\geq m_{\ast }$,
and by definition, $\Phi _{\varepsilon }\left( y\right) =\Phi \left(
y\right) +\varepsilon \left( y\ln \left( y\right) -y\right) ,$ for every $%
\varepsilon >0.$ The first part of assertion (i) follows from the fact that $%
\Phi _{\varepsilon }\in C^{1}$ and the regularity of the bounded solution $%
u\left( t\right) $ on the intervals $[\tau ,\infty )$, for every $\tau >0$
(indeed, $u\in C^{\alpha /2}\left( \tau ,\infty ;C^{\alpha }\left( \overline{%
\Omega }\right) \right) $). The second part is a consequence of Remark \ref%
{esse2}. By (\ref{nonincr}), $\mathcal{E}_{\varepsilon }\left( u\left(
t\right) \right) $ is nonincreasing on $\mathcal{Z}_{DN}.$ The elementary
inequality $y\ln \left( y\right) -y\geq -1$, for all $y\geq 0$, and the
assumption (H2) on $A$, yields on $\mathcal{Z}_{DN}$ that $\Phi
_{\varepsilon }\left( u\left( t\right) \right) \geq -C_{\varepsilon },$ for
some constant $C_{\varepsilon }>0$ depending only on the physical parameters
of the problem and the $L^{\infty }$-bound of $u$. Rewriting the energy $%
\mathcal{E}_{\varepsilon }\left( t\right) $ in the following equivalent form%
\begin{align*}
\mathcal{E}_{\varepsilon }\left( u\left( t\right) \right) & =\int_{\Omega
}\Phi _{\varepsilon }\left( u\left( t\right) \right) dx+\frac{1}{4}\int
\int_{\Omega \times \Omega }\mathcal{K}\left( x-y\right) \left( u\left(
x,t\right) -u\left( y,t\right) \right) ^{2}dxdy \\
& -\frac{1}{2}\int_{\Omega }a\left( x\right) \left( u\left( x,t\right)
\right) ^{2}dx,
\end{align*}%
where%
\begin{equation*}
a\left( x\right) =\int_{\Omega }\mathcal{K}\left( x-y\right) dy\in C\left( 
\overline{\Omega }\right) ,
\end{equation*}%
the second part of assertion (ii) is also immediate. First, by Theorem \ref%
{well-nondeg}, problem (\ref{bed1})-(\ref{bed3}) defines a nonlinear $C_{0}$%
-semigroup $S_{\varepsilon }\left( t\right) :\mathcal{Z}_{DN}\rightarrow 
\mathcal{Z}_{DN}$, $u\left( t\right) =S_{\varepsilon }\left( t\right) u_{0}$
with $u_{0}\in \mathcal{Z}_{DN}$. Second, by the results of Theorems \ref%
{holderthm} and \ref{h1lemma}, we know that for any $u_{0}\in \mathcal{Z}%
_{DN}$, there is $t_{+}>0$ such that $\cup _{t\geq t_{+}}S_{\varepsilon
}\left( t\right) u_{0}$ is bounded in $H^{1}\left( \Omega \right) \cap
C^{\alpha }\left( \overline{\Omega }\right) ,$ and hence relatively compact
in $\mathcal{Z}_{DN}$ (when endowed with the metric topology of $L^{\infty
}\left( \Omega \right) $). Third, it can be seen from (i) that the function $%
\mathcal{E}_{\varepsilon }\left( u\left( t\right) \right) :$ $\mathcal{F}%
_{DN}\rightarrow \mathbb{R}$ is a Lyapunov function on $\mathcal{F}%
_{DN}\subset \mathcal{Z}_{DN}$, for any (positively invariant) subset $%
\mathcal{F}_{DN}\subset H^{1}\left( \Omega \right) \cap C^{\alpha }\left( 
\overline{\Omega }\right) $. In particular, by Theorem \ref{att_nondeg} we
can take $\mathcal{F}_{DN}=\mathcal{A}_{\varepsilon ,M},$ where $\mathcal{A}%
_{\varepsilon ,M}$ is the global attractor for $\left( S_{\varepsilon
}\left( t\right) ,\mathcal{Z}_{DN}\right) .$ Moreover, $\mathcal{E}%
_{\varepsilon }$ satisfies: if for $t>0$, $\mathcal{E}_{\varepsilon }\left(
S_{\varepsilon }\left( t\right) u_{\ast }\right) =\mathcal{E}_{\varepsilon
}\left( u_{\ast }\right) ,$ then $u_{\ast }$ is an equilibrium point of $%
S_{\varepsilon }\left( t\right) $. In conclusion, by \cite[Chapter 10,
Definition 10.1]{R}, $\left( S_{\varepsilon }\left( t\right) ,\mathcal{A}%
_{\varepsilon ,M}\right) $ is a gradient system. Thus, by \cite[Propositions
10.3 and 10.12]{R} we immediately conclude (iii), i.e., $\omega \left(
u_{0}\right) $ is a nonempty, compact, connected invariant set, and $\omega
\left( u_{0}\right) $ consists only of stationary solutions. The final part
of (ii), $\inf_{t>0}\mathcal{E}_{\varepsilon }\left( u\left( t\right)
\right) =\mathcal{E}_{\varepsilon }\left( u_{\ast }\right) =\mathcal{E}%
_{\varepsilon ,\infty }$ is satisfied owing to (\ref{nonincr}), (\ref%
{bounded_bel}) and (iii).

Finally, for (iv) we observe that if $u_{\ast }\in \omega \left(
u_{0}\right) $ is a solution of (\ref{stat}), then for any $\psi \in
H^{1}\left( \Omega \right) \cap L^{\infty }\left( \Omega \right) $ it
follows from (\ref{stat}) that%
\begin{align*}
0& =\int_{\Omega }\text{div}\left( u_{\ast }\nabla \left( \Phi _{\varepsilon
}^{^{\prime }}\left( u_{\ast }\right) -\mathcal{K}\ast u_{\ast }\right)
\right) \psi dx-\int_{\Gamma }u_{\ast }\nabla \left( \Phi _{\varepsilon
}^{^{\prime }}\left( u_{\ast }\right) -\mathcal{K}\ast u_{\ast }\right)
\cdot \overrightarrow{n}\psi d\sigma \\
& =-\int_{\Omega }\left( \nabla A_{\varepsilon }\left( u_{\ast }\right)
-u_{\ast }\nabla \mathcal{K}\ast u_{\ast }\right) \cdot \nabla \psi dx,
\end{align*}%
which, by straightforward computations, is just the following:%
\begin{equation*}
\frac{d}{d\delta }\mathcal{E}_{\varepsilon }\left( u_{\ast }+\delta \psi
\right) _{\mid \delta =0}=0,
\end{equation*}%
i.e., $u_{\ast }$ is also a critical point of $\mathcal{E}_{\varepsilon }$
in $\mathcal{Z}_{DN}.$ In fact, we easily see that the statements (iii) and
(iv) are equivalent to each other. The proof is complete.
\end{proof}

Next, we show that $\omega \left( u_{0}\right) $ has a positive bound from
below depending only on the physical parameters of the problem and $u_{0}.$

\begin{proposition}
\label{omega_bb}Let $u\left( t\right) =S_{\varepsilon }\left( t\right)
u_{0}, $ $u_{0}\in \mathcal{Z}_{DN},$ be the unique solution of (\ref{bed1}%
)-(\ref{bed3}) such that $\left\langle u_{0}\right\rangle >0$. Suppose that $%
u_{\ast }\in \omega \left( u_{0}\right) $. Then, there exists a constant $%
\underline{u}>0$, depending only on $u_{0},$ $\Omega ,$ $\mathcal{K}$ and $%
\varepsilon >0,$ such that $u_{\ast }\left( x\right) \geq \underline{u}>0,$
for all $x\in \Omega .$
\end{proposition}

\begin{proof}
The proof follows a similar argument used in \cite[Section 2]{FLP}, \cite[%
Proposition 3.3]{ZZ}. Indeed, owing to the homogeneous Neumann boundary
condition of (\ref{stat}), the first equation of (\ref{stat}) also reads%
\begin{equation*}
\int_{\Omega }u_{\ast }\nabla \left( \Phi _{\varepsilon }^{^{\prime }}\left(
u_{\ast }\right) -\mathcal{K}\ast u_{\ast }\right) \cdot \nabla \psi dx=0,
\end{equation*}%
for any $\psi \in C^{1}\left( \overline{\Omega }\right) $, whence $\nabla
(\Phi _{\varepsilon }^{^{\prime }}\left( u_{\ast }\right) -\mathcal{K}\ast
u_{\ast })=0$ in each connected component of the open set where $x\in \Omega 
$ such that $u_{\ast }\left( x\right) >0.$ Therefore, $\Phi _{\varepsilon
}^{^{\prime }}\left( u_{\ast }\right) -\mathcal{K}\ast u_{\ast }$ is
constant on each such connected component$.$ If $\left\langle
u_{0}\right\rangle >0$, by the last equation of (\ref{stat}), there exists $%
x_{0}\in \Omega $ such that $u_{\ast }\left( x_{0}\right) >0.$ The previous
statement implies that there exists $\gamma \in \mathbb{R}$ such that 
\begin{equation}
\gamma =\Phi _{\varepsilon }^{^{\prime }}\left( u_{\ast }\left( x\right)
\right) -\left( \mathcal{K}\ast u_{\ast }\right) \left( x\right) =\Phi
^{^{\prime }}\left( u_{\ast }\left( x\right) \right) +\varepsilon \ln \left(
u_{\ast }\left( x\right) \right) -\left( \mathcal{K}\ast u_{\ast }\right)
\left( x\right) ,  \label{co_co}
\end{equation}%
for $x\in \Omega \left( x_{0}\right) ,$ which is the same connected
component of $\left\{ x\in \Omega :u_{\ast }\left( x\right) >0\right\} $ as $%
x_{0}.$ Since $u_{\ast }\in \omega \left( u_{0}\right) $ is bounded and $%
\mathcal{K}\in W_{\text{loc}}^{1,1}\left( \mathbb{R}^{d}\right) $, we
observe from (\ref{co_co}) that $u_{\ast }\left( x\right) $ satisfies the
inequality%
\begin{equation*}
u_{\ast }\left( x\right) \geq e^{\frac{1}{\varepsilon }\left( \gamma
-\left\Vert \mathcal{K}\ast u_{\ast }\right\Vert _{L^{\infty }\left( \Omega
\right) }-\left\Vert \Phi ^{^{\prime }}\left( u_{\ast }\left( x\right)
\right) \right\Vert _{L^{\infty }\left( \Omega \right) }\right) }>0,
\end{equation*}%
for all $x\in \Omega \left( x_{0}\right) .$ In particular, this yields that $%
\Omega \left( x_{0}\right) =\Omega $ and the claim is proved.
\end{proof}

\begin{remark}
Some results on properties of the steady states for the aggregation equation
with nonlinear diffusion (\ref{be1}) in the case $\Omega =\mathbb{R}^{d}$
can be found \cite{Be, BF, BF2}.
\end{remark}

As a result of the proof of Proposition \ref{omegalimit}, we can now
conclude the following

\begin{theorem}
\label{unst_gl}Let the assumptions of Theorem \ref{att_nondeg} be satisfied.
The global attractor $\mathcal{A}=\mathcal{A}_{\varepsilon ,M}$ of problem (%
\ref{bed1})-(\ref{bed3}) consists entirely of unstable manifolds of the
equilibria, which are bounded solutions of (\ref{stat}).
\end{theorem}

Under an additional assumption on the kernel which ensures that $\mathcal{K}$
is reasonably smooth at the origin, we can show that $\mathcal{A}%
_{\varepsilon ,M}$ is also finite dimensional.

\begin{theorem}
\label{fin_gl}Let the assumptions of Theorem \ref{well-nondeg} be satisfied.
Assume%
\begin{equation}
\left( D^{2}\mathcal{K}\right) 1_{B_{1}\left( 0\right) }\in L^{1}\left( 
\mathbb{R}^{d}\right) ,\text{ }d\geq 2.  \label{smoothK}
\end{equation}%
The global attractor $\mathcal{A}_{\varepsilon ,M}$ of (\ref{bed1})-(\ref%
{bed3}) has finite fractal dimension:%
\begin{equation*}
\dim _{F}\left( \mathcal{A}_{\varepsilon ,M},H^{1-}\left( \Omega \right)
\cap L^{p}\left( \Omega \right) \right) \leq C_{\varepsilon ,M}<\infty ,
\end{equation*}%
for any $1<p<\infty .$
\end{theorem}

\begin{remark}
Recall that fractal dimension of a compact set $\mathcal{Y},$ $\dim
_{F}\left( \mathcal{Y},X\right) $ is defined as%
\begin{equation*}
\dim _{F}\left( \mathcal{Y},X\right) =\lim_{\delta \rightarrow 0^{+}}\frac{%
\log _{2}N_{\delta }\left( \mathcal{Y},X\right) }{\log _{2}\left( 1/\delta
\right) },
\end{equation*}%
where $N_{\delta }\left( \mathcal{Y},X\right) $ is the minimal number of
balls $B_{\delta }$ that can be used to cover the compact set $\mathcal{Y}$
in the metric of $X.$
\end{remark}

The statement of Theorem \ref{fin_gl} is in fact a consequence of a much
stronger result which states that (\ref{bed1})-(\ref{bed3}) admits an
exponential attractor $\mathcal{M}_{\varepsilon ,M}$ provided that (\ref%
{smoothK}) is also satisfied. The precise statement is given by the
following.

\begin{theorem}
\label{expo} Let the assumptions of Theorem \ref{fin_gl} be satisfied.
Assume (\ref{smoothK}). For every fixed $\varepsilon >0$, there exists an
exponential attractor $\mathcal{M}=\mathcal{M}_{\varepsilon ,M}$ bounded in $%
C^{\alpha }\left( \overline{\Omega }\right) $ for the dynamical system $%
\left( S_{\varepsilon }\left( t\right) ,\mathcal{Z}_{DN}\right) $ which
satisfies the following properties:

(i) Semi-invariance: $S_{\varepsilon }\left( t\right) \mathcal{M}\subset 
\mathcal{M}$, for every $t\geq 0.$

(ii) Exponential attraction: For every bounded subset $\mathcal{B}\subset 
\mathcal{Z}_{DN},$ 
\begin{equation}
dist_{H^{1-}\left( \Omega \right) \cap L^{p}\left( \Omega \right) }\left(
S_{\varepsilon }\left( t\right) \mathcal{B},\mathcal{M}\right) \leq
Ce^{-\kappa t},\quad \forall t\geq 0,  \label{abs1}
\end{equation}%
for some positive constants $C=C\left( \varepsilon ,M\right) $ and $\kappa $%
, for any $1<p<\infty .$

(iii) Finite dimensionality:%
\begin{equation}
\dim _{F}\left( \mathcal{M},H^{1-}\left( \Omega \right) \cap L^{p}\left(
\Omega \right) \right) \leq C_{\varepsilon ,M}<\infty ,  \label{abs2}
\end{equation}%
for any $1<p<\infty .$ The constants $C,C_{\varepsilon ,M},\kappa $ can be
computed explicitly in terms of the physical parameters of the problem.
\end{theorem}

Here, and everywhere else, we denote $H^{1-}\left( \Omega \right)
:=H^{1-\delta }\left( \Omega \right) $, for any $\delta \in (0,1]$. Since
the global attractor $\mathcal{A}_{\varepsilon ,M}$ is always contained in $%
\mathcal{M}_{\varepsilon ,M}$, the above theorem immediately implies that
the fractal dimension of the global attractor $\mathcal{A}_{\varepsilon ,M}$
is also finite. It is worth to recall that, in the global attractors theory,
it is usually extremely difficult (if not impossible) to estimate and to
express the rate of convergence in (\ref{gl_conv}) in terms of the physical
parameters of the system considered. This constitutes the main drawback of
the theory. Simple examples show that the rate of convergence in (\ref%
{gl_conv}) can be arbitrarily slow and non-uniform with respect to the
parameters of the system considered. As a consequence, the global attractor
becomes sensitive to small perturbations and, in some sense, cannot even be
observed in experiments. The concept of exponential attractor overcomes this
difficulty (see, e.g., the survey article \cite{MZ}). Indeed, in contrast to
the global attractors theory, the constants $C,C_{\varepsilon ,M},\kappa $
in (\ref{abs1})-(\ref{abs2}) can be explicitly found in terms of the
physical parameters.

We report for the reader's convenience the following abstract result on the
existence of exponential attractors \cite[Proposition 4.1]{EZ} (cf. also 
\cite[Proposition 2.17]{GM}) which will be used in the proof of Theorem \ref%
{expo}.

\begin{proposition}
\label{abstract}Let $\mathcal{H}$,$\mathcal{V}$,$\mathcal{V}_{1}$ be Banach
spaces such that the embedding $\mathcal{V}_{1}\subset \mathcal{V}$ is
compact. Let $B$ be a closed bounded subset of $\mathcal{H}$ and let $%
\mathbb{S}:B\rightarrow B$ be a map. Assume also that there exists a
uniformly Lipschitz continuous map $\mathbb{T}:B\rightarrow \mathcal{V}_{1}$%
, i.e.,%
\begin{equation}
\left\Vert \mathbb{T}b_{1}-\mathbb{T}b_{2}\right\Vert _{\mathcal{V}_{1}}\leq
L\left\Vert b_{1}-b_{2}\right\Vert _{\mathcal{H}},\quad \forall
b_{1},b_{2}\in B,  \label{gl1}
\end{equation}%
for some $L\geq 0$, such that%
\begin{equation}
\left\Vert \mathbb{S}b_{1}-\mathbb{S}b_{2}\right\Vert _{\mathcal{H}}\leq
\theta \left\Vert b_{1}-b_{2}\right\Vert _{\mathcal{H}}+K\left\Vert \mathbb{T%
}b_{1}-\mathbb{T}b_{2}\right\Vert _{\mathcal{V}},\quad \forall
b_{1},b_{2}\in B,  \label{gl2}
\end{equation}%
for some $\theta <\frac{1}{2}$ and $K\geq 0$. Then, there exists a
(discrete) exponential attractor $\mathcal{M}_{d}\subset B$ of the semigroup 
$\{\mathbb{S}(n):=\mathbb{S}^{n},n\in Z+\}$ with discrete time in the phase
space $\mathcal{H}$.
\end{proposition}

We delay the proof of Theorem \ref{expo} until the end of the section. The
idea is to verify (\ref{gl1})-(\ref{gl2}) for a suitable choice of maps. We
begin by showing that the semigroup $S_{\varepsilon }\left( t\right) $ is
strongly (Lipschitz) continuous with respect to the $\left( H^{1}\left(
\Omega \right) \right) ^{\ast }$-metric.

\begin{proposition}
\label{uniq} Let $u_{i},$ $i=1,2$, be a pair of weak solutions according to
the assumptions of Theorem \ref{fin_gl}. Then the following estimate holds: 
\begin{align}
& \left\Vert u_{1}\left( t\right) -u_{2}\left( t\right) \right\Vert _{\left(
H^{1}\left( \Omega \right) \right) ^{\ast }}^{2}+\varepsilon
\int_{0}^{t}\left\Vert u_{1}\left( s\right) -u_{2}\left( s\right)
\right\Vert _{L^{2}\left( \Omega \right) }^{2}ds  \label{Lipschitz} \\
& \leq \left\Vert u_{1}\left( 0\right) -u_{2}\left( 0\right) \right\Vert
_{\left( H^{1}\left( \Omega \right) \right) ^{\ast }}^{2}e^{\kappa t}, 
\notag
\end{align}%
for all $t\geq 0$, for some positive constants $\kappa ,C$ which depend on $%
\varepsilon >0$ and $\mathcal{K}$ but are independent of $u_{i}\left(
0\right) .$
\end{proposition}

\begin{proof}
Following Lemma \ref{cdp}, we have that $u:=u_{1}-u_{2}$ and $\eta \left(
t\right) :=\left\Vert u\left( t\right) \right\Vert _{\left( H^{1}\right)
^{\ast }}^{2}$ satisfies the problem%
\begin{equation}
\frac{1}{2}\frac{d}{dt}\eta \left( t\right) =-\left\langle \partial
_{t}u\left( t\right) ,\phi \left( t\right) \right\rangle =I_{1}+I_{2}+I_{3},
\notag
\end{equation}%
with $I_{1},I_{2},I_{3}$ given in (\ref{Is}). The integral $I_{1}$ is
estimated in (\ref{esteps}), whereas for $I_{2},I_{3}$ we have (\ref{estu2}%
). The estimate for $I_{2}$ can be improved in (\ref{estu2}), using (\ref%
{smoothK}), Young's inequality for convolutions and Lemma \ref{kest}, (b),
as follows:%
\begin{align*}
I_{2}& \leq C\int_{\Omega }\left\vert D^{2}\mathcal{K}\ast u_{1}\right\vert
\left\vert \nabla \phi \right\vert ^{2}dx\leq C\left\Vert D^{2}\mathcal{K}%
\ast u_{1}\right\Vert _{L^{\infty }\left( \Omega \right) }\left\Vert \nabla
\phi \right\Vert _{L^{2}\left( \Omega \right) }^{2} \\
& \leq C\left( \left\Vert D^{2}\mathcal{K}\right\Vert _{L^{1}\left(
B_{1}\left( 0\right) \right) }\left\Vert u\right\Vert _{L^{\infty }\left(
\Omega \right) }+\left\Vert D^{2}\mathcal{K}\right\Vert _{L^{p}\left( 
\mathbb{R}^{d}\backslash B_{1}\left( 0\right) \right) }\left\Vert
u\right\Vert _{L^{p/\left( p-1\right) }\left( \Omega \right) }\right)
\left\Vert \nabla \phi \right\Vert _{L^{2}\left( \Omega \right) }^{2} \\
& \leq C\eta \left( t\right) .
\end{align*}

Thus, we get%
\begin{equation}
\frac{d}{dt}\left\Vert u\left( t\right) \right\Vert _{\left( H^{1}\right)
^{\ast }}^{2}+2\varepsilon \left\Vert u\left( t\right) \right\Vert
_{L^{2}}^{2}\leq C\left\Vert u\left( t\right) \right\Vert _{\left(
H^{1}\right) ^{\ast }}^{2},  \label{decay1b}
\end{equation}%
which yields the desired inequality (\ref{Lipschitz}) by application of the
Gronwall's inequality.
\end{proof}

\begin{remark}
A crucial point in the proof of Theorem \ref{expo} is that we need the
global Lipschitz continuity of $S_{\varepsilon }\left( t\right) $ in the
norm of $\left( H^{1}\left( \Omega \right) \right) ^{\ast }.$ The assumption
(\ref{smoothK}) plays an essential role with respect to this issue (see (\ref%
{decay1b})). While Newtonian potentials do not satisfy (\ref{smoothK}), in
population dynamics the non-local effects are generally modelled with
smooth, fast-decaying kernels $\mathcal{K}$ which obey (\ref{smoothK}), see
e.g., \cite{TBL}.
\end{remark}

The step needed to establish the existence of an exponential attractor is
the validity of so-called smoothing property for the difference of two
solutions of (\ref{bed1})-(\ref{bed3}). In the present case, such a property
is a consequence of the following two lemmas. The first result establishes
that the semigroup $S_{\varepsilon }\left( t\right) $ is some kind of
contraction map, up to the term $\left\Vert u_{1}-u_{2}\right\Vert _{L^{2}(%
\left[ 0,t\right] ;\left( H^{1}\right) ^{\ast })}$.

\begin{lemma}
\label{dec} Let the assumptions of Proposition \ref{uniq} hold. Then, for
every $t\geq 0$, we have:%
\begin{equation}
\left\Vert u_{1}\left( t\right) -u_{2}\left( t\right) \right\Vert _{\left(
H^{1}\right) ^{\ast }}^{2}\leq e^{-\kappa t}\left\Vert u_{1}\left( 0\right)
-u_{2}\left( 0\right) \right\Vert _{\left( H^{1}\right) ^{\ast
}}^{2}+C_{\varepsilon ,M}\int_{0}^{t}\left\Vert u_{1}\left( s\right)
-u_{2}\left( s\right) \right\Vert _{\left( H^{1}\right) ^{\ast }}^{2}ds,
\label{decay1d}
\end{equation}%
for some positive constants $C_{\varepsilon ,M}$, $\kappa $ which depend on $%
\varepsilon >0,$ $\Omega $ and $\mathcal{K}.$
\end{lemma}

\begin{proof}
Recall that $u:=u_{1}-u_{2}$. Combining (\ref{decay1b}) together with Poincar%
\'{e}'s inequality 
\begin{equation*}
\left\Vert A_{N}^{-1/2}\left( u-\left\langle u\right\rangle \right)
\right\Vert _{L^{2}}^{2}+\left\langle u\right\rangle ^{2}\leq C_{\Omega
}\left\Vert u\right\Vert _{L^{2}}^{2},
\end{equation*}%
we deduce from (\ref{decay1b}) the following inequality:%
\begin{equation*}
\frac{d}{dt}\left\Vert u\left( t\right) \right\Vert _{\left( H^{1}\right)
^{\ast }}^{2}+\frac{2\varepsilon }{C_{\Omega }}\left\Vert u\left( t\right)
\right\Vert _{\left( H^{1}\right) ^{\ast }}^{2}\leq C\left\Vert u\left(
t\right) \right\Vert _{\left( H^{1}\right) ^{\ast }}^{2},
\end{equation*}%
for all $t\geq 0$. Thus, Gronwall's inequality entails the desired estimate (%
\ref{decay1d}).
\end{proof}

We now need some compactness for the term $\left\Vert u_{1}-u_{2}\right\Vert
_{L^{2}\left( \left[ 0,t\right] ;\left( H^{1}\right) ^{\ast }\right) }$ on
the right-hand side of (\ref{decay1d}). This is given by

\begin{lemma}
\label{lipdif} Let the assumptions of Proposition \ref{uniq} hold. Then, for
every $t\geq 0$, the following estimate holds:%
\begin{align}
& \left\Vert \partial _{t}u_{1}-\partial _{t}u_{2}\right\Vert _{L^{2}\left( 
\left[ 0,t\right] ;D\left( A_{N}\right) ^{^{\prime }}\right)
}^{2}+\varepsilon \int_{0}^{t}\left\Vert u_{1}\left( s\right) -u_{2}\left(
s\right) \right\Vert _{L^{2}\left( \Omega \right) }^{2}ds  \label{comp1} \\
& \leq C_{\varepsilon ,M}e^{\kappa t}\left\Vert u_{1}\left( 0\right)
-u_{2}\left( 0\right) \right\Vert _{\left( H^{1}\right) ^{\ast }}^{2}, 
\notag
\end{align}%
where $C_{\varepsilon ,M}>0$ and $\kappa >0$ also depend on $\varepsilon ,$ $%
\Omega $ and $\mathcal{K}.$
\end{lemma}

\begin{proof}
The second term on the left-hand side of (\ref{comp1}) can be easily
controlled by (\ref{Lipschitz}). Thus we only need to estimate the time
derivative. Recall that each $\partial _{t}u_{i},$ $i=1,2$, satisfies (\ref%
{iden-bed}). Furthermore, in light of Theorem \ref{well-nondeg}, recall that
we have%
\begin{equation}
\sup_{t\geq 0}\left\Vert u_{i}\left( t\right) \right\Vert _{L^{\infty
}\left( \Omega \right) }\leq C_{\varepsilon ,M}\text{, }i=1,2.
\label{linfbis}
\end{equation}%
Thus, for any test function $w\in D(A_{N})$, using the weak formulation (\ref%
{iden-bed}), for $\partial _{t}u:=\partial _{t}u_{1}-\partial _{t}u_{2}$,
there holds 
\begin{equation*}
\left\langle \partial _{t}u\left( t\right) ,w\right\rangle =I_{4}+I_{5},
\end{equation*}%
where 
\begin{align*}
I_{4}& :=-\left\langle \nabla \left( A_{\varepsilon }\left( u_{1}\right)
-A_{\varepsilon }\left( u_{2}\right) \right) ,\nabla w\right\rangle ,\text{ }
\\
I_{5}& :=\left\langle u\left( \nabla \mathcal{K}\ast u_{1}\right)
-u_{2}\left( \nabla \mathcal{K}\ast u\right) ,\nabla w\right\rangle .
\end{align*}%
First, for every $w\in D\left( A_{N}\right) $ we have%
\begin{align*}
& I_{4}=\left\langle A_{\varepsilon }\left( u_{1}\right) -A_{\varepsilon
}\left( u_{2}\right) ,\Delta w\right\rangle \leq Q\left( \left\Vert
u_{i}\right\Vert _{L^{\infty }\left( \Omega \right) }\right) \left\Vert
u_{1}-u_{2}\right\Vert _{L^{2}}\left\Vert \Delta w\right\Vert _{L^{2}} \\
& \leq C\left\Vert u_{1}-u_{2}\right\Vert _{L^{2}}\left\Vert w\right\Vert
_{D\left( A_{N}\right) }.
\end{align*}%
On the other hand, it is easy to show, on account of (\ref{linfbis}), that%
\begin{equation*}
I_{5}\leq C\left\Vert u_{1}-u_{2}\right\Vert _{L^{2}}\left\Vert \nabla
w\right\Vert _{L^{2}}.
\end{equation*}%
These estimates together with (\ref{Lipschitz}) gives the desired estimate
on the time derivative in (\ref{comp1}).
\end{proof}

We now show that the semigroup $S_{\varepsilon }\left( t\right) $ is
actually uniformly H\"{o}lder continuous in the $H^{1-}\cap L^{p}$-norm with
respect to the initial data.

\begin{lemma}
\label{holder} Let $u_{i}\left( t\right) =S\left( t\right) u_{i}\left(
0\right) $, with $u_{i}(0)\in \mathcal{Z}_{DN}$. Then, for any $1<p<\infty $%
, the following estimate is valid:%
\begin{equation}
\left\Vert u_{1}\left( t\right) -u_{2}\left( t\right) \right\Vert
_{H^{1-}\left( \Omega \right) \cap L^{p}\left( \Omega \right) }\leq
C_{\varepsilon ,M}e^{\kappa t}\left\Vert u_{1}\left( 0\right) -u_{2}\left(
0\right) \right\Vert _{\left( H^{1}\right) ^{\ast }}^{\gamma },
\label{hclinf}
\end{equation}%
for all $t\geq t_{\ast }$, where the constants $C_{\varepsilon ,M}>0,$ $%
\kappa >0$ and $\gamma =\gamma \left( p\right) <1$ are independent of the
initial data and time.
\end{lemma}

\begin{proof}
Using the interpolation $[H^{1},\left( H^{1}\right) ^{\ast }]_{1/2,2}=L^{2}$%
, we deduce from estimates (\ref{Lipschitz}) and (\ref{h1est}) that 
\begin{equation}
\left\Vert u_{1}\left( t\right) -u_{2}\left( t\right) \right\Vert
_{L^{2}\left( \Omega \right) }\leq C_{\varepsilon ,M}e^{\kappa t}\left\Vert
u_{1}\left( 0\right) -u_{2}\left( 0\right) \right\Vert _{\left( H^{1}\right)
^{\ast }}^{1/2},  \label{hc2}
\end{equation}%
for all $t\geq t_{\ast }=t_{+}+1$. By interpolation in the spaces $L^{\infty
}\subset L^{p}\subset L^{2}$, for $2<p<\infty $, $H^{1}\subset H^{1-\delta
}\subset \left( H^{1}\right) ^{\ast }$, $\delta \in (0,1)$, the estimate (%
\ref{hclinf}) also holds for the difference of solutions $u=u_{1}-u_{2}$.
\end{proof}

The last ingredient we need is the uniform H\"{o}lder continuity of $%
t\mapsto S_{\varepsilon }(t)u_{0}$ in the $H^{1-}\cap L^{p}$-norm, namely,

\begin{lemma}
\label{hctime} Let the assumptions of Proposition \ref{uniq} be satisfied.
Consider $u\left( t\right) =S_{\varepsilon }\left( t\right) u_{0}$ with $%
u_{0}\in \mathcal{Z}_{DN}$. The following estimate holds 
\begin{equation}
\left\Vert u\left( t\right) -u\left( s\right) \right\Vert _{H^{1-}\left(
\Omega \right) \cap L^{p}\left( \Omega \right) }\leq C_{\varepsilon
,M}\left\vert t-s\right\vert ^{\gamma },\quad \forall t,s\geq t_{\ast },
\label{hclinf2}
\end{equation}%
where $\gamma =\gamma \left( p\right) <1$ and the positive constant $%
C_{\varepsilon ,M}$ is independent of initial data, $u$ and $t,s$.
\end{lemma}

\begin{proof}
According to (\ref{supbound}) and (\ref{h1est}), and recalling that $%
\mathcal{K}\in W_{\text{loc}}^{1,1}\left( \mathbb{R}^{d}\right) $, the
following bound holds: 
\begin{equation*}
\sup_{t\geq t_{\ast }}\left\Vert \Delta A_{\varepsilon }u\left( t\right)
-\nabla \cdot \left( u\left( \nabla \mathcal{K}\ast u\left( t\right) \right)
\right) \right\Vert _{\left( H^{1}\right) ^{\ast }}\leq C_{\varepsilon ,M}.
\end{equation*}%
Consequently, by comparison in (\ref{iden-bed}), we have that 
\begin{equation*}
\sup_{t\geq t_{\ast }}\left\Vert \partial _{t}u\left( t\right) \right\Vert
_{\left( H^{1}\right) ^{\ast }}\leq C_{\varepsilon ,M},
\end{equation*}%
which entails 
\begin{equation}
\left\Vert u\left( t\right) -u\left( s\right) \right\Vert _{\left(
H^{1}\right) ^{\ast }}\leq C_{\varepsilon ,M}\left\vert t-s\right\vert
,\quad \forall t,s\geq t_{\ast }.  \label{hcb3}
\end{equation}%
Estimate (\ref{hclinf2}) now follows from (\ref{hcb3}), and standard
interpolation inequalities, as in the proof of Lemma \ref{holder}.
\end{proof}

\textit{Proof of Theorem \ref{expo}}. In order to apply Proposition \ref%
{abstract}, it is sufficient to verify the existence of an exponential
attractor for the restriction of $S(t)$ on some properly chosen
semi-invariant absorbing set in $\mathcal{Z}_{DN}$. Recall that, by (\ref%
{supbound}) and Lemma \ref{h1lemma}, the ball $B_{0}:=B_{C^{\alpha }\left( 
\overline{\Omega }\right) \cap H^{1}\left( \Omega \right) }\left(
C_{\varepsilon ,M}\right) $ will be absorbing for $S_{\varepsilon }\left(
t\right) $, provided that $C_{\varepsilon ,M}>0$ is sufficiently large.
Since we want this ball to be semi-invariant with respect to the semigroup,
we push it forward by the semigroup, by defining first the set $B_{1}=\left[
\cup _{t\geq 0}S_{\varepsilon }\left( t\right) B_{0}\right] _{\left(
H^{1}\right) ^{\ast }}$, where $\left[ \cdot \right] _{\left( H^{1}\right)
^{\ast }}$ denotes closure in the space $\left( H^{1}\left( \Omega \right)
\right) ^{\ast }.$ Then set $\mathbb{B}=S\left( 1\right) B_{1}$. Thus, $%
\mathbb{B}$ is a semi-invariant compact (for the metric of $\left(
H^{1}\left( \Omega \right) \right) ^{\ast }$) subset of the phase space $%
\mathcal{Z}_{DN}$. On the other hand, due to the results proven in this
section, we have%
\begin{equation}
\sup_{t\geq 0}\left( \left\Vert u\left( t\right) \right\Vert _{C^{\alpha
}\left( \overline{\Omega }\right) \cap H^{1}\left( \Omega \right)
}+\left\Vert A_{\varepsilon }u\left( t\right) \right\Vert _{H^{1}\left(
\Omega \right) }+\left\Vert \partial _{t}u\left( t\right) \right\Vert
_{\left( H^{1}\right) ^{\ast }}\right) \leq C_{\varepsilon ,M},
\label{reg_sec}
\end{equation}%
for every trajectory $u$ originating from $u_{0}\in \mathbb{B}$, for some
positive constant $C_{\varepsilon ,M}$ which is independent of the choice of 
$u_{0}\in \mathbb{B}$. We can now apply the abstract result above to the map 
$\mathbb{S}=S_{\varepsilon }\left( T\right) $ and $\mathcal{H}=\left(
H^{1}\left( \Omega \right) \right) ^{\ast }$, for a fixed $T>0$ such that $%
e^{-\kappa T}<\frac{1}{2}$, where $\kappa >0$ is the same as in Lemma \ref%
{dec}. To this end, we introduce the functional spaces%
\begin{equation}
\mathcal{V}_{1}:=L^{2}\left( \left[ 0,T\right] ;L^{2}\left( \Omega \right)
\right) \cap H^{1}\left( \left[ 0,T\right] ;D(A_{N})^{\prime }\right) ,\quad 
\mathcal{V}:=L^{2}\left( \left[ 0,T\right] ;\left( H^{1}\left( \Omega
\right) \right) ^{\ast }\right) ,  \label{fsp}
\end{equation}%
and note that $\mathcal{V}_{1}$ is compactly embedded into $\mathcal{V}$.
Finally, we introduce the operator $\mathbb{T}:\mathbb{B}\rightarrow 
\mathcal{V}_{1}$, by $\mathbb{T}u_{0}:=u\in \mathcal{V}_{1},$ where $u$
solves (\ref{bed1})-(\ref{bed3}) with $u\left( 0\right) =u_{0}\in \mathbb{B}$%
. We claim that the maps $\mathbb{S}$, $\mathbb{T}$, the spaces $\mathcal{H}$%
,$\mathcal{V}$,$\mathcal{V}_{1}$ thus defined satisfy all the assumptions of
Proposition \ref{abstract}. Indeed, the global Lipschitz continuity (\ref%
{gl1}) of $\mathbb{T}$ is an immediate corollary of Lemma \ref{lipdif}, and
estimate (\ref{gl2}) follows from estimate (\ref{decay1d}). Therefore, due
to Proposition \ref{abstract}, the semigroup $\mathbb{S}(n)=S_{\varepsilon
}\left( nT\right) $ generated by the iterations of the operator $\mathbb{S}:%
\mathbb{B\rightarrow B}$ possesses a (discrete) exponential attractor $%
\mathcal{M}_{d}$ in $\mathbb{B}$ endowed by the topology of $\left(
H^{1}\left( \Omega \right) \right) ^{\ast }$. In order to construct the
exponential attractor $\mathcal{E}$ for the semigroup $S_{\varepsilon }(t)$
with continuous time, we note that, due to Lemma \ref{uniq}, this semigroup
is Lipschitz continuous with respect to the initial data in the topology of $%
\left( H^{1}\left( \Omega \right) \right) ^{\ast }$. Moreover, by (\ref%
{hclinf}) and (\ref{hclinf2}) the map $\left( t,u_{0}\right) \mapsto
S_{\varepsilon }\left( t\right) u_{0}$ is also uniformly H\"{o}lder
continuous on $\left[ 0,T\right] \times \mathbb{B}$, where $\mathbb{B}$ is
endowed with the metric topology of $\left( H^{1}\left( \Omega \right)
\right) ^{\ast }$. Hence, the desired exponential attractor $\mathcal{M}$
for the continuous semigroup $S_{\varepsilon }(t)$ can be obtained by the
standard formula%
\begin{equation}
\mathcal{M}=\bigcup_{t\in \left[ 0,T\right] }S_{\varepsilon }\left( t\right) 
\mathcal{M}_{d}.  \label{st}
\end{equation}%
In order to finish the proof of the theorem, we only need to verify that $%
\mathcal{M}$ defined as above will be the exponential attractor for $%
S_{\varepsilon }(t)$ restricted to $\mathbb{B}$ not only with respect to the 
$\left( H^{1}\left( \Omega \right) \right) ^{\ast }$-metric, but also in
with respect to a stronger metric. This is an immediate corollary of the
following facts: $\mathbb{B}$ is bounded in $C^{\alpha }\left( \overline{%
\Omega }\right) \cap H^{1}\left( \Omega \right) $ and standard interpolation
inequalities between the following spaces: $L^{\infty }\subset L^{p}\subset
L^{2}$, $H^{1}\subset H^{1-\delta }\subset \left( H^{1}\right) ^{\ast },$
for $2<p<\infty $ and $\delta \in (0,1].$ Theorem \ref{expo} is now proved.

\section{Convergence to steady states}

\label{sta}

In this section, we show that any global-in-time bounded solution to the
model (\ref{bed1})-(\ref{bed3}) converges to a single equilibrium of (\ref%
{stat}) as time tends to infinity. The proof of the main result is based on
a suitable version of the Lojasiewicz--Simon theorem and Propositions \ref%
{omegalimit}, \ref{omega_bb}. The question of such convergence is usually a
delicate matter since it is well known that the topology of the set of
stationary solutions of (\ref{stat}) can be non-trivial. In particular,
there may be a continuum of stationary solutions for (\ref{stat}) even in
the simplest cases, for instance when $\Omega $ is a disk, $\mathcal{K}$ is
either a Newtonian or Bessel potential and $A\left( y\right) \equiv y$ (see 
\cite{Har, SS}).

The main result of this section reads as follows.

\begin{theorem}
\label{conv_equil}Let the assumptions of Theorem \ref{well-nondeg} be
satisfied. Assume that $\Phi $ is a real analytic function on $\mathbb{R}%
_{+} $. For any $u_{0}\in \mathcal{Z}_{DN}$ with $\left\langle
u_{0}\right\rangle >0$, the corresponding positive solution $u\left(
t\right) =S_{\varepsilon }\left( t\right) u_{0}$ to the non-degenerate
aggregation equations (\ref{bed1})-(\ref{bed3}) converges to a single
stationary state $u_{\ast }$ of (\ref{stat}) in the sense that%
\begin{equation}
\lim_{t\rightarrow \infty }\left\Vert u\left( t\right) -u_{\ast }\right\Vert
_{L^{p}\left( \Omega \right) }=0,  \label{convlinftime}
\end{equation}%
for any $p>1$. Moreover, there exist constants $C>0$, $\rho =\rho \left(
p,u_{\ast }\right) \in \left( 0,1\right) $ such that the convergence rate
holds:%
\begin{equation}
\left\Vert u\left( t\right) -u_{\ast }\right\Vert _{L^{p}\left( \Omega
\right) }\leq C\left( 1+t\right) ^{-\rho },  \label{conv_rate}
\end{equation}%
for all $t\geq 0.$
\end{theorem}

Owing to (\ref{oml}), Proposition \ref{omegalimit}-(iii), Proposition \ref%
{omega_bb}, and the regularity properties of $u\in C^{\alpha /2}\left(
(0;\infty );C^{\alpha }\left( \overline{\Omega }\right) \right) $, we may
then assume without loss of generality that%
\begin{equation}
\inf_{x\in \Omega }u\left( t,x\right) \geq \underline{u}>0,\text{ for all }%
t>0.  \label{infsol}
\end{equation}

We employ a generalized version of the {\L }ojasiewicz-Simon theorem proved
in \cite[Theorem 6]{GG2} (cf. also \cite[Lemma 2.20]{GM}). The version that
applies to our case is formulated in the following.

\begin{lemma}
\label{LS_lemma} Under the assumptions of Theorem \ref{conv_equil}, there
exist constants $\theta \in (0,\frac{1}{2}],$ $C>0,$ $\delta >0$ such that
the following inequality holds:%
\begin{equation}
\left\vert \mathcal{E}_{\varepsilon }\left( u\right) -\mathcal{E}%
_{\varepsilon }\left( u_{\ast }\right) \right\vert ^{1-\theta }\leq
C\left\Vert \mu -\left\langle \mu \right\rangle \right\Vert _{L^{2}\left(
\Omega \right) },  \label{LSineq}
\end{equation}%
for all $u\in L^{\infty }\left( \Omega \right) \cap H^{1}\left( \Omega
\right) $ provided that $\left\Vert u-u_{\ast }\right\Vert _{L^{2}\left(
\Omega \right) }\leq \delta .$ Here $\mu =\mu \left( u\right) $ denotes $%
\Phi _{\varepsilon }^{^{\prime }}\left( u\right) -\mathcal{K}\ast u.$
\end{lemma}

\begin{proof}
We will now apply the abstract result \cite[Theorem 6]{GG2} to the energy
functional $\mathcal{E}_{\varepsilon }\left( u\right) $, which according to (%
\ref{ener_nond}) is the sum of entropy and an interface energy term. In
contrast to this feature, we shall split $\mathcal{E}_{\varepsilon }\left(
\varphi \right) $ into the sum of a convex (entropy) functional $\Sigma
:L^{2}\left( \Omega \right) \rightarrow \mathbb{R}\cup \left\{ \infty
\right\} $, with a suitable effective domain, and a non-local interaction
functional $\Psi :L^{2}\left( \Omega \right) \rightarrow \mathbb{R}$. Let $%
\underline{u}>0$ be the lower bound from (\ref{infsol}) and recall that $%
\varepsilon >0$. We define the lower-semicontinuous and strongly convex
functional $\Sigma =\Sigma _{\varepsilon }$ by%
\begin{equation*}
\Sigma \left( u\right) :=\left\{ 
\begin{array}{ll}
\int_{\Omega }\left( \Phi _{\varepsilon }\left( u\right) -\varepsilon \ln
\left( \underline{u}\right) u\right) dx, & \text{if }u\in L^{\infty }\left(
\Omega \right) ,\text{ }u\geq 0, \\ 
+\infty , & \text{otherwise,}%
\end{array}%
\right.
\end{equation*}%
with closed effective domain dom$\left( \Sigma \right) =\mathcal{Z}_{DN}\cap
H^{1}\left( \Omega \right) $, and the quadratic functional $\Psi =\Psi
_{\varepsilon }:L^{2}\left( \Omega \right) \rightarrow \mathbb{R}$, given by%
\begin{equation*}
\Psi \left( u\right) :=-\frac{1}{2}\int_{\Omega }\left[ u\left( \mathcal{K}%
\ast u\right) -2\varepsilon \ln \left( \underline{u}\right) u\right] dx.
\end{equation*}%
We have that $\Sigma $ is Fr\'{e}chet differentiable on any open subset $%
\overline{U}$ of%
\begin{equation*}
U_{M}:=\left\{ \psi \in L^{\infty }(\Omega ):\left\vert \left\langle \psi
\right\rangle \right\vert \leq M;\text{ }\underline{u}\leq \psi \left(
x\right) \leq C_{M}\right\} ,
\end{equation*}%
(where $C_{M}>0$ is such that $\left\Vert \psi \right\Vert _{L^{\infty
}\left( \Omega \right) }\leq C_{M},$ since $\psi $ is bounded) with Fr\'{e}%
chet derivative $D\Sigma :\overline{U}\rightarrow L^{\infty }\left( \Omega
\right) $ having the form%
\begin{equation*}
\left\langle D\Sigma \left( u\right) ,\xi \right\rangle =\int_{\Omega
}\left( \Phi _{\varepsilon }^{^{\prime }}\left( u\right) -\varepsilon \ln
\left( \underline{u}\right) \right) \cdot \xi dx,
\end{equation*}%
for all $u\in \overline{U}$ and $\xi \in L^{\infty }\left( \Omega \right) $.
The analyticity of $D\Sigma $ as a mapping on $L^{\infty }\left( \Omega
\right) $ is standard owing to the analyticity of $\Phi $ (see, e.g., \cite[%
Remark 3]{GG2}). Moreover, due to assumptions on $A$ in Theorem \ref%
{well-nondeg}, we have $\Phi ^{^{\prime }}\left( y\right) \sim A^{^{\prime
}}\left( y\right) \sim y^{m-1}$ and recalling that $\Phi _{\varepsilon
}^{^{\prime }}\left( y\right) =\Phi ^{^{\prime }}\left( y\right)
+\varepsilon \ln \left( y\right) $, one has $\Phi ^{^{\prime }}\left( 
\underline{u}\right) >0$ and%
\begin{equation*}
\left\langle D\Sigma \left( u_{1}\right) -D\Sigma \left( u_{2}\right)
,u_{1}-u_{2}\right\rangle \geq \Phi ^{^{\prime }}\left( \underline{u}\right)
\left\Vert u_{1}-u_{2}\right\Vert _{L^{2}\left( \Omega \right) }^{2},
\end{equation*}%
for all $u_{1},u_{2}\in \overline{U}$, and%
\begin{equation*}
\left\Vert D\Sigma \left( u_{1}\right) -D\Sigma \left( u_{2}\right)
\right\Vert _{\left( L^{2}\left( \Omega \right) \right) ^{\ast }}\leq \gamma
\left\Vert u_{1}-u_{2}\right\Vert _{L^{2}\left( \Omega \right) },
\end{equation*}%
for some positive constant $\gamma =\gamma \left( C_{M},\underline{u}%
,\varepsilon \right) .$ Moreover, computing the second Fr\'{e}chet
derivative $D^{2}\Sigma $\ of $\Sigma ,$%
\begin{equation*}
\left\langle D^{2}\Sigma \left( u\right) \xi _{1},\xi _{2}\right\rangle
=\int_{\Omega }\Phi _{\varepsilon }^{\prime \prime }\left( u\right) \xi
_{1}\cdot \xi _{2}dx
\end{equation*}%
yields that $D^{2}\Sigma \in \mathcal{L}\left( L^{\infty }\left( \Omega
\right) ,L^{\infty }\left( \Omega \right) \right) $ is an isomorphism for
every $u\in \overline{U}$, owing to the fact that $\Phi _{\varepsilon
}^{\prime \prime }\left( u\right) =A_{\varepsilon }^{^{\prime }}\left(
u\right) /u\geq \zeta =\zeta \left( C_{M},\underline{u}\right) >0$.
Concerning the (quadratic) function $\Psi $, we see that%
\begin{equation*}
\Psi \left( u\right) =\frac{1}{2}\left\langle -\mathcal{K}\ast
u,u\right\rangle _{L^{2}\left( \Omega \right) }+\left\langle \varepsilon \ln
\left( \underline{u}\right) ,u\right\rangle _{L^{2}\left( \Omega \right) }%
\text{, }\forall u\in L^{2}\left( \Omega \right) .
\end{equation*}%
We recall that the linear operator $\psi \mapsto \mathcal{K}\ast \psi $ is
self-adjoint and compact from $L^{2}\left( \Omega \right) $ to itself and is
also compact from $L^{\infty }(\Omega )$ to $C^{0}(\overline{\Omega })$
(since $\mathcal{K}\in W_{\text{loc}}^{1,1}$). On the other hand, we also
have the following (orthogonal) sum decomposition of $L^{2}\left( \Omega
\right) =L_{0}^{2}\left( \Omega \right) \oplus H_{1}$, where%
\begin{equation*}
L_{0}^{2}\left( \Omega \right) :=\left\{ u\in L^{2}\left( \Omega \right)
:\left\langle u\right\rangle =0\right\} \text{, }H_{1}:=\left\{ u\in
L^{2}\left( \Omega \right) :u=\text{const.}\right\}
\end{equation*}%
Then, the annihilator of $L_{0}^{2}\left( \Omega \right) $ is the
one-dimensional subspace 
\begin{equation*}
L_{0}^{0}:=\left\{ ch\in \left( L^{2}\left( \Omega \right) \right) ^{\ast
}:c\in \mathbb{R}\right\}
\end{equation*}%
of constant functions, where $h\in \left( L^{2}\left( \Omega \right) \right)
^{\ast }\simeq L^{2}\left( \Omega \right) $ is given by $\left\langle
h,u\right\rangle =$ $\frac{1}{\left\vert \Omega \right\vert }\int_{\Omega
}udx$, $u\in L^{2}\left( \Omega \right) .$ Hence, the hypotheses of \cite[%
Theorem 6]{GG2} are satisfied and the sum%
\begin{equation*}
\mathcal{E}_{\varepsilon }=\Sigma +\Psi :L^{2}\left( \Omega \right)
\rightarrow \mathbb{R}\cup \left\{ \infty \right\}
\end{equation*}%
is a well defined, bounded from below functional with nonempty, closed, and
convex effective domain dom$\left( \mathcal{E}_{\varepsilon }\right) =$dom$%
\left( \Sigma \right) .$ Unravelling notation in \cite[Theorem 6]{GG2}, and
observing that the Fr\'{e}chet derivative 
\begin{equation*}
D\mathcal{E}_{\varepsilon }\left( u\right) =\Phi _{\varepsilon }^{^{\prime
}}\left( u\right) -\mathcal{K}\ast u=:\mu ,
\end{equation*}%
we have%
\begin{align*}
\left\vert \mathcal{E}_{\varepsilon }\left( u\right) -\mathcal{E}%
_{\varepsilon }\left( u_{\ast }\right) \right\vert ^{1-\theta }& \leq
C\inf_{u\in L^{2}\left( \Omega \right) }\left\{ \left\Vert D\mathcal{E}%
_{\varepsilon }\left( u\right) -\mu _{\ast }\right\Vert _{L^{2}\left( \Omega
\right) }:\mu _{\ast }\in L_{0}^{0}\right\} \\
& =C\left\Vert \mu -\left\langle \mu \right\rangle \right\Vert _{L^{2}\left(
\Omega \right) },
\end{align*}%
from which (\ref{LSineq}) follows.
\end{proof}

We are now in a position to complete the proof of Theorem \ref{conv_equil}.

\begin{proof}[Proof of Theorem \protect\ref{conv_equil}]
We argue as in the proof of \cite[Theorem 2.21]{GM}. First, we note that by
virtue of the regularity results proven in the previous section (see, e.g., (%
\ref{reg_sec})), all $u_{\ast }\in \omega \left[ u_{0}\right] $ are bounded
in $C^{\alpha }\left( \overline{\Omega }\right) \cap H^{1}\left( \Omega
\right) $. Besides, recalling Proposition \ref{omegalimit}-(ii), we have 
\begin{equation*}
\mathcal{E}_{\varepsilon }\left( u\left( t\right) \right) \rightarrow 
\mathcal{E}_{\varepsilon ,\infty },\text{ as }t\rightarrow \infty ,
\end{equation*}%
and the limit energy $\mathcal{E}_{\varepsilon ,\infty }$ is the same for
every steady-state solution $u_{\ast }\in $ $\omega \left[ u_{0}\right] $.
Moreover, we can integrate (\ref{nonincr}) over $\left( t,\infty \right) $
to get 
\begin{align}
& \underline{u}\int_{t}^{\infty }\int_{\Omega }\left\vert \nabla \Phi
_{\varepsilon }^{^{\prime }}\left( u\left( t\right) \right) -\nabla \mathcal{%
K}\ast u\left( t\right) \right\vert ^{2}dxds  \label{convls0} \\
& \leq \int_{t}^{\infty }\int_{\Omega }u\left( t\right) \left\vert \nabla
\Phi _{\varepsilon }^{^{\prime }}\left( u\left( t\right) \right) -\nabla 
\mathcal{K}\ast u\left( t\right) \right\vert ^{2}dxds  \notag \\
& =\mathcal{E}_{\varepsilon }\left( u\left( t\right) \right) -\mathcal{E}%
_{\varepsilon ,\infty }=\mathcal{E}_{\varepsilon }\left( u\left( t\right)
\right) -\mathcal{E}_{\varepsilon }\left( u_{\ast }\right) .  \notag
\end{align}%
By virtue of Lemma \ref{LS_lemma} and recalling that $\mu \left( t\right)
=\Phi _{\varepsilon }^{^{\prime }}\left( u\left( t\right) \right) -\nabla 
\mathcal{K}\ast u\left( t\right) $, we have%
\begin{equation}
\left\vert \mathcal{E}_{\varepsilon }\left( u\left( t\right) \right) -%
\mathcal{E}_{\varepsilon }\left( u_{\ast }\right) \right\vert ^{1-\theta
}\leq C\left\Vert \mu \left( t\right) -\left\langle \mu \left( t\right)
\right\rangle \right\Vert _{L^{2}\left( \Omega \right) }\leq C\left\Vert
\nabla \mu \left( t\right) \right\Vert _{L^{2}\left( \Omega \right) }
\label{convls00}
\end{equation}%
exploiting Poincare's inequality, provided that%
\begin{equation}
\left\Vert u-u_{\ast }\right\Vert _{L^{2}\left( \Omega \right) }\leq \delta .
\label{small}
\end{equation}%
This, combined with the previous identity, yields%
\begin{equation}
\int_{t}^{\infty }\left\Vert \nabla \mu \left( s\right) \right\Vert
_{L^{2}\left( \Omega \right) }^{2}ds\leq C\left\Vert \nabla \mu \left(
t\right) \right\Vert _{L^{2}\left( \Omega \right) }^{\frac{1}{\left(
1-\theta \right) }},  \label{convls1}
\end{equation}%
for all $t>0$, for as long as (\ref{small}) holds. Note that, in general,
the quantities $\theta ,$ $C$ and $\delta $ above may depend on $u_{\ast }>0$
and $\varepsilon >0$. Finally, let us set%
\begin{equation*}
\hat{W}=\cup \left\{ \mathcal{I}:\mathcal{I}\text{ is an open interval on
which (\ref{small}) holds}\right\} .
\end{equation*}%
Clearly, $\hat{W}$ is nonempty since $u_{\ast }\in $ $\omega \left[ u_{0}%
\right] $. We can now use (\ref{convls1}), the fact that $\left\Vert \nabla
\mu \left( t\right) \right\Vert _{L^{2}\left( \Omega \right) }\in
L^{2}\left( 0,\infty \right) $, and exploit \cite[Lemma 5.1]{FLP} (with $%
\alpha =2\left( 1-\theta \right) $) to deduce that $\left\Vert \nabla \mu
\left( \cdot \right) \right\Vert _{L^{2}\left( \Omega \right) }\in L^{1}(%
\hat{W})$ and%
\begin{equation}
\int_{\hat{W}}\left\Vert \nabla \mu \left( s\right) \right\Vert
_{L^{2}\left( \Omega \right) }ds\leq C\left( u_{\ast },\underline{u}\right)
<\infty .  \label{convls2bis}
\end{equation}%
Consequently, using the bound (\ref{convls2bis}) and the main equation (\ref%
{bed1}), which also reads $\partial _{t}u\left( t\right) =div\left( u\left(
t\right) \nabla \mu \left( t\right) \right) $, we obtain%
\begin{equation}
\int_{\hat{W}}\left\Vert \partial _{t}u\left( s\right) \right\Vert _{\left(
H^{1}\left( \Omega \right) \right) ^{\ast }}ds\leq C<\infty .
\label{convls2tris}
\end{equation}%
In order to finish the proof of the convergence result in (\ref{convlinftime}%
) it suffices to show that it holds in $L^{2}$-norm. Indeed, in this case (%
\ref{convlinftime}) will become an immediate consequence of the $L^{2}$-$%
(C^{\alpha }\cap H^{1})$ smoothing property of the bounded solutions $%
u\left( t\right) $ and all $u_{\ast }\in \omega \left[ u_{0}\right] $. We
claim that we can find a sufficiently large time $\tau >0$ such that $\left(
\tau ,\infty \right) \subset \hat{W}$. To this end, recalling (\ref{convls0}%
) and the above bounds, we also have that $\partial _{t}\varphi \in
L^{2}(0,\infty ;\left( H^{1}\left( \Omega \right) \right) ^{\ast })$, $%
\nabla \mu \in L^{2}(0,\infty ;\left( L^{2}\left( \Omega \right) \right)
^{d})$ and, furthermore, for any $k>0$ there exists a time $t_{\ast
}=t_{\ast }\left( k\right) >0$ such that 
\begin{equation}
\left\Vert \partial _{t}u\right\Vert _{L^{1}(\hat{W}\cap \left( t_{\ast
},\infty \right) ;\left( H^{1}\right) ^{\ast })}\leq k,\text{ }\left\Vert
\partial _{t}u\right\Vert _{L^{2}(\left( t_{\ast },\infty \right) ;\left(
H^{1}\right) ^{\ast })}\leq k,\text{ }\left\Vert \nabla \mu \right\Vert
_{L^{2}(\left( t_{\ast },\infty \right) ;\left( L^{2}\right) ^{d})}\leq k.
\label{convls3}
\end{equation}%
Next, observe that by the regularity properties of $u$ (see Section \ref%
{optim}), there is a time $t_{\#}>0$ such that%
\begin{equation}
\sup_{t\geq t_{\#}}\left\Vert u\left( t\right) \right\Vert _{H^{1}\cap
C^{\alpha }\left( \overline{\Omega }\right) }\leq C.  \label{bls}
\end{equation}%
Now, let $\left( t_{0},t_{2}\right) \subset \hat{W}$, for some $%
t_{2}>t_{0}\geq t_{\ast }\left( k\right) ,$ $\left\vert
t_{0}-t_{2}\right\vert \geq 1$ such that (\ref{bls}) holds (without loss of
generality, we can assume that $t_{\ast }\geq t_{\#}$). This claim is an
immediate consequence of the aforementioned $L^{2}$-$(H^{1}\cap C^{\alpha
}\left( \overline{\Omega }\right) )$ smoothing property and bounds (\ref%
{convls3}). Using (\ref{convls3}) and (\ref{bls}), we obtain%
\begin{align}
\left\Vert u\left( t_{0}\right) -u\left( t_{2}\right) \right\Vert
_{L^{2}\left( \Omega \right) }^{2}& =2\int_{t_{0}}^{t_{2}}\left\langle
\partial _{t}u\left( s\right) ,u\left( s\right) -u\left( t_{0}\right)
\right\rangle _{L^{2}\left( \Omega \right) }ds  \label{bls2b} \\
& \leq 2\int_{t_{0}}^{t_{2}}\left\Vert \partial _{t}u\left( s\right)
\right\Vert _{\left( H^{1}\left( \Omega \right) \right) ^{\ast }}\left(
\left\Vert u\left( s\right) \right\Vert _{H^{1}\left( \Omega \right)
}+\left\Vert u\left( t_{0}\right) \right\Vert _{H^{1}\left( \Omega \right)
}\right) ds  \notag \\
& \leq C\left\Vert \partial _{t}u\right\Vert _{L^{1}(t_{0},t_{2};\left(
H^{1}\right) ^{\ast })}\left( \left\Vert u\right\Vert _{L^{\infty }(t_{\ast
},\infty ;H^{1})}+1\right) \leq Ck.  \notag
\end{align}%
Therefore we can choose a time $t_{\ast }\left( k\right) =\tau <t_{0}<t_{2}$%
, such that%
\begin{equation}
\left\Vert u\left( t_{0}\right) -u\left( t_{2}\right) \right\Vert
_{L^{2}\left( \Omega \right) }<\frac{\delta }{3}  \label{bls3}
\end{equation}%
provided that (\ref{small}) holds for all $t\in \left( t_{0},t_{2}\right) $.
Since $u_{\ast }\in \omega \left[ u_{0}\right] $, a large (redefined) $\tau $
can be chosen such that 
\begin{equation}
\left\Vert u\left( \tau \right) -u_{\ast }\right\Vert _{L^{2}\left( \Omega
\right) }<\frac{\delta }{3},  \label{bls4}
\end{equation}%
whence, (\ref{bls3}) yields $\left( \tau ,\infty \right) \subset \hat{W}$.
Indeed, taking%
\begin{equation*}
\overline{t}=\inf \left\{ t>\tau :\left\Vert u\left( t\right) -u_{\ast
}\right\Vert _{L^{2}\left( \Omega \right) }\geq \delta \right\} ,
\end{equation*}%
we have $\overline{t}>\tau $ and $\left\Vert u\left( \overline{t}\right)
-u_{\ast }\right\Vert _{L^{2}\left( \Omega \right) }\geq \varepsilon $ if $%
\overline{t}$ is finite. On the other hand, in view of (\ref{bls3}) and (\ref%
{bls4}), we have%
\begin{equation*}
\left\Vert u\left( t\right) -u_{\ast }\right\Vert _{L^{2}\left( \Omega
\right) }\leq \left\Vert u\left( t\right) -u\left( \tau \right) \right\Vert
_{L^{2}\left( \Omega \right) }+\left\Vert u\left( \tau \right) -u_{\ast
}\right\Vert _{L^{2}\left( \Omega \right) }<\frac{2\delta }{3},
\end{equation*}%
for all $\overline{t}>t\geq \tau $, and this leads to a contradiction.
Therefore, $\overline{t}=\infty $ and by (\ref{convls3}) the integrability
of $\partial _{t}u$ in $L^{1}(\tau ,\infty ;\left( H^{1}\left( \Omega
\right) \right) ^{\ast })$ follows. Hence, $\omega \left[ u_{0}\right]
=\{u_{\ast }\}$ and \eqref{convlinftime} holds on account of the $L^{2}$-$%
(H^{1}\cap C^{\alpha }\left( \overline{\Omega }\right) )$ smoothing
property. The proof is finished.
\end{proof}

\begin{proof}[Proof of (\protect\ref{conv_rate})]
Without loss of generality, suppose now that, for all $t\geq t_{\ast }>0,$
we have $\mathcal{E}_{\varepsilon }\left( u\left( t\right) \right) >\mathcal{%
E}_{\varepsilon }\left( u_{\ast }\right) $ (otherwise, there is nothing to
prove). Define the function%
\begin{equation*}
\Xi \left( t\right) :=\mathcal{E}_{\varepsilon }\left( u\left( t\right)
\right) -\mathcal{E}_{\varepsilon }\left( u_{\ast }\right)
\end{equation*}%
and observe that by (\ref{convls0})-(\ref{convls00}) and (\ref{nonincr}), it
satisfies%
\begin{equation*}
\frac{d}{dt}\Xi \left( t\right) +C\Xi \left( t\right) ^{2\left( 1-\theta
\right) }\leq 0,\text{ for all }t\geq t_{\ast }
\end{equation*}%
for some positive constant $C=C\left( \underline{u}\right) .$ Integration of
the preceding inequality yields%
\begin{equation}
\Xi \left( t\right) \leq \Xi \left( 0\right) \left( 1+C\Xi \left( 0\right)
^{1-2\theta }t\right) ^{-\frac{1}{1-2\theta }},  \label{bls5}
\end{equation}%
for all $t\geq t_{\ast }.$ On the other hand, we have on account (\ref%
{convls0}) that%
\begin{equation*}
-\frac{d}{dt}\Xi \left( t\right) ^{\theta }=-\theta \Xi \left( t\right)
^{\theta -1}\frac{d}{dt}\Xi \left( t\right) \geq C\theta \left\Vert \nabla
\mu \left( t\right) \right\Vert _{L^{2}\left( \Omega \right) },
\end{equation*}%
for all $t\geq t_{\ast }$, provided that $\left\Vert u\left( t\right)
-u_{\ast }\right\Vert \leq \delta .$ Integrating this inequality over $%
\left( t,\infty \right) ,$ we also get%
\begin{equation*}
\int_{t}^{\infty }\left\Vert \nabla \mu \left( s\right) \right\Vert
_{L^{2}\left( \Omega \right) }ds\leq C<\infty ,
\end{equation*}%
for all $t\geq t_{\ast }.$ As above, we obtain%
\begin{equation*}
\int_{t}^{\infty }\left\Vert \partial _{t}u\left( s\right) \right\Vert
_{\left( H^{1}\left( \Omega \right) \right) ^{\ast }}ds\leq C<\infty ,
\end{equation*}%
and combining with (\ref{convls0}) and (\ref{bls5}) yields%
\begin{align*}
\left\Vert u\left( t\right) -u_{\ast }\right\Vert _{\left( H^{1}\left(
\Omega \right) \right) ^{\ast }}& \leq \int_{t}^{\infty }\left\Vert \partial
_{t}u\left( s\right) \right\Vert _{\left( H^{1}\left( \Omega \right) \right)
^{\ast }}ds\leq C\int_{t}^{\infty }\left\Vert \nabla \mu \left( s\right)
\right\Vert _{L^{2}\left( \Omega \right) }ds \\
& \leq C\left( 1+t\right) ^{-\frac{\theta }{1-2\theta }},
\end{align*}%
for some positive constant $C$, which depends on $\Xi \left( 0\right) ,$ $%
\theta ,$ $\underline{u}$ and $\varepsilon .$ By using standard
interpolation inequalities (see, Section \ref{non_deg}) one can deduce the
convergence rate estimate in the stronger norm in (\ref{conv_rate}). Of
course, the convergence exponent deteriorates. The proof is complete.
\end{proof}

\begin{acknowledgement}
The author thanks the anonymous referees for their careful reading of the
manuscript and many useful comments.
\end{acknowledgement}

\end{document}